\documentclass[12pt,article]{article}
\usepackage[cp1250]{inputenc}
\usepackage{lmodern}
\usepackage[T1]{fontenc}
\usepackage{ifthen}
\usepackage{pgfpages}
\usepackage{amssymb}
\usepackage{epstopdf}
\usepackage{pgf,tikz}
\usetikzlibrary{arrows}

\def\CC{\mathbb C}
\def\RR{\mathbb R}
\def\QQ{\mathbb Q}
\def\ZZ{\mathbb Z}
\def\NN{\mathbb N}

\newtheorem{prop}{Proposition}

\newtheorem{defin}{Definition}
\newtheorem{remark}{Remark}
\newtheorem{example}{Example}

\newtheorem{prob}{Problem}
\newtheorem{algo}{Algorithm}

\date{}

\renewcommand{\baselinestretch}{1.1}\normalsize
\newcommand{\qed}{\hfill ~$\square$\bigskip}

\begin{document}

\begin{center}
\begin{Large}
\textbf{Platonic configurations of points and lines} \vspace{10mm}
\end{Large}

{\large {Jurij Kovi\v c\footnote{\textit{This work is supported in part by the Slovenian Research Agency (research program P1-0294 and research projects N1-0032, J1-9187)\\
\indent E-mail addresses}: jurij.kovic@siol.net (Jurij Kovi\v c), aleks.simonic@gmail.com (Aleksander Simoni\v c).}}} \vspace{2mm}

\begin{small}
\textit{IMFM, Ljubljana and UP FAMNIT, Koper, Slovenia} \vspace{5mm}
\end{small}

{\large {Aleksander Simoni\v c}} \vspace{2mm}

\begin{small}
\textit{School of Science. The University of New South Wales (Canberra), ACT, Australia} \vspace{5mm}
\end{small}

\end{center}

\parindent0,5cm

\noindent {\sc abstract:} We present some methods for constructing connected spatial geometric configurations $(p_{q}, n_{k})$ of points and lines, preserved by the same rotations (and reflections) of Euclidean space $E^{3}$ as the chosen Platonic solid. In this paper we are primarily interested in balanced configurations $(n_{3}), (n_{4})$ and $(n_{5})$, but also in unbalanced configurations $(p_{3},n_{4}), (p_{3}, n_{5})$ and $(p_{4}, n_{5})$.  \\

\noindent {\sc keywords:} Configuration of points and lines, symmetry group, Platonic solid, centrally symmetric solid, concentric solid,  projection from a point.\\


\section{Basic notions, formulation of the problem}

Platonic configurations are "configurations with the symmetries of a given Platonic solid". Before giving a precise definition of this concept we briefly review the concepts of configurations and their symmetries. Basic notions, recent results and wider context of configurations are explained in the books \cite{Grunbaum} and \cite{PisaServ}. In this section we explain the terminology and notation we use in this paper. 

An \emph{incidence structure} consists of "points", "blocks", and "incidences" between them (given by a 0-1 matrix or by a bipartite graph). A \emph{configuration} satisfies an additional condition: no two blocks are both incident to more than one common point and every point is incident to at least two blocks. We distinguish between abstract (or combinatorial) configurations and geometric configurations. A \textit{combinatorial configuration} $X$ is completely determined by the incidences of its \textit{points} and \textit{lines}. A \textit{geometric configuration} $Y = Y(X)$ (called also a \emph{geometric realisation} of $X$) consists of  points and (straight) lines in the Euclidean space $E^{3}$. Points incident with $q$ lines are called $q$-\textit{valent}, and likewise $k$-\textit{valent} lines are lines incident with $k$ points.

For denoting configurations we use the same notation as Gr\"unbaum in \cite{Grunbaum}. Let $(p_{q}, n_{k})$ be a configuration of $p$ $q$-valent points and $n$ $k$-valent lines; it is called also a $(q,k)$-configuration; if $p = n$ and $q = k$, it is called a \textit{balanced} configuration and denoted simply as $(n_{k})$; it is called also a $k$-configuration. Thus, a configuration $(n_{k})$ consists of $n$ \textit{points} and $n$ \textit{lines}, each point being incident with $k$ lines and each line incident with $k$ points.  %Note that in general holds: $pq = nk$. 
For a configuration containing $p_{1}, p_{2}, p_{3} \dots$ points of valence $q_{1}, q_{2},$ $ q_{3}\dots$, respectively, and $n_{1}, n_{2}, n_{3} \dots$ lines of valence $k_{1}, k_{2}, k_{3}\dots$, respectively, we use the natural generalization of this notation, as in the configuration $(6_{4}12_{3},15_{4})$ shown in Figure 1.

\begin{center}
\definecolor{xdxdff}{rgb}{0.49019607843137253,0.49019607843137253,1.0}
\definecolor{uuuuuu}{rgb}{0.26666666666666666,0.26666666666666666,0.26666666666666666}
\definecolor{zzttqq}{rgb}{0.6,0.2,0.0}
\definecolor{qqqqff}{rgb}{0.0,0.0,1.0}
\begin{tikzpicture}[line cap=round,line join=round,>=triangle 45,x=0.3cm,y=0.3cm]
\clip(8.321719999999997,0.8406600000000017) rectangle (16.041519999999995,8.666939999999995);
\fill[color=zzttqq,fill=zzttqq,fill opacity=0.1] (-0.2,2.64) -- (2.52,2.68) -- (3.8453589838486226,5.055589098293674) -- (2.450717967697245,7.391178196587347) -- (-0.2692820323027547,7.351178196587348) -- (-1.5946410161513782,4.975589098293676) -- cycle;
\fill[color=zzttqq,fill=zzttqq,fill opacity=0.1] (0.637606315720525,4.141335454421212) -- (1.6443316611283283,4.146471588246776) -- (2.1432463114620566,5.0208893789163795) -- (1.63543561638798,5.89017103576042) -- (0.6287102709801761,5.885034901934857) -- (0.12979562064644834,5.010617111265253) -- cycle;
\fill[color=zzttqq,fill=zzttqq,fill opacity=0.1] (10.98,2.86) -- (13.1,2.84) -- (14.177320508075685,4.665973856023011) -- (13.134641016151374,6.511947712046019) -- (11.014641016151376,6.531947712046018) -- (9.937320508075688,4.705973856023009) -- cycle;
\draw [color=zzttqq] (-0.2,2.64)-- (2.52,2.68);
\draw [color=zzttqq] (2.52,2.68)-- (3.8453589838486226,5.055589098293674);
\draw [color=zzttqq] (3.8453589838486226,5.055589098293674)-- (2.450717967697245,7.391178196587347);
\draw [color=zzttqq] (2.450717967697245,7.391178196587347)-- (-0.2692820323027547,7.351178196587348);
\draw [color=zzttqq] (-0.2692820323027547,7.351178196587348)-- (-1.5946410161513782,4.975589098293676);
\draw [color=zzttqq] (-1.5946410161513782,4.975589098293676)-- (-0.2,2.64);
\draw (-0.2,2.64)-- (2.450717967697245,7.391178196587347);
\draw (2.52,2.68)-- (-0.2692820323027547,7.351178196587348);
\draw [color=zzttqq] (0.637606315720525,4.141335454421212)-- (1.6443316611283283,4.146471588246776);
\draw [color=zzttqq] (1.6443316611283283,4.146471588246776)-- (2.1432463114620566,5.0208893789163795);
\draw [color=zzttqq] (2.1432463114620566,5.0208893789163795)-- (1.63543561638798,5.89017103576042);
\draw [color=zzttqq] (1.63543561638798,5.89017103576042)-- (0.6287102709801761,5.885034901934857);
\draw [color=zzttqq] (0.6287102709801761,5.885034901934857)-- (0.12979562064644834,5.010617111265253);
\draw [color=zzttqq] (0.12979562064644834,5.010617111265253)-- (0.637606315720525,4.141335454421212);
\draw (-1.5946410161513782,4.975589098293676)-- (3.8453589838486226,5.055589098293674);
\draw (0.6287102709801761,5.885034901934857)-- (-0.52,5.88);
\draw (0.12979562064644834,5.010617111265253)-- (-0.52,5.88);
\draw (0.637606315720525,4.141335454421212)-- (1.12,3.46);
\draw (1.6443316611283283,4.146471588246776)-- (1.12,3.46);
\draw (1.63543561638798,5.89017103576042)-- (2.64,5.92);
\draw (2.1432463114620566,5.0208893789163795)-- (2.64,5.92);
\draw (0.6287102709801761,5.885034901934857)-- (1.16,6.84);
\draw (1.63543561638798,5.89017103576042)-- (1.16,6.84);
\draw (0.12979562064644834,5.010617111265253)-- (-0.3,4.1);
\draw (0.637606315720525,4.141335454421212)-- (-0.3,4.1);
\draw (1.6443316611283283,4.146471588246776)-- (2.6,4.16);
\draw (2.1432463114620566,5.0208893789163795)-- (2.6,4.16);
\draw (-0.2692820323027547,7.351178196587348)-- (1.02,9.32);
\draw (1.02,9.32)-- (2.450717967697245,7.391178196587347);
\draw (-0.2692820323027547,7.351178196587348)-- (-2.76,7.28);
\draw (-2.76,7.28)-- (-1.5946410161513782,4.975589098293676);
\draw (2.450717967697245,7.391178196587347)-- (5.34,7.46);
\draw (3.8453589838486226,5.055589098293674)-- (5.34,7.46);
\draw (1.1,0.4)-- (2.52,2.68);
\draw (-0.2,2.64)-- (1.1,0.4);
\draw (-1.5946410161513782,4.975589098293676)-- (-2.98,2.64);
\draw (-2.98,2.64)-- (-0.2,2.64);
\draw (2.52,2.68)-- (5.08,2.68);
\draw (3.8453589838486226,5.055589098293674)-- (5.08,2.68);
\draw (-2.76,7.28)-- (-0.52,5.88);
\draw (2.6,4.16)-- (-0.52,5.88);
\draw (2.6,4.16)-- (5.08,2.68);
\draw (5.34,7.46)-- (2.64,5.92);
\draw (-0.3,4.1)-- (2.64,5.92);
\draw (-0.3,4.1)-- (-2.98,2.64);
\draw (1.02,9.32)-- (1.16,6.84);
\draw (1.16,6.84)-- (1.12,3.46);
\draw (1.12,3.46)-- (1.1,0.4);
\draw [color=zzttqq] (10.98,2.86)-- (13.1,2.84);
\draw [color=zzttqq] (13.1,2.84)-- (14.177320508075685,4.665973856023011);
\draw [color=zzttqq] (14.177320508075685,4.665973856023011)-- (13.134641016151374,6.511947712046019);
\draw [color=zzttqq] (13.134641016151374,6.511947712046019)-- (11.014641016151376,6.531947712046018);
\draw [color=zzttqq] (11.014641016151376,6.531947712046018)-- (9.937320508075688,4.705973856023009);
\draw [color=zzttqq] (9.937320508075688,4.705973856023009)-- (10.98,2.86);
\draw (15.219999999999997,2.819999999999998)-- (8.894641016151377,6.55194771204602);
\draw (8.86,2.879999999999998)-- (15.254641016151373,6.491947712046021);
\draw (12.091961524227067,8.357921568069024)-- (12.022679491924317,1.0140261439769869);
\draw (9.937320508075688,4.705973856023009)-- (13.134641016151374,6.511947712046019);
\draw (13.134641016151374,6.511947712046019)-- (13.1,2.84);
\draw (13.1,2.84)-- (9.937320508075688,4.705973856023009);
\draw (11.014641016151376,6.531947712046018)-- (10.98,2.86);
\draw (14.177320508075685,4.665973856023011)-- (10.98,2.86);
\draw (14.177320508075685,4.665973856023011)-- (11.014641016151376,6.531947712046018);
\draw (8.86,2.879999999999998)-- (12.091961524227067,8.357921568069024);
\draw (15.219999999999997,2.819999999999998)-- (12.091961524227067,8.357921568069024);
\draw (8.86,2.879999999999998)-- (15.219999999999997,2.819999999999998);
\draw (8.894641016151377,6.55194771204602)-- (12.022679491924317,1.0140261439769869);
\draw (12.022679491924317,1.0140261439769869)-- (15.254641016151373,6.491947712046021);
\draw (15.254641016151373,6.491947712046021)-- (8.894641016151377,6.55194771204602);
\draw (7.04396,0.46798)-- (-6.82506,-0.06442);
\begin{scriptsize}
\draw [fill=qqqqff] (-0.2,2.64) circle (1.5pt);
\draw [fill=qqqqff] (2.52,2.68) circle (1.5pt);
\draw [fill=uuuuuu] (3.8453589838486226,5.055589098293674) circle (1.5pt);
\draw [fill=uuuuuu] (2.450717967697245,7.391178196587347) circle (1.5pt);
\draw [fill=uuuuuu] (-0.2692820323027547,7.351178196587348) circle (1.5pt);
\draw [fill=uuuuuu] (-1.5946410161513782,4.975589098293676) circle (1.5pt);
\draw [fill=xdxdff] (0.637606315720525,4.141335454421212) circle (1.5pt);
\draw [fill=xdxdff] (1.6443316611283283,4.146471588246776) circle (1.5pt);
\draw [fill=uuuuuu] (2.1432463114620566,5.0208893789163795) circle (1.5pt);
\draw [fill=uuuuuu] (1.63543561638798,5.89017103576042) circle (1.5pt);
\draw [fill=uuuuuu] (0.6287102709801761,5.885034901934857) circle (1.5pt);
\draw [fill=uuuuuu] (0.12979562064644834,5.010617111265253) circle (1.5pt);
\draw [fill=qqqqff] (-0.52,5.88) circle (1.5pt);
\draw [fill=qqqqff] (1.12,3.46) circle (1.5pt);
\draw [fill=qqqqff] (2.64,5.92) circle (1.5pt);
\draw [fill=qqqqff] (1.16,6.84) circle (1.5pt);
\draw [fill=qqqqff] (-0.3,4.1) circle (1.5pt);
\draw [fill=qqqqff] (2.6,4.16) circle (1.5pt);
\draw [fill=qqqqff] (1.02,9.32) circle (1.5pt);
\draw [fill=qqqqff] (-2.76,7.28) circle (1.5pt);
\draw [fill=qqqqff] (5.34,7.46) circle (1.5pt);
\draw [fill=qqqqff] (1.1,0.4) circle (1.5pt);
\draw [fill=qqqqff] (-2.98,2.64) circle (1.5pt);
\draw [fill=qqqqff] (5.08,2.68) circle (1.5pt);
\draw [fill=qqqqff] (10.98,2.86) circle (1.5pt);
\draw [fill=qqqqff] (13.1,2.84) circle (1.5pt);
\draw [fill=uuuuuu] (14.177320508075685,4.665973856023011) circle (1.5pt);
\draw [fill=uuuuuu] (13.134641016151374,6.511947712046019) circle (1.5pt);
\draw [fill=uuuuuu] (11.014641016151376,6.531947712046018) circle (1.5pt);
\draw [fill=uuuuuu] (9.937320508075688,4.705973856023009) circle (1.5pt);
\draw [fill=uuuuuu] (12.091961524227067,8.357921568069024) circle (1.5pt);
\draw [fill=uuuuuu] (8.894641016151377,6.55194771204602) circle (1.5pt);
\draw [fill=uuuuuu] (8.86,2.879999999999998) circle (1.5pt);
\draw [fill=uuuuuu] (12.022679491924317,1.0140261439769869) circle (1.5pt);
\draw [fill=uuuuuu] (15.219999999999997,2.819999999999998) circle (1.5pt);
\draw [fill=uuuuuu] (15.254641016151373,6.491947712046021) circle (1.5pt);
\draw [fill=uuuuuu] (11.003094010767587,5.307965141364012) circle (1.5pt);
\draw [fill=uuuuuu] (12.068867513459482,5.909956426705015) circle (1.5pt);
\draw [fill=uuuuuu] (13.123094010767582,5.287965141364014) circle (1.5pt);
\draw [fill=uuuuuu] (13.100000000000001,2.84) circle (1.5pt);
\draw [fill=uuuuuu] (12.0457735026919,3.461991285341008) circle (1.5pt);
\draw [fill=uuuuuu] (10.991547005383792,4.083982570682006) circle (1.5pt);
\draw [fill=uuuuuu] (13.111547005383793,4.063982570682006) circle (1.5pt);
\draw [fill=qqqqff] (7.04396,0.46798) circle (1.5pt);
\draw [fill=qqqqff] (-6.82506,-0.06442) circle (1.5pt);
\draw [fill=qqqqff] (2.73152,0.97376) circle (1.5pt);
\end{scriptsize}
\end{tikzpicture}

\end{center}

\begin{center}
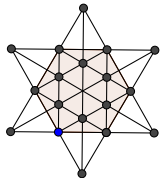
Figure 1. A configuration $(6_{4}12_{3},15_{4} )$ having 6 points of valence 4, 12 points of valence 3 and 15 lines of valence 4.
\end{center}

We are especially interested in symmetrical configurations. \textit{Symmetry} (or, more precisely $E^{n}$-symmetry) 
of any combinatorial object $X$ (e.g. solid, configuration, graph, etc.), geometrically realised as $Y = Y(X)$ in the Euclidean space $E^{n}$, is an isometry (rotation or reflection) of the Euclidean space $E^{n}$ preserving $Y$. The \textit{(full)} \textit{symmetry group} $Sym(Y)$ 
of $Y = Y(X) \subset E^{n}$ 
consists of all the rotations and reflections of the Euclidean space $E^{n}$ preserving $Y$; it is a subgroup of the group of automorphisms $Aut(X)$ of the corresponding $X$. The (rotational) symmetry group $Sym_{R}(Y)$ 
of $Y = Y(X) \subset E^{n}$ consists of all the rotations of the Euclidean space $E^{n}$ preserving $Y$.

Symmetries depend both on $X$ and on the geometry of the space (its dimension $n$ and its type -- Euclidean, projective, etc.) in which $X$ is embedded as a geometrical object $Y = Y(X)$. The number of symmetries may increase when the dimension $n$ of the space $E^{n}$ in which $X$ is embedded $n$ increases.

Various tools and techniques have been used to construct symmetrical planar geometric configurations $(n_{k})$ of points and lines (e. g. see \cite{BermBoko, BBGT}). It is well known that a \textit{planar} geometric configuration $Y = Y(X) \subset E^{2}$ may have only "cyclical" or "dihedral" symmetry.  However, the same underlying configuration $X$, which has a planar realisation $Y = Y(X)$, may have (and reveal) more (hidden) symmetries, if it is realised as a spatial configuration $Z = Z(X)$ in some higher-dimensional space.   %Thus, to obtain geometric configurations with richer groups of symmetry, we have to study spatial  configurations $Z = Z(X)$. 
Recently some authors started to investigate spatial configurations and symmetrical configurations more systematically (e.g. \cite{Gevayspatial, Gevaysymmetry}). 
For example, G\'abor G\'evay constructed a $(100_{4})$ configuration with the symmetry of a right pentagonal prism (\cite{PisaServ}, p. 249). 

The aim of this paper is to show how to construct \textit{spatial configurations} $Z \subset E^{3}$ of points and lines having the same symmetries %(i.e. being preserved by the same rotations and reflections of the Euclidean space $E^{3}$) 
as the chosen Platonic solid $P$. Such configurations we call \textit{Platonic configurations}. 

\begin{defin}
A \textit{Platonic configuration} is a geometric configuration $Z \subset E^{3}$ with symmetries of a Platonic solid $P \in \lbrace T,C,O,D,I\rbrace$. A "full" Platonic configuration $Z$ is preserved by all rotations and reflections of $E^{3}$ preserving $P$, thus $Sym(Z) = Sym(P)$. A "rotational" Platonic configurations $Z$ is preserved only by \textit{rotational symmetries} of $P$, thus: $Sym(Z) = Sym_{R}(P)$. For any $P \in \lbrace T, C, O, D, I \rbrace$ let $P(p_{q},n_{k})$ and $P_{R}(p_{q},n_{k})$, respectively, denote the classes of %"full" Platonic configurations and "rotational" Platonic configurations, respectively, belonging to the class of $(p_{q},n_{k})$ configurations.
"full" and "rotational" $(p_{q},n_{k})$ Platonic configurations, respectively.
\end{defin}

For example, the configurations obtained from the skeletons (vertices and edges) of the five Platonic solids $T, C, O, D, I$, respectively,  belong to the classes: $T (4_{3})$, $C (8_{3}, 12_{2})$, $O (6_{4},12_{2})$, $D(20_{3}, 30_{2})$, $I(12_{5},30_{2})$, respectively.

Now we can formulate the problem, on which we focus in this paper:

\begin{prob}\label{basicproblem}
%To construct examples of connected balanced geometric configurations $(n_{3})$, $(n_{4})$, $(n_{5})$, whose symmetry group $Sym(C)$ is isomorphic either to the full symmetry group $Sym(P)$ of rotations and reflections preserving the chosen (Platonic) solid $P$ or to the symmetry group $Sym_{R}(P)$, consisting only of rotations preserving $P$. 
For each of the five Platonic solids $P$ construct examples of connected Platonic configurations $Z \in P(p_{q},n_{k})$ and $Z \in P_{R}(p_{q}, n_{k})$ of points and lines, such that $3 \leq q \leq k \leq 5$. 
\end{prob}

In other words, we are looking for Platonic configurations of the following six valence-types: $(n_{3})$, $(n_{4})$, $(n_{5})$, $(p_{3},n_{4}),(p_{3},n_{5}), (p_{4},n_{5})$. Our problem splits into $5 \times 6 \times 2 = 60$ cases (there are 5 Platonic solids, 6 classes of configurations, and 2 types of Platonic configurations - "full" and "rotational").

In section \ref{twomethods} we present some methods 
for solving this problem. Using them, we may produce as many of Platonic (or even more general polyhedral) configurations as we want, and in many different ways as well. 

The first idea how to construct a Platonic configuration is 1) to place a copy of some smaller configuration $A$ "symetrically" around (or inside) a chosen Platonic solid $P$ (or around its "concentric copies") and 2) to use some additional lines and points (again placed "symmetrically" with respect to $P$) to link these building blocks in such a way that we obtain a \textit{connected}, \textit{symmetrical} and \textit{balanced} configuration $Z$ of the desired "target class" $(p_{q},n_{k})$. 

This idea is explained in more detail in Algorithm 1. Note that the "symmetrical places" (in which the identical copies of $A$ may be placed) may be the faces, edges, vertices, diagonal or axes of $P$ (and this list may still be not complete!), as it is shown in the examples in section 3.

In section 3 we give examples of Platonic configurations. Some of them are presented in the form of problems (with solutions), whose purpose is to help the reader to get the idea how to construct more examples of Platonic configurations. 

Most of the constructions may be used to produce infinite families of examples and may be easily generalised into theorems. But to make the constructions easier to understand we use the approach "from bottom up", with concrete examples and their visual representation.% (a more theoretical approach to Platonic configurations will be presented in the next paper). 

In the last section we present some open problems.

\section{How to construct Platonic configurations}\label{twomethods}

To find examples of each class $P(p_{q},n_{k})$ with $p$ as little as possible (the \textit{minimal} Platonic configurations of the given class)
we always try first to find example of a 1-\textit{layer Platonic configuration} $Z$ whose points all lie on the boundary of a chosen Platonic solid $P$. But sometimes the points of $Z$ lie on the faces of $s$ "concentric copies" of $P$, projected radially from the central point (centroid) of $P$; such $Z$ are called $s$-\textit{layer Platonic configurations}. 

In most of the examples in this paper the points of Platonic configurations are placed into vertices (represented as yellow points in the Figures), edges (red points) or faces (blue points) of a given Platonic solid (or of its concentric copies). The concepts of single-layer and multi-layer Platonic configurations are useful also for classification purposes. 

However, the subdivision of points of a Platonic configuration $Z$ into layers is not always the most appropriate description of their position in space. For example, some points of $Z$ may lie in a plane determined by two axes %or two diagonals 
of $P$. Such configurations we call \textit{spider-web configurations}.
Another interesting class %of Platonic configurations 
are \textit{helical configurations}, whose points and lines form (double or multiple) "helices" along each edge of a solid. 

To construct examples of connected Platonic configurations we combine some known methods with some new ones. Our approach is based on the very general idea (admitting many variations) of connecting isomorphic copies of some initial (planar) configuration $A$, placed "symmetrically" around the solid $P$.  The idea of obtaining bigger configurations by linking together small isomorphic building blocks has been explored by various authors, e.g. in \cite{BBGT}. Our approach of connecting smaller configurations using additional points and lines to obtain connected bigger configurations may be regarded also as the dual concept of the "splittable configurations" studied in \cite{Basicetal}. The idea of systematically constructing $(q',k')$-configurations from $(q,k)$-configurations, where $q' \leq q$ and $k' \leq k$, has been used e.g. in \cite{GevayPascal}.
Gr\"unbaum developed an \textit{incidence calculus} (see \cite{PisaServ}, pp. 243--263), composed of various operations on configurations (e.g. "parallel shift" connecting copies of the same configuration with parallel lines, or connecting copies of the same configuration, projected from some point, using lines through the projecting point, etc.). 
The idea to construct symmetric spatial configurations from polyhedra has been explored in \cite{Gevaysymmetry}. In $\cite{PisaServ}$, pp. 252--255,  the Grey configuration $(27_{3})$ is represented as a spatial geometric configurations with $3 \times 3 \times 3$ points in the form of a cube, and the generalized Grey configuration $(256_{4})$ is represented by four copies of cube with $4 \times 4 \times 4$ points, connected via parallel shift. Instead of parallel shift we will use "radial projection" and "antipodal lines" (as explained in Algorithm 1 below).

We will use the ``vector'' notation $P = (v,e,f,d,m)$ for the polyhedron $P$ with $v$ vertices of valence $d$, with $e$ edges and $f$ $m$-gonal faces. Thus we have the following parameters for the five Platonic solids:
% $(v_{d},e_{2})$:

\begin{center}
$\begin{array}{cccccc}
 & v & e & f & d & m\\
tetrahedron  & 4 & 6 & 4 & 3 & 3 \\ 
cube  & 8 & 12 & 6 & 3 & 4 \\ 
octahedron  & 6 & 12 & 8 & 4 & 3 \\ 
dodecahedron  & 20 & 30 & 12 & 3 & 5 \\ 
icosahedron  & 12 & 30 & 20 & 5 & 3 \\
\end{array} $
\end{center}

\begin{center}
Table 1: Parameters of Platonic solids.
\end{center}

To construct an example of a Platonic configuration $Z \in (p_{q},n_{k})$  with the symmetry group of the Platonic solid $P = (v,e,f,d,m)$ the following general procedure may be applied (at least for the cases $k = 3,4$ and 5 and $3 \leq q \leq k$):

\begin{algo}\label{construction}
Step 1. Start with a planar configuration $A$ with a cyclical $C_{m}$ or dihedral $D_{m}$ symmetry; all of its lines must have valence $k$ and each of its vertices must have a valence at most $k$.

Step 2. Place a copy of $A$ on each of its $f$ faces $P$ "symmetrically" -- in such a way that the symmetries of $P$ preserve also $B = fA$ (the set of copies of $A$).  

Step 3. Identifying some points of different copies of $A$ or connecting them "in a symmetrical way" try to get a connected configuration $B_{c}$ with the symmetries of $P$. 

Step 4. If necessary, adding new points and lines increase the valences of points to get a desired configuration $Z \in (p_{q},n_{k})$. 

\end{algo}

\begin{center}

\begin{remark}
In Step 2 you obtain a configuration $B = fA$ which may not be connected, but has the desired symmetry group $Sym(P)$ or $Sym_{R}(P)$ of $P$.

Note that if you (in Step 3) identify a 1-valent point of $A$ with a vertex of $P$, its valence increases to $d$. If you connect all the 1-valent points of each copy of $A$ with different points on the interiors of the edges of a given face, all these edge points have valence at least 2 (since each edge of $P$ is incident to two faces) and at most 3 (this happens if the edges of $P$ are included among the lines of the configuration $B_{c}$).

Step 4 may be divided into substeps (for example, we first increase the valence of 3-valent points to 4, then to 5, always taking care that the symmetry of $P$ is preserved). For the added lines we may usually use lines perpendicular to faces, or axes of $P$, or diagonals of $P$. Sometimes we use lines connecting pairs of antipodal points of $B_{c}$, or lines of the radial projection from the center of $P$ connecting points in the concentric copies of $B_{c}$,  as symbolically shown in Figure 2; in such cases we get multi-layer configurations.
\end{remark}

\definecolor{uuuuuu}{rgb}{0.26666666666666666,0.26666666666666666,0.26666666666666666}
\definecolor{zzttqq}{rgb}{0.6,0.2,0.0}
\definecolor{xdxdff}{rgb}{0.49019607843137253,0.49019607843137253,1.0}
\definecolor{qqqqff}{rgb}{0.0,0.0,1.0}
\begin{tikzpicture}[line cap=round,line join=round,>=triangle 45,x=0.25cm,y=0.25cm]
\clip(-14.375021130749278,-6.027465949450025) rectangle (19.134128610291505,9.336687769989792);
\fill[color=zzttqq,fill=zzttqq,fill opacity=0.1] (-4.487116383681567,1.34) -- (-6.75075626763347,2.6736140464892273) -- (-6.7738799687609506,0.046437378722607114) -- cycle;
\fill[color=zzttqq,fill=zzttqq,fill opacity=0.1] (-2.7626556637206616,1.3172029956387876) -- (-7.612218374825833,4.180932431047255) -- (-7.667499459902249,-1.450776791719792) -- cycle;
\fill[color=zzttqq,fill=zzttqq,fill opacity=0.1] (-1.19313458318255,1.2964542603942681) -- (-8.382073903249925,5.527964856221762) -- (-8.452199915585876,-2.8135945191351612) -- cycle;
\fill[color=zzttqq,fill=zzttqq,fill opacity=0.1] (8.116883195997438,-0.08628143059347736) -- (10.61159049470153,-0.08628143059347736) -- (10.611590494701531,2.408425868110613) -- (8.11688319599744,2.4084258681106143) -- cycle;
\fill[color=zzttqq,fill=zzttqq,fill opacity=0.1] (6.788180395600693,-1.550566149398052) -- (12.075875213506105,-1.5234497657164858) -- (12.048758829824541,3.7642450521889246) -- (6.761064011919128,3.7371286685073604) -- cycle;
\draw(-6.0,1.36) circle (1.5130158084186498cm);
\draw [color=zzttqq] (-4.487116383681567,1.34)-- (-6.75075626763347,2.6736140464892273);
\draw [color=zzttqq] (-6.75075626763347,2.6736140464892273)-- (-6.7738799687609506,0.046437378722607114);
\draw [color=zzttqq] (-6.7738799687609506,0.046437378722607114)-- (-4.487116383681567,1.34);
\draw [domain=-14.375021130749278:-6.0] plot(\x,{(--6.860655754953845--1.3136140464892272*\x)/-0.75075626763347});
\draw [domain=-14.375021130749278:-6.0] plot(\x,{(-8.933852485179251-1.313562621277393*\x)/-0.7738799687609506});
\draw [color=zzttqq] (-2.7626556637206616,1.3172029956387876)-- (-7.612218374825833,4.180932431047255);
\draw [color=zzttqq] (-7.612218374825833,4.180932431047255)-- (-7.667499459902249,-1.450776791719792);
\draw [color=zzttqq] (-7.667499459902249,-1.450776791719792)-- (-2.7626556637206616,1.3172029956387876);
\draw [color=zzttqq] (-1.19313458318255,1.2964542603942681)-- (-8.382073903249925,5.527964856221762);
\draw [color=zzttqq] (-8.382073903249925,5.527964856221762)-- (-8.452199915585876,-2.8135945191351612);
\draw [color=zzttqq] (-8.452199915585876,-2.8135945191351612)-- (-1.19313458318255,1.2964542603942681);
\draw (-6.0,1.36)-- (0.6869940672482949,1.2966541371663987);
\draw [color=zzttqq] (8.116883195997438,-0.08628143059347736)-- (10.61159049470153,-0.08628143059347736);
\draw [color=zzttqq] (10.61159049470153,-0.08628143059347736)-- (10.611590494701531,2.408425868110613);
\draw [color=zzttqq] (10.611590494701531,2.408425868110613)-- (8.11688319599744,2.4084258681106143);
\draw [color=zzttqq] (8.11688319599744,2.4084258681106143)-- (8.116883195997438,-0.08628143059347736);
\draw(9.364236845349485,1.1610722187585687) circle (1.7640244479892377cm);
\draw (8.11688319599744,2.4084258681106143)-- (10.61159049470153,-0.08628143059347736);
\draw (10.61159049470153,-0.08628143059347736)-- (12.075875213506105,-1.5234497657164858);
\draw (8.116883195997438,-0.08628143059347736)-- (6.788180395600693,-1.550566149398052);
\draw [color=zzttqq] (6.788180395600693,-1.550566149398052)-- (12.075875213506105,-1.5234497657164858);
\draw [color=zzttqq] (12.075875213506105,-1.5234497657164858)-- (12.048758829824541,3.7642450521889246);
\draw [color=zzttqq] (12.048758829824541,3.7642450521889246)-- (6.761064011919128,3.7371286685073604);
\draw [color=zzttqq] (6.761064011919128,3.7371286685073604)-- (6.788180395600693,-1.550566149398052);
\draw (12.048758829824541,3.7642450521889246)-- (10.611590494701531,2.408425868110613);
\draw (8.11688319599744,2.4084258681106143)-- (6.761064011919128,3.7371286685073604);
\draw (9.364236845349485,1.1610722187585687)-- (8.116883195997438,-0.08628143059347736);
\draw (9.364236845349485,1.1610722187585687)-- (10.611590494701531,2.408425868110613);
\draw (6.761064011919128,3.7371286685073604)-- (5.188313758388286,5.445460840446027);
\draw (12.048758829824541,3.7642450521889246)-- (13.62150908335538,5.445460840446027);
\draw (6.788180395600693,-1.550566149398052)-- (5.269662909432985,-3.2860147050182884);
\draw (12.075875213506105,-1.5234497657164858)-- (13.594392699673813,-3.231781937655156);
\begin{scriptsize}
\draw [fill=qqqqff] (-4.487116383681567,1.34) circle (1.5pt);
\draw [fill=xdxdff] (-6.75075626763347,2.6736140464892273) circle (1.5pt);
\draw [fill=uuuuuu] (-6.7738799687609506,0.046437378722607114) circle (1.5pt);
\draw [fill=xdxdff] (-2.7626556637206616,1.3172029956387876) circle (1.5pt);
\draw [fill=xdxdff] (-7.612218374825833,4.180932431047255) circle (1.5pt);
\draw [fill=uuuuuu] (-7.667499459902249,-1.450776791719792) circle (1.5pt);
\draw [fill=xdxdff] (-1.19313458318255,1.2964542603942681) circle (1.5pt);
\draw [fill=xdxdff] (-8.382073903249925,5.527964856221762) circle (1.5pt);
\draw [fill=uuuuuu] (-8.452199915585876,-2.8135945191351612) circle (1.5pt);
\draw [fill=qqqqff] (8.116883195997438,-0.08628143059347736) circle (1.5pt);
\draw [fill=qqqqff] (10.61159049470153,-0.08628143059347736) circle (1.5pt);
\draw [fill=uuuuuu] (10.611590494701531,2.408425868110613) circle (1.5pt);
\draw [fill=uuuuuu] (8.11688319599744,2.4084258681106143) circle (1.5pt);
\draw [fill=qqqqff] (12.075875213506105,-1.5234497657164858) circle (1.5pt);
\draw [fill=qqqqff] (6.788180395600693,-1.550566149398052) circle (1.5pt);
\draw [fill=uuuuuu] (12.048758829824541,3.7642450521889246) circle (1.5pt);
\draw [fill=uuuuuu] (6.761064011919128,3.7371286685073604) circle (1.5pt);
\end{scriptsize}
\end{tikzpicture}

\end{center}
\begin{center}
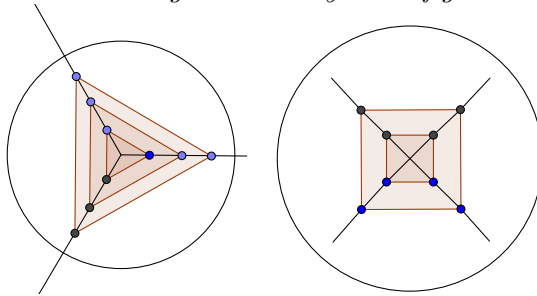
Figure 2. Radial projection from the center of $P$ through the concentric copies of $P$ increases the valences of those points, which are connected either by ``radial rays'' (left) or by ``antipodal lines'' (right), by 1.
\end{center}

Although we use the same general procedure (described in Algorithm 1) for all the five Platonic solids, and in some cases (subproblems of Problem 1) there are "universal" constructions (applicable to all the Platonic solids), some "special" constructions (applicable only to some of them) depend on special values of their parameters $(v,e,f,d,m)$, as it will be shown by examples in Section 3. 

The following simple calculation helps us to determine the number of points and lines in a 1-layer Platonic configuration $Z$.
 
\begin{prop}\label{calculation}
Let $Z$ be a 1-layer Platonic configuration, whose points and lines all lie on the boundary of the chosen Platonic solid $P = (v,e,f,d,m)$. Let $x,y,z,u,v$, respectively, denote the number of points of $Z$ in each vertex (represented by yellow points in our Figures), in the interior of each edge (red points) and in the interior of each face (blue points) of $P$, respectively. 
Then $x \in \lbrace 0,1 \rbrace$ and the number of points of $Z$ is $$p(Z) = xe + ye + zc.$$ Let $u,v$, respectively, denote the number of lines of $Z$ incident with each edge and with each face of $P$, respectively. Then $u \in \lbrace 0,1\rbrace$ and the number of lines in $Z$ is $$l(Z) = ue + vf.$$
Thus, for $P \in \lbrace T,C,O,D,I\rbrace$, respectively, we obtain the formulas:
 \begin{center}
$\begin{array}{ccc}
 P    & p(Z) & l(Z) \\
tetrahedron  & 4x + 6y + 4z & 6u + 4v \\ 
cube  & 8x + 12y + 6z & 12u + 6v \\ 
octahedron  & 6x + 12y + 8z & 12u + 8v \\ 
dodecahedron  & 20x + 30y + 12z & 30u + 12v  \\ 
icosahedron  & 12x + 30y + 20z & 30u + 20v \\
\end{array} $
\end{center}
 \end{prop}
 
 \begin{center}
Table 2: Numbers of points and lines in a 1-layer Platonic configuration $Z$.
\end{center}

\begin{remark}\label{similar}
Similar formulas may be obtained for the numbers of points and lines of $s$-layer Platonic configurations whose points and lines lie on $s$ "concentric copies"  of a Platonic solid $P$, connected by "radial" or "antipodal" lines. %The reader may find amusing to write down such formulas after going through the examples in the next section. 
Note that the number of points in a $s$-layer configuration $s \times Z$ must be $s$-times bigger than in $Z$, but the number of added "radial" (in $T$) or "antipodal" (in $C,O,D,I$) lines may depend on the number of points whose valence must be increased (some of the points may already have the desired valence $k$). 
In each construction with $s$-layer configurations special care must be devoted to the question which points of $Z$ have their antipodal points in $Z$ and which points do not have such antipodal points.

\end{remark}

\section{Examples of Platonic $(n_{k})$ configurations}\label{SectionExamples}

%In this section we solve Problem 1 by applying Algorithm 1 to get examples of Platonic $(n_{k})$ configurations of points and lines for $k = 3,4$ and 5.

In this section we apply Algorithm 1 to get examples of most of the Platonic $(p_{q},n_{k})$ configurations of points and lines for $k = 3,4$ and 5 and $3 \le q \le k$.

For any Platonic solid $P \in \lbrace T, C, O, D, I\rbrace$ let $P_{k}$ denote the class of all Platonic configuration $(n_{k})$ whose group of symmetry is the same as the full symmetry group of $P$; likewise let $P_{k,R}$ denote the class of all Platonic configurations $(n_{k})$ with the rotational symmetry group of $P$. 

Let $Cyc_{m}$ denote the class of all (planar) configurations, preserved by the cyclical group of rotations for the multiples of $\frac{2\pi}{m}$ and let $Dih_{m}$ denote the class of all (planar) configurations, preserved by the dihedral group $D_{m}$.

\begin{example}\label{exampletretji} 
To construct a Platonic configuration $Z \in T_{3,R}$ (with the rotational symmetries of a tetrahedron) place the copies of the planar configuration $A = (3_{2}3_{1},3_{3}) \in Cyc_{3}$  with 3-fold rotational symmetry (Figure 3 left) on the faces of $T$. You get a configuration $B = 4A$, which is not connected, but whose symmetry group is isomorphic to the group $Sym_{R}(T)$, the group of rotations of $T$. %Step 3. 
Place the 1-valent vertices of each copy of $A$ into vertices of the corresponding face. Now we have a connected configuration $B_{c} = (4_{3}12_{2}, 12_{3})$.
%Step 4. 
To get a balanced configuration $Z$, we have to increase the valence of 12 2-valent points in $B_{c}$ to 3. Take 3 concentric copies of $B_{c}$ to get $B_{1} = (12_{3}36_{2},27_{3})$. Add the 12 lines from the central point of $T$ through the $12$ triples of 2-valent points. Now we have $36 + 12 = 48$ 3-valent lines and $12 + 36 = 48$ 3-valent points. Hence we solved the problem of finding (a 3-layer) $Z \in T_{3,R}$ for $n = 48$. 
\end{example}

\vspace{-0mm}

\begin{center}

\definecolor{zzttqq}{rgb}{0.6,0.2,0.0}
\definecolor{ffffqq}{rgb}{1.0,1.0,0.0}
\definecolor{uuuuuu}{rgb}{0.26666666666666666,0.26666666666666666,0.26666666666666666}
\definecolor{ffqqqq}{rgb}{1.0,0.0,0.0}
\definecolor{xdxdff}{rgb}{0.49019607843137253,0.49019607843137253,1.0}
\definecolor{qqqqff}{rgb}{0.0,0.0,1.0}


\end{center}

\vspace{-0mm}

\begin{center}
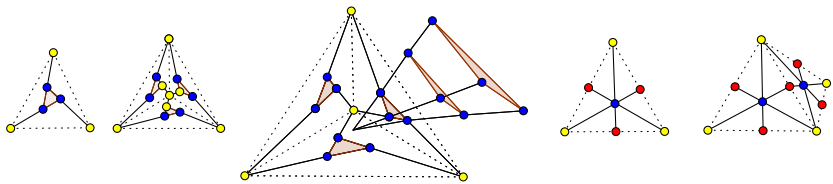
Figure 3. Construction of two tetrahedral Platonic configurations.
The valence of 1-valent points, placed into $m$-valent vertices, increases to $m$. 
\end{center}

\begin{prob} \label{TR} To construct a Platonic configuration $Z \in T_{3}$ of points and lines with the full symmetry group of the tetrahedron $T$.
\end{prob}

\textbf{Solution.} Subdivide baricentricaly each face of $T$ by the three symmetrals of its edges (see Figure 3 right). These $3 f = 12$ symmetrals are the lines of the configuration $B_{C}$ with $v + e + f = 4 + 6 + 4 = 14$ points. The four $v$-points and the 4 $f$-points have valence 3, the 6 $e$-points have valence 2. To increase their valence, take 3 concentric copies of $T$ and connect the triples of 2-valent points by 6 radial lines from the center of $T$. Thus we get $(14 \times 3)$ 3-valent points and  $12 \times 3 + 6 = 42$ 3-valent lines, hence a configuration $Z \in T(42_{3})$.

Many configurations may be obtained from the Pappus configuration $(9_{3})$.

\begin{prob} \label{pappusTOI} To construct a Platonic configuration $Z \in (n_{3})$ of points and lines with the full symmetry group of the tetrahedron $T$, or octahedron $O$, or icosahedron $I$. Do this by placing the copies of the same planar configuration $A$ on each of the faces of these solids.
\end{prob} 

\textbf{Solution}. \textit{Analysis}. Since $T$, $O$ and $I$ have different valences of vertices ($3,4$ and 5, respectively), no vertex of the solid $P$ can be included among the points of $Z$ ($x = 0$). Hence, to obtain a connected configuration $X$, the ``connecting points'' (shared by at least two copies of $A$) should lie on the edges of $P$.  If we want that all the edges of $P$ are included among the lines of $X$ (note that other options are possible, too),  there should be 3 connecting points on each edge (y = 3). We need at least 6 additional 3-valent points in the interior of each face (z = 6). Now we just ``connect the dots'' in the right ("symmetrical") way, and we find the Pappus configuration as a suitable starting point for building our desired Platonic configuration!

\begin{center}
\definecolor{ffqqqq}{rgb}{1.0,0.0,0.0}
\definecolor{uuuuuu}{rgb}{0.26666666666666666,0.26666666666666666,0.26666666666666666}
\definecolor{xdxdff}{rgb}{0.49019607843137253,0.49019607843137253,1.0}
\definecolor{qqqqff}{rgb}{0.0,0.0,1.0}

\end{center}

\begin{center}
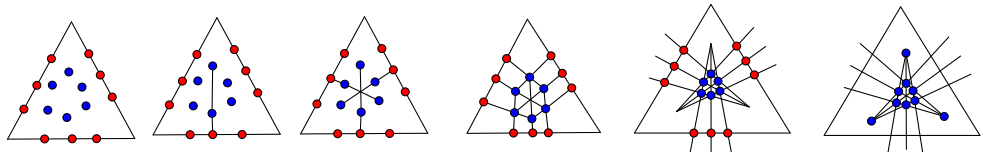
Figure 4. Analysis: connecting the dots we find the Pappus configuration as a possible starting configuration for our construction.
\end{center}

\textit{Synthesis.} Draw on each face of $P$ a copy of the Pappus configuration $A = (9_{3})$ in the 3-fold rotational symmetry form. The lines of these copies of $A$ intersect each edge of $P$ in 3 points. Include the lines through the edges of $P$ among the lines of configuration $X$. Exclude the vertices of $P$ and the 3 ``projecting '' points of each copy from the points of $X$ (as in Figure 5 right). 

\vspace{-2mm}

\begin{center}
\definecolor{uuuuuu}{rgb}{0.26666666666666666,0.26666666666666666,0.26666666666666666}
\definecolor{ffqqqq}{rgb}{1.0,0.0,0.0}
\definecolor{zzttqq}{rgb}{0.6,0.2,0.0}
\definecolor{qqqqff}{rgb}{0.0,0.0,1.0}
\begin{tikzpicture}[line cap=round,line join=round,>=triangle 45,x=0.5cm,y=0.5cm]
\clip(-3.84,-0.3000000000000012) rectangle (16.4,5.26);
\fill[color=zzttqq,fill=zzttqq,fill opacity=0.1] (3.7,0.62) -- (8.44,0.66) -- (6.035358983848624,4.74496041393824) -- cycle;
\fill[color=zzttqq,fill=zzttqq,fill opacity=0.1] (9.66,0.62) -- (14.4,0.66) -- (11.995358983848625,4.74496041393824) -- cycle;
\fill[color=zzttqq,fill=zzttqq,fill opacity=0.1] (11.995358983848625,4.744960413938242) -- (14.4,0.6600000000000001) -- (15.840000000000003,4.0600000000000005) -- cycle;
\fill[color=zzttqq,fill=zzttqq,fill opacity=0.1] (11.995358983848625,4.744960413938242) -- (14.4,0.6600000000000001) -- (15.840000000000003,4.0600000000000005) -- cycle;
\draw (-3.28,0.52)-- (0.25767949192431183,2.60248020696912);
\draw (-3.28,0.52)-- (-0.34348076211353223,3.6237203104536797);
\draw (0.858839745962156,1.58124010348456)-- (-3.28,0.52);
\draw (-0.9446410161513763,4.644960413938239)-- (-0.9099999999999999,0.54);
\draw (-0.9446410161513763,4.644960413938239)-- (-2.0949999999999998,0.53);
\draw (0.275,0.55)-- (-0.9446410161513763,4.644960413938239);
\draw (1.46,0.56)-- (-2.696160254037844,1.5512401034845598);
\draw (-1.5284807621135323,3.6137203104536795)-- (1.46,0.56);
\draw (1.46,0.56)-- (-2.112320508075688,2.5824802069691195);
\draw [color=zzttqq] (3.7,0.62)-- (8.44,0.66);
\draw [color=zzttqq] (8.44,0.66)-- (6.035358983848624,4.74496041393824);
\draw [color=zzttqq] (6.035358983848624,4.74496041393824)-- (3.7,0.62);
\draw (3.7,0.62)-- (7.237679491924312,2.70248020696912);
\draw (3.7,0.62)-- (6.636519237886468,3.72372031045368);
\draw (7.838839745962156,1.68124010348456)-- (3.7,0.62);
\draw (6.035358983848624,4.74496041393824)-- (6.07,0.64);
\draw (6.035358983848624,4.74496041393824)-- (4.885,0.63);
\draw (7.255,0.65)-- (6.035358983848624,4.74496041393824);
\draw (8.44,0.66)-- (4.283839745962156,1.65124010348456);
\draw (5.451519237886468,3.71372031045368)-- (8.44,0.66);
\draw (8.44,0.66)-- (4.867679491924312,2.68248020696912);
\draw [color=zzttqq] (9.66,0.62)-- (14.4,0.66);
\draw [color=zzttqq] (14.4,0.66)-- (11.995358983848625,4.74496041393824);
\draw [color=zzttqq] (11.995358983848625,4.74496041393824)-- (9.66,0.62);
\draw (9.66,0.62)-- (13.197679491924312,2.70248020696912);
\draw (9.66,0.62)-- (12.596519237886469,3.72372031045368);
\draw (13.798839745962155,1.68124010348456)-- (9.66,0.62);
\draw (11.995358983848625,4.74496041393824)-- (12.030000000000001,0.64);
\draw (11.995358983848625,4.74496041393824)-- (10.845,0.63);
\draw (13.215,0.65)-- (11.995358983848625,4.74496041393824);
\draw (14.4,0.66)-- (10.243839745962156,1.65124010348456);
\draw (11.41151923788647,3.71372031045368)-- (14.4,0.66);
\draw (14.4,0.66)-- (10.827679491924313,2.68248020696912);
\draw (11.995358983848625,4.744960413938242)-- (15.84,4.06);
\draw (15.84,4.06)-- (14.4,0.6600000000000001);
\draw (13.197679491924312,2.70248020696912)-- (15.84,4.06);
\draw (13.798839745962155,1.68124010348456)-- (15.84,4.06);
\draw (12.596519237886469,3.72372031045368)-- (15.84,4.06);
\draw (11.995358983848625,4.744960413938242)-- (15.48,3.21);
\draw (11.995358983848625,4.744960413938242)-- (15.120000000000001,2.36);
\draw (11.995358983848625,4.744960413938242)-- (14.760000000000002,1.51);
\draw (14.4,0.6600000000000001)-- (14.878839745962157,4.23124010348456);
\draw (14.4,0.6600000000000001)-- (13.917679491924313,4.402480206969121);
\draw (14.4,0.6600000000000001)-- (12.956519237886468,4.573720310453681);
\draw [color=ffqqqq] (3.6999999999999997,0.6199999999999999)-- (8.440000000000001,0.6599999999999999);
\draw [color=ffqqqq] (8.440000000000001,0.6599999999999999)-- (6.035358983848625,4.744960413938241);
\draw [color=ffqqqq] (6.035358983848625,4.744960413938241)-- (3.6999999999999997,0.6199999999999999);
\draw [color=ffqqqq] (9.66,0.6200000000000002)-- (14.4,0.6600000000000001);
\draw [color=ffqqqq] (14.4,0.6600000000000001)-- (15.84,4.06);
\draw [color=ffqqqq] (9.66,0.6200000000000002)-- (11.995358983848625,4.744960413938242);
\draw [color=ffqqqq] (11.995358983848625,4.744960413938242)-- (14.4,0.6600000000000001);
\draw [color=ffqqqq] (11.995358983848625,4.744960413938242)-- (15.84,4.06);
\draw (0.7,5.68)-- (4.44,6.2);
\draw [color=zzttqq] (11.995358983848625,4.744960413938242)-- (14.4,0.6600000000000001);
\draw [color=zzttqq] (14.4,0.6600000000000001)-- (15.840000000000003,4.0600000000000005);
\draw [color=zzttqq] (15.840000000000003,4.0600000000000005)-- (11.995358983848625,4.744960413938242);
\draw [color=zzttqq] (11.995358983848625,4.744960413938242)-- (14.4,0.6600000000000001);
\draw [color=zzttqq] (14.4,0.6600000000000001)-- (15.840000000000003,4.0600000000000005);
\draw [color=zzttqq] (15.840000000000003,4.0600000000000005)-- (11.995358983848625,4.744960413938242);
\begin{scriptsize}
\draw [fill=qqqqff] (-1.6019890069220177,2.293554463116389) circle (1.5pt);
\draw [fill=qqqqff] (-1.8649282032302748,1.3529920827876476) circle (1.5pt);
\draw [fill=qqqqff] (-0.914948716593054,1.1264229162768915) circle (1.5pt);
\draw [fill=qqqqff] (0.031071796769724675,1.3689920827876478) circle (1.5pt);
\draw [fill=qqqqff] (-0.2477032926363041,2.3049830345449593) circle (1.5pt);
\draw [fill=qqqqff] (-0.9307846096908263,3.002976248362944) circle (1.5pt);
\draw [fill=ffqqqq] (7.237679491924312,2.70248020696912) circle (1.5pt);
\draw [fill=ffqqqq] (6.07,0.64) circle (1.5pt);
\draw [fill=ffqqqq] (4.867679491924312,2.68248020696912) circle (1.5pt);
\draw [fill=ffqqqq] (4.867679491924312,2.68248020696912) circle (1.5pt);
\draw [fill=ffqqqq] (5.451519237886468,3.71372031045368) circle (1.5pt);
\draw [fill=ffqqqq] (4.867679491924312,2.68248020696912) circle (1.5pt);
\draw [fill=ffqqqq] (4.867679491924312,2.68248020696912) circle (1.5pt);
\draw [fill=ffqqqq] (4.283839745962156,1.65124010348456) circle (1.5pt);
\draw [fill=ffqqqq] (4.885,0.63) circle (1.5pt);
\draw [fill=ffqqqq] (6.07,0.64) circle (1.5pt);
\draw [fill=ffqqqq] (7.255,0.65) circle (1.5pt);
\draw [fill=ffqqqq] (7.237679491924312,2.70248020696912) circle (1.5pt);
\draw [fill=ffqqqq] (7.838839745962156,1.68124010348456) circle (1.5pt);
\draw [fill=ffqqqq] (7.237679491924312,2.70248020696912) circle (1.5pt);
\draw [fill=ffqqqq] (7.237679491924312,2.70248020696912) circle (1.5pt);
\draw [fill=ffqqqq] (6.636519237886468,3.72372031045368) circle (1.5pt);
\draw [fill=qqqqff] (5.378010993077981,2.393554463116388) circle (1.5pt);
\draw [fill=qqqqff] (5.115071796769725,1.4529920827876484) circle (1.5pt);
\draw [fill=qqqqff] (6.065051283406945,1.2264229162768912) circle (1.5pt);
\draw [fill=qqqqff] (7.011071796769725,1.4689920827876486) circle (1.5pt);
\draw [fill=qqqqff] (6.732296707363698,2.4049830345449585) circle (1.5pt);
\draw [fill=qqqqff] (6.049215390309174,3.102976248362944) circle (1.5pt);
\draw [fill=uuuuuu] (-3.2799999999999994,0.52) circle (1.5pt);
\draw [fill=uuuuuu] (1.4600000000000002,0.5599999999999999) circle (1.5pt);
\draw [fill=uuuuuu] (-0.9446410161513765,4.644960413938239) circle (1.5pt);
\draw [fill=ffqqqq] (13.197679491924312,2.70248020696912) circle (1.5pt);
\draw [fill=ffqqqq] (12.030000000000001,0.64) circle (1.5pt);
\draw [fill=ffqqqq] (10.827679491924313,2.68248020696912) circle (1.5pt);
\draw [fill=ffqqqq] (10.827679491924313,2.68248020696912) circle (1.5pt);
\draw [fill=ffqqqq] (11.41151923788647,3.71372031045368) circle (1.5pt);
\draw [fill=ffqqqq] (10.827679491924313,2.68248020696912) circle (1.5pt);
\draw [fill=ffqqqq] (10.827679491924313,2.68248020696912) circle (1.5pt);
\draw [fill=ffqqqq] (10.243839745962156,1.65124010348456) circle (1.5pt);
\draw [fill=ffqqqq] (10.845,0.63) circle (1.5pt);
\draw [fill=ffqqqq] (12.030000000000001,0.64) circle (1.5pt);
\draw [fill=ffqqqq] (13.215,0.65) circle (1.5pt);
\draw [fill=ffqqqq] (13.197679491924312,2.70248020696912) circle (1.5pt);
\draw [fill=ffqqqq] (13.798839745962155,1.68124010348456) circle (1.5pt);
\draw [fill=ffqqqq] (13.197679491924312,2.70248020696912) circle (1.5pt);
\draw [fill=ffqqqq] (13.197679491924312,2.70248020696912) circle (1.5pt);
\draw [fill=ffqqqq] (12.596519237886469,3.72372031045368) circle (1.5pt);
\draw [fill=qqqqff] (11.338010993077983,2.393554463116388) circle (1.5pt);
\draw [fill=qqqqff] (11.075071796769725,1.4529920827876475) circle (1.5pt);
\draw [fill=qqqqff] (12.025051283406944,1.2264229162768912) circle (1.5pt);
\draw [fill=qqqqff] (12.971071796769726,1.4689920827876484) circle (1.5pt);
\draw [fill=qqqqff] (12.692296707363699,2.404983034544959) circle (1.5pt);
\draw [fill=qqqqff] (12.009215390309173,3.102976248362944) circle (1.5pt);
\draw [fill=ffqqqq] (13.917679491924313,4.402480206969121) circle (1.5pt);
\draw [fill=ffqqqq] (12.956519237886468,4.573720310453681) circle (1.5pt);
\draw [fill=ffqqqq] (14.878839745962157,4.23124010348456) circle (1.5pt);
\draw [fill=ffqqqq] (12.030000000000001,0.64) circle (1.5pt);
\draw [fill=ffqqqq] (15.120000000000001,2.36) circle (1.5pt);
\draw [fill=ffqqqq] (14.760000000000002,1.51) circle (1.5pt);
\draw [fill=ffqqqq] (15.48,3.21) circle (1.5pt);
\draw [fill=qqqqff] (13.245215390309175,3.7909762483629432) circle (1.5pt);
\draw [fill=qqqqff] (13.98658242164941,3.8678401774021007) circle (1.5pt);
\draw [fill=qqqqff] (14.783071796769724,3.516992082787648) circle (1.5pt);
\draw [fill=qqqqff] (14.673622711978375,2.7007086305626045) circle (1.5pt);
\draw [fill=qqqqff] (14.207071796769725,2.15699208278765) circle (1.5pt);
\draw [fill=uuuuuu] (13.57515385022084,2.896411605973534) circle (1.5pt);
\draw [fill=qqqqff] (0.7,5.68) circle (1.5pt);
\draw [fill=qqqqff] (4.44,6.2) circle (1.5pt);
\end{scriptsize}
\end{tikzpicture}
\end{center}

\vspace{-3mm}

\begin{center}
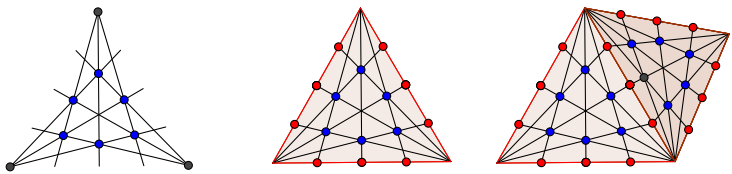
Figure 5. Synthesis: starting with the Pappus configuration $A = (9_{3})$ placed on each face of $P$ we get a Platonic configuration.
\end{center}

On each face there are $z = 6$ interior points of valence 3, on each edge there are $y = 3$ points of valence 3. So (directly or by Proposition 1) we get $6f + 3e = 4e + 3e = 7e$ points of valence 3, and $9f + e = 6e + e = 7e$ lines of valence 3, and the obtained Platonic configuration is balanced $(7e_{3}) = (\frac{21f}{2})$.  The same construction works for any solid $P$ such that all of its faces are equilateral triangles, hence $3f = 2e$. 

Very often a solution of one problem leads to a solution of another problem, yet some modifications in the construction are needed. 
The Pappus configuration lies in the background of the following construction, too.

\begin{example}\label{examplecetrti}
To construct a Platonic configuration $Z \in D_{3}$ divide each pentagonal face of a dodecahedron into 5 congruent triangles and draw into each triangle the same configuration as in Figure 5.
The obtained configuration has $(6 + 3) \times (5 \times f) = 9 \times 60 = 540$ 3-valent points in the interiors of faces of $D$ and $3 \times e = 90$ 3-valent points on the edges of $D$, and $(9  + 5)\times f  + e = 540 + 60  + 30 = 630$ 3-valent lines. Thus we get a $(630_{3})$ configuration $ Z \in D_{3}$  with the full symmetry group of the dodecahedron. 
\end{example}

\begin{center}

\definecolor{xdxdff}{rgb}{0.49019607843137253,0.49019607843137253,1.0}
\definecolor{uuuuuu}{rgb}{0.26666666666666666,0.26666666666666666,0.26666666666666666}
\definecolor{zzttqq}{rgb}{0.6,0.2,0.0}
\definecolor{qqqqff}{rgb}{0.0,0.0,1.0}
\definecolor{ffqqqq}{rgb}{1.0,0.0,0.0}
\begin{tikzpicture}[line cap=round,line join=round,>=triangle 45,x=0.5cm,y=0.5cm]
\clip(-3.3727272727272717,-2.6145454545454547) rectangle (7.627272727272728,6.112727272727272);
\fill[color=zzttqq,fill=zzttqq,fill opacity=0.1] (-1.64,1.1) -- (1.8400000000000003,1.1) -- (0.3199999999999999,3.6000000000000005) -- cycle;
\fill[color=zzttqq,fill=zzttqq,fill opacity=0.1] (10.427272727272726,1.2763636363636361) -- (12.572727272727272,1.2581818181818178) -- (13.253001124591618,3.293011853244603) -- (11.527978941247047,4.568787794424364) -- (9.781582748728205,3.3224306530400574) -- cycle;
\draw (-1.64,1.1)-- (1.84,1.1);
\draw (1.84,1.1)-- (3.42,3.92);
\draw (0.5,5.9)-- (3.42,3.92);
\draw (-2.76,4.04)-- (0.5,5.9);
\draw (-2.76,4.04)-- (-1.64,1.1);
\draw (-1.64,1.1)-- (0.32,3.6);
\draw (0.32,3.6)-- (1.84,1.1);
\draw (0.32,3.6)-- (3.42,3.92);
\draw (0.32,3.6)-- (0.5,5.9);
\draw (0.32,3.6)-- (-2.76,4.04);
\draw (3.42,3.92)-- (6.14,3.94);
\draw (6.14,3.94)-- (6.56,0.34);
\draw (6.56,0.34)-- (3.48,-0.86);
\draw (1.84,1.1)-- (3.48,-0.86);
\draw (-1.64,1.1)-- (-2.02,-1.64);
\draw (1.04,-2.48)-- (3.48,-0.86);
\draw (-2.02,-1.64)-- (1.04,-2.48);
\draw (0.7209345357693566,2.9405681977477687)-- (0.28,2.54);
\draw (-0.027272727272726238,2.3127272727272725)-- (0.28,2.54);
\draw (-0.17899242934717985,2.963530064608189)-- (0.28,2.54);
\draw (0.5909090909090919,2.203636363636363)-- (0.28,2.54);
\draw (0.2272727272727283,1.5854545454545452)-- (-0.3,1.78);
\draw (-0.3,1.78)-- (-0.917109675373578,2.022053985492885);
\draw (0.2272727272727283,1.5854545454545452)-- (0.7181818181818191,1.7309090909090905);
\draw (0.7181818181818191,1.7309090909090905)-- (1.3267753843278354,1.944119433671324);
\draw (0.28,2.54)-- (0.2272727272727283,1.5854545454545452);
\draw (0.2272727272727283,1.5854545454545452)-- (0.19090909090909203,1.1);
\draw (-0.5799999999999998,1.1)-- (-0.3,1.78);
\draw (-0.3,1.78)-- (-0.027272727272726238,2.3127272727272725);
\draw (0.5909090909090919,2.203636363636363)-- (0.7181818181818191,1.7309090909090905);
\draw (0.7181818181818191,1.7309090909090905)-- (0.9000000000000001,1.1);
\draw (-0.5761353997383964,2.4569701533949027)-- (-0.027272727272726238,2.3127272727272725);
\draw (-0.027272727272726238,2.3127272727272725)-- (0.5909090909090919,2.203636363636363);
\draw (1.044025980094388,2.409167795897388)-- (0.5909090909090919,2.203636363636363);
\draw (0.5909090909090919,2.203636363636363)-- (-0.3,1.78);
\draw (-0.027272727272726238,2.3127272727272725)-- (0.7181818181818191,1.7309090909090905);
\draw (-1.64,1.1)-- (0.36,-0.44);
\draw (0.36,-0.44)-- (-2.02,-1.6399999999999995);
\draw (0.36,-0.44)-- (1.0399999999999998,-2.48);
\draw (0.36,-0.44)-- (3.48,-0.86);
\draw (0.36,-0.44)-- (1.8400000000000003,1.1);
\draw [color=zzttqq] (-1.64,1.1)-- (1.8400000000000003,1.1);
\draw [color=zzttqq] (1.8400000000000003,1.1)-- (0.3199999999999999,3.6000000000000005);
\draw [color=zzttqq] (0.3199999999999999,3.6000000000000005)-- (-1.64,1.1);
\draw (0.7209345357693566,2.9405681977477687)-- (1.3,3.149090909090909);
\draw (1.3,3.149090909090909)-- (1.9,3.44);
\draw (1.3267753843278354,1.944119433671324)-- (1.82,2.34);
\draw (1.82,2.34)-- (2.3363636363636373,2.8218181818181827);
\draw (1.4272727272727281,2.7127272727272733)-- (1.82,2.34);
\draw (1.82,2.34)-- (2.2853046195822406,1.894784194444252);
\draw (2.732522777735242,2.6929836919072043)-- (2.3363636363636373,2.8218181818181827);
\draw (2.3363636363636373,2.8218181818181827)-- (1.3,3.149090909090909);
\draw (3.0448914533900373,3.2505024674429777)-- (2.40909090909091,3.367272727272728);
\draw (2.40909090909091,3.367272727272728)-- (1.9,3.44);
\draw (2.3363636363636373,2.8218181818181827)-- (2.40909090909091,3.367272727272728);
\draw (2.40909090909091,3.367272727272728)-- (2.459907746797908,3.820893702895268);
\draw (1.105147009820623,3.6810474332718064)-- (1.3,3.149090909090909);
\draw (1.3,3.149090909090909)-- (1.4272727272727281,2.7127272727272733);
\draw (1.82,2.34)-- (1.9,3.44);
\draw (1.9,3.44)-- (1.9235583377949836,3.765528602611095);
\draw (1.4272727272727281,2.7127272727272733)-- (1.044025980094388,2.409167795897388);
\draw (1.4272727272727281,2.7127272727272733)-- (2.40909090909091,3.367272727272728);
\draw [color=zzttqq] (10.427272727272726,1.2763636363636361)-- (12.572727272727272,1.2581818181818178);
\draw [color=zzttqq] (12.572727272727272,1.2581818181818178)-- (13.253001124591618,3.293011853244603);
\draw [color=zzttqq] (13.253001124591618,3.293011853244603)-- (11.527978941247047,4.568787794424364);
\draw [color=zzttqq] (11.527978941247047,4.568787794424364)-- (9.781582748728205,3.3224306530400574);
\draw [color=zzttqq] (9.781582748728205,3.3224306530400574)-- (10.427272727272726,1.2763636363636361);
\draw (0.36224710656846537,4.139824139485946)-- (1.04,4.22);
\draw (1.04,4.22)-- (1.105147009820623,3.6810474332718064);
\draw (0.4076055914624981,4.7194047797985865)-- (0.82,4.6);
\draw (0.4530858259431835,5.300541109274011)-- (0.92,5.1);
\draw (0.92,5.1)-- (1.2895645467108001,5.3646103416139095);
\draw (1.46,4.92)-- (1.7352077642446249,5.062427611916316);
\draw (1.86,4.46)-- (2.1710916862165375,4.76686248674358);
\draw (1.86,4.46)-- (2.459907746797908,3.820893702895268);
\draw (1.54,4.06)-- (1.9235583377949836,3.765528602611095);
\draw (1.54,4.06)-- (1.86,4.46);
\draw (1.86,4.46)-- (1.46,4.92);
\draw (1.46,4.92)-- (0.92,5.1);
\draw (0.92,5.1)-- (0.82,4.6);
\draw (0.82,4.6)-- (1.04,4.22);
\draw (1.54,4.06)-- (1.04,4.22);
\draw (0.92,5.1)-- (1.54,4.06);
\draw (1.04,4.22)-- (1.46,4.92);
\draw (0.82,4.6)-- (1.86,4.46);
\draw (-0.5761353997383964,2.4569701533949027)-- (-0.9,2.6581818181818178);
\draw (-0.917109675373578,2.022053985492885)-- (-1.46,2.34);
\draw (-0.17899242934717985,2.963530064608189)-- (-0.609090909090908,3.3127272727272725);
\draw (-2.07659405940594,2.246059405940594)-- (-1.46,2.34);
\draw (-2.2774653465346533,2.773346534653465)-- (-1.7181818181818171,2.8218181818181813);
\draw (-2.4758019801980193,3.293980198019802)-- (-1.76,3.34);
\draw (-1.9451999999999996,3.9236)-- (-1.76,3.34);
\draw (-1.3264,3.8352)-- (-1.12,3.54);
\draw (-0.7887999999999995,3.7584)-- (-0.609090909090908,3.3127272727272725);
\draw (-1.318181818181817,4.349090909090909)-- (-0.7545454545454535,4.076363636363636);
\draw (-0.19090909090908986,4.0218181818181815)-- (-0.7545454545454535,4.076363636363636);
\draw (-0.19090909090908986,4.0218181818181815)-- (-0.027272727272726238,4.6036363636363635);
\draw (-0.027272727272726238,4.6036363636363635)-- (-0.17272727272727167,5.058181818181818);
\draw (-0.17272727272727167,5.058181818181818)-- (-0.9727272727272717,4.694545454545454);
\draw (-0.9727272727272717,4.694545454545454)-- (-1.318181818181817,4.349090909090909);
\begin{scriptsize}
\draw [fill=ffqqqq] (1.2895645467108001,5.3646103416139095) circle (1.5pt);
\draw [fill=ffqqqq] (1.7352077642446249,5.062427611916316) circle (1.5pt);
\draw [fill=ffqqqq] (2.1710916862165375,4.76686248674358) circle (1.5pt);
\draw [fill=ffqqqq] (1.105147009820623,3.6810474332718064) circle (1.5pt);
\draw [fill=ffqqqq] (1.9235583377949836,3.765528602611095) circle (1.5pt);
\draw [fill=ffqqqq] (2.459907746797908,3.820893702895268) circle (1.5pt);
\draw [fill=ffqqqq] (0.36224710656846537,4.139824139485946) circle (1.5pt);
\draw [fill=ffqqqq] (0.4076055914624981,4.7194047797985865) circle (1.5pt);
\draw [fill=ffqqqq] (0.4530858259431835,5.300541109274011) circle (1.5pt);
\draw [fill=ffqqqq] (-1.8184155829405415,4.577223010960305) circle (1.5pt);
\draw [fill=ffqqqq] (-1.2345073541938782,4.910373104662389) circle (1.5pt);
\draw [fill=ffqqqq] (-0.7195610199329887,5.204176841387927) circle (1.5pt);
\draw [fill=ffqqqq] (-1.9451999999999996,3.9236) circle (1.5pt);
\draw [fill=ffqqqq] (-1.3264,3.8352) circle (1.5pt);
\draw [fill=ffqqqq] (-0.7887999999999995,3.7584) circle (1.5pt);
\draw [fill=ffqqqq] (-0.17899242934717985,2.963530064608189) circle (1.5pt);
\draw [fill=ffqqqq] (-0.5761353997383964,2.4569701533949027) circle (1.5pt);
\draw [fill=ffqqqq] (-0.917109675373578,2.022053985492885) circle (1.5pt);
\draw [fill=ffqqqq] (-2.4758019801980193,3.293980198019802) circle (1.5pt);
\draw [fill=ffqqqq] (-2.2774653465346533,2.773346534653465) circle (1.5pt);
\draw [fill=ffqqqq] (-2.07659405940594,2.246059405940594) circle (1.5pt);
\draw [fill=ffqqqq] (-0.5799999999999998,1.1) circle (1.5pt);
\draw [fill=ffqqqq] (0.19090909090909203,1.1) circle (1.5pt);
\draw [fill=ffqqqq] (0.9000000000000001,1.1) circle (1.5pt);
\draw [fill=ffqqqq] (0.7209345357693566,2.9405681977477687) circle (1.5pt);
\draw [fill=ffqqqq] (1.044025980094388,2.409167795897388) circle (1.5pt);
\draw [fill=ffqqqq] (1.3267753843278354,1.944119433671324) circle (1.5pt);
\draw [fill=ffqqqq] (2.2853046195822406,1.894784194444252) circle (1.5pt);
\draw [fill=ffqqqq] (2.732522777735242,2.6929836919072043) circle (1.5pt);
\draw [fill=ffqqqq] (3.0448914533900373,3.2505024674429777) circle (1.5pt);
\draw [fill=ffqqqq] (2.3148848603625676,0.5324546790788829) circle (1.5pt);
\draw [fill=ffqqqq] (2.6419206271435574,0.14160705536501705) circle (1.5pt);
\draw [fill=ffqqqq] (3.009936305732484,-0.2982165605095539) circle (1.5pt);
\draw [fill=ffqqqq] (4.417891345731439,-0.49458778737736137) circle (1.5pt);
\draw [fill=ffqqqq] (5.102987260213794,-0.22766730121540474) circle (1.5pt);
\draw [fill=ffqqqq] (5.781317908917851,0.03661736711085112) circle (1.5pt);
\draw [fill=ffqqqq] (6.245201425047958,3.0382734995889287) circle (1.5pt);
\draw [fill=ffqqqq] (6.329416278432447,2.316431899150452) circle (1.5pt);
\draw [fill=ffqqqq] (6.422570567278706,1.517966566182516) circle (1.5pt);
\draw [fill=ffqqqq] (4.259660485484132,3.9261739741579715) circle (1.5pt);
\draw [fill=ffqqqq] (5.019619397740174,3.9317619073363246) circle (1.5pt);
\draw [fill=ffqqqq] (5.539885386819484,3.935587392550143) circle (1.5pt);
\draw [fill=ffqqqq] (-1.73947935180345,0.38270151594354407) circle (1.5pt);
\draw [fill=ffqqqq] (-1.8032033455305803,-0.07678201777313132) circle (1.5pt);
\draw [fill=ffqqqq] (-1.8816664924202822,-0.6425426032409827) circle (1.5pt);
\draw [fill=ffqqqq] (-1.3942152306042186,-1.8117840543439399) circle (1.5pt);
\draw [fill=ffqqqq] (-0.700600643546657,-2.002188058634251) circle (1.5pt);
\draw [fill=ffqqqq] (-0.13899892742223807,-2.156353235609582) circle (1.5pt);
\draw [fill=ffqqqq] (1.6129932385171368,-2.099570062951737) circle (1.5pt);
\draw [fill=ffqqqq] (2.1998675681977153,-1.7099239916064348) circle (1.5pt);
\draw [fill=ffqqqq] (2.814504080205176,-1.3018456516670551) circle (1.5pt);
\draw [fill=qqqqff] (0.28,2.54) circle (1.5pt);
\draw [fill=qqqqff] (-0.027272727272726238,2.3127272727272725) circle (1.5pt);
\draw [fill=qqqqff] (-0.3,1.78) circle (1.5pt);
\draw [fill=qqqqff] (0.2272727272727283,1.5854545454545452) circle (1.5pt);
\draw [fill=qqqqff] (0.7181818181818191,1.7309090909090905) circle (1.5pt);
\draw [fill=qqqqff] (0.5909090909090919,2.203636363636363) circle (1.5pt);
\draw [fill=qqqqff] (1.3,3.149090909090909) circle (1.5pt);
\draw [fill=qqqqff] (1.4272727272727281,2.7127272727272733) circle (1.5pt);
\draw [fill=qqqqff] (1.82,2.34) circle (1.5pt);
\draw [fill=qqqqff] (2.3363636363636373,2.8218181818181827) circle (1.5pt);
\draw [fill=qqqqff] (2.40909090909091,3.367272727272728) circle (1.5pt);
\draw [fill=qqqqff] (1.9,3.44) circle (1.5pt);
\draw [fill=qqqqff] (0.92,5.1) circle (1.5pt);
\draw [fill=qqqqff] (0.82,4.6) circle (1.5pt);
\draw [fill=qqqqff] (1.04,4.22) circle (1.5pt);
\draw [fill=qqqqff] (1.54,4.06) circle (1.5pt);
\draw [fill=qqqqff] (1.86,4.46) circle (1.5pt);
\draw [fill=qqqqff] (1.46,4.92) circle (1.5pt);
\draw [fill=qqqqff] (-0.19090909090908986,4.0218181818181815) circle (1.5pt);
\draw [fill=qqqqff] (-0.027272727272726238,4.6036363636363635) circle (1.5pt);
\draw [fill=qqqqff] (-0.17272727272727167,5.058181818181818) circle (1.5pt);
\draw [fill=qqqqff] (-0.9727272727272717,4.694545454545454) circle (1.5pt);
\draw [fill=qqqqff] (-1.318181818181817,4.349090909090909) circle (1.5pt);
\draw [fill=qqqqff] (-0.7545454545454535,4.076363636363636) circle (1.5pt);
\draw [fill=qqqqff] (-0.609090909090908,3.3127272727272725) circle (1.5pt);
\draw [fill=qqqqff] (-1.12,3.54) circle (1.5pt);
\draw [fill=qqqqff] (-1.76,3.34) circle (1.5pt);
\draw [fill=qqqqff] (-1.7181818181818171,2.8218181818181813) circle (1.5pt);
\draw [fill=qqqqff] (-1.46,2.34) circle (1.5pt);
\draw [fill=qqqqff] (-0.9,2.6581818181818178) circle (1.5pt);
\draw [fill=qqqqff] (10.427272727272726,1.2763636363636361) circle (1.5pt);
\draw [fill=qqqqff] (12.572727272727272,1.2581818181818178) circle (1.5pt);
\draw [fill=uuuuuu] (13.253001124591618,3.293011853244603) circle (1.5pt);
\draw [fill=uuuuuu] (11.527978941247047,4.568787794424364) circle (1.5pt);
\draw [fill=uuuuuu] (9.781582748728205,3.3224306530400574) circle (1.5pt);
\draw [fill=xdxdff] (-1.0127735005906884,0.6170355954548302) circle (1.5pt);
\draw [fill=xdxdff] (-0.6552780235020167,0.34176407809655296) circle (1.5pt);
\draw [fill=xdxdff] (-0.23808582402593248,0.020526084499968178) circle (1.5pt);
\draw [fill=xdxdff] (0.8067074249730986,0.02481718544498074) circle (1.5pt);
\draw [fill=xdxdff] (1.1364352158144357,0.3679123191582641) circle (1.5pt);
\draw [fill=xdxdff] (1.4385578892830102,0.6822832091188078) circle (1.5pt);
\draw [fill=xdxdff] (-1.354816733292044,-1.3046134789707782) circle (1.5pt);
\draw [fill=xdxdff] (-0.9193119756770446,-1.0850312482405264) circle (1.5pt);
\draw [fill=xdxdff] (-0.5126784425528857,-0.8800059374216229) circle (1.5pt);
\draw [fill=xdxdff] (0.5390909090909093,-0.9772727272727277) circle (1.5pt);
\draw [fill=xdxdff] (0.6427272727272728,-1.2881818181818185) circle (1.5pt);
\draw [fill=xdxdff] (0.7681818181818183,-1.664545454545455) circle (1.5pt);
\draw [fill=xdxdff] (1.299135488557937,-0.5664220849981838) circle (1.5pt);
\draw [fill=xdxdff] (1.6837737344384647,-0.618200310405178) circle (1.5pt);
\draw [fill=xdxdff] (2.224327840702706,-0.6909672093253643) circle (1.5pt);
\draw [fill=qqqqff] (-1.3545454545454534,0.4036363636363634) circle (1.5pt);
\draw [fill=qqqqff] (-1.4636363636363625,-0.12363636363636384) circle (1.5pt);
\draw [fill=qqqqff] (-1.409090909090908,-0.5963636363636365) circle (1.5pt);
\draw [fill=qqqqff] (-0.9181818181818171,-0.7236363636363639) circle (1.5pt);
\draw [fill=qqqqff] (-0.609090909090908,-0.30545454545454565) circle (1.5pt);
\draw [fill=qqqqff] (-0.9545454545454535,0.11272727272727251) circle (1.5pt);
\draw [fill=qqqqff] (-0.6272727272727262,0.767272727272727) circle (1.5pt);
\draw [fill=qqqqff] (-0.2636363636363626,0.34909090909090884) circle (1.5pt);
\draw [fill=qqqqff] (0.28181818181818286,0.07636363636363615) circle (1.5pt);
\draw [fill=qqqqff] (0.8090909090909101,0.44) circle (1.5pt);
\draw [fill=qqqqff] (0.9727272727272737,0.7854545454545452) circle (1.5pt);
\draw [fill=qqqqff] (0.19090909090909194,0.8581818181818179) circle (1.5pt);
\draw [fill=qqqqff] (0.009090909090910124,-0.8872727272727274) circle (1.5pt);
\draw [fill=qqqqff] (-0.8636363636363625,-1.5781818181818184) circle (1.5pt);
\draw [fill=qqqqff] (0.35454545454545555,-1.9054545454545455) circle (1.5pt);
\draw [fill=qqqqff] (-0.6272727272727262,-1.232727272727273) circle (1.5pt);
\draw [fill=qqqqff] (0.26363636363636467,-1.36) circle (1.5pt);
\draw [fill=qqqqff] (-0.46363636363636257,-1.76) circle (1.5pt);
\draw [fill=qqqqff] (1.0090909090909101,-0.8327272727272729) circle (1.5pt);
\draw [fill=qqqqff] (1.20909090909091,-1.7781818181818183) circle (1.5pt);
\draw [fill=qqqqff] (2.427272727272728,-0.9963636363636366) circle (1.5pt);
\draw [fill=qqqqff] (1.7181818181818191,-0.8509090909090911) circle (1.5pt);
\draw [fill=qqqqff] (0.9545454545454556,-1.3054545454545456) circle (1.5pt);
\draw [fill=qqqqff] (1.9363636363636374,-1.5236363636363637) circle (1.5pt);
\draw [fill=qqqqff] (1.2636363636363646,-0.23272727272727292) circle (1.5pt);
\draw [fill=qqqqff] (1.80909090909091,0.6763636363636362) circle (1.5pt);
\draw [fill=qqqqff] (2.6272727272727283,-0.39636363636363653) circle (1.5pt);
\draw [fill=qqqqff] (1.390909090909092,0.18545454545454523) circle (1.5pt);
\draw [fill=qqqqff] (1.790909090909092,-0.36) circle (1.5pt);
\draw [fill=qqqqff] (2.3181818181818192,0.09454545454545434) circle (1.5pt);
\end{scriptsize}
\end{tikzpicture}
\end{center}

\begin{center}
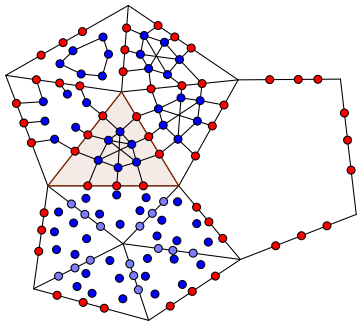
Figure 6. "Connecting the dots" in the triangulated pentagonal faces of $D$; as a result we obtain a configuration $Z \in D(630_{3})$.
\end{center}

\begin{prob}\label{D3}
To construct a Platonic configuration $Z \in D_{3,R}$ with the rotational symmetry group $Sym_{R}(D)$ of the dodecahedron.
\end{prob}

\textbf{Solution.} Draw on each face of $D$ the configuration $(15_{3}15_{1},20_{3})$ as in Figure 7; include the edges of $D$ among the lines of configuration, exclude the vertices of $D$ from the points of the configuration. Thus we get a configuration with $12 \times 15 + 30 \times 3 = 270$ 3-valent points and $12 \times 20 + 30 = 270$ 3-valent lines, hence a $(270_{3}$ configuration. The same result we get using Proposition 1, using the fact that $x = 0, y = 3, z = 20$ and dodecahedral formula for points $20x + 30y + 12z$. Note that on each edge the position of three points remains the same if the edge is reflected about its middle point.

\vspace{-1mm}

\begin{center}
\definecolor{uuuuuu}{rgb}{0.26666666666666666,0.26666666666666666,0.26666666666666666}
\definecolor{zzttqq}{rgb}{0.6,0.2,0.0}
\definecolor{ffffff}{rgb}{1.0,1.0,1.0}
\definecolor{qqqqff}{rgb}{0.0,0.0,1.0}
\definecolor{ffqqqq}{rgb}{1.0,0.0,0.0}

\end{center}

\begin{center}

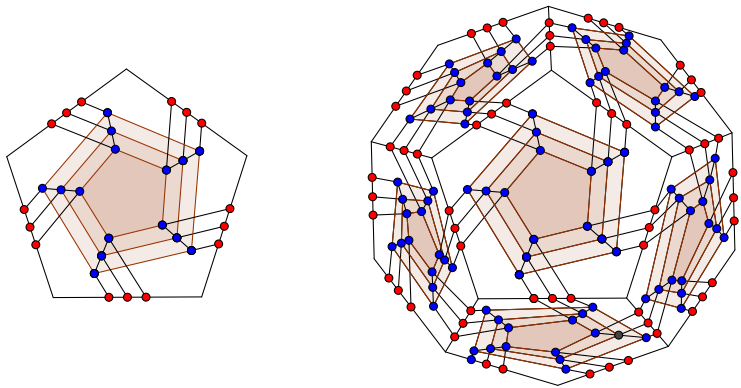
Figure 7. Placing the copies of $(15_{3}15_{1},20_{3})$ to each face raises the valences of 1-valent points to 2; adding the edges of $D$ raises them to 3.

\end{center}

\begin{prob}\label{examplepeti} To construct examples of Platonic configurations $Z \in C(p_{3}, n_{4})$ and $Z \in C(n_{4})$ with the full symmetry of the cube.
\end{prob}
\textbf{Solution.} We have already mentioned the Grey configuration $(27_{3})$,
% with $3 \times 3 \times 3$ points and three sets of 9 parallel 3-valent lines,
 represented as a cubic $C_{3}$ configuration with the full symmetry group of the cube (see \cite{PisaServ}, pp. 252--254). Similarly, subdivide each face of two "concentric cubes" into a square grid with $4 \times 4 \times 4$ points. You have $2 \times (8 + 12 \times 2) = 64$ 3-valent points (in the vertices and on the edges) and $2 \times 6 \times 4 = 48$ 2-valent points (in the interiors of faces). You have $2 \times 3 \times 12 = 72$ 4-valent lines parallel with one of the faces. Add to them $24/2 = 12$ "antipodal lines" (each of them connecting 2 "antipodal pairs" of 2-valent points) and you get a configuration with $64 + 48 = 112$ 3-valent points and $72 + 12$ = 84  4-valent points, thus a $(112_{3}, 84_{4})$ configuration $B$. To obtain from it a $(n_{4})$ configuration with the full symmetry group of the cube, place $3 \times 4$ copies of $B$ around the origin of the Cartesian coordinate system in $E^{3}$ in such a way that they form a "centrally symmetric" pattern and that the coordinate lines go through the centers of two opposite faces of each cube. Add $3 \times 112 = 336$ "antipodal lines". 
Now you have a configuration with $12 \times 112 = 1344$ 4-valent points and $12 \times 84 + 336 = 1008 + 336 = 1344$ 4-valent lines, thus a $(1344_{4})$ configuration with the full symmetry group of the cube. Similar iterative constructions produce $C(n_{k})$ configurations for all even $k \geq 2$.

\begin{prob} \label{C3R} To construct a Platonic configuration $Z \in C_{3,R}$ of points and lines with the symmetry group of the rotations of the cube $C$.
\end{prob} 

\textbf{Solution.} Place on each face of the cube a copy of $A = (4_{3}1_{2}4_{1},6_{3})$ and identify 1-valent points with the 8 vertices of the cube. We obtain a configuration $B = ((8+24)_{3}6_{2},36_{3})$. Adding the central point of the cube and 3 lines connecting the antipodal 2-valent points of $B$ we get $X = (39_{3})$.

\begin{center}
\definecolor{ffffqq}{rgb}{1.0,1.0,0.0}
\definecolor{uuuuuu}{rgb}{0.26666666666666666,0.26666666666666666,0.26666666666666666}
\definecolor{zzttqq}{rgb}{0.6,0.2,0.0}
\definecolor{qqqqff}{rgb}{0.0,0.0,1.0}
\begin{tikzpicture}[line cap=round,line join=round,>=triangle 45,x=0.4cm,y=0.4cm]
\clip(4.809190120849255,0.06511515329807492) rectangle (15.388107140440338,5.324913148085403);
\fill[color=zzttqq,fill=zzttqq,fill opacity=0.1] (-1.84,2.34) -- (-0.56,1.48) -- (0.2999999999999996,2.76) -- (-0.98,3.6199999999999997) -- cycle;
\fill[color=zzttqq,fill=zzttqq,fill opacity=0.1] (11.725584837749954,1.855605775930031) -- (12.463397851523762,1.489037328276599) -- (12.796766655931773,2.196514866839123) -- (12.156766655931772,2.6165148668391227) -- cycle;
\fill[color=zzttqq,fill=zzttqq,fill opacity=0.1] (14.10828957881953,2.328368019834027) -- (14.436766655931773,2.236514866839122) -- (14.576766655931774,2.7565148668391215) -- (14.256766655931774,2.9965148668391235) -- cycle;
\fill[color=zzttqq,fill=zzttqq,fill opacity=0.1] (12.436766655931773,3.9765148668391217) -- (12.977238120125843,3.75357398133016) -- (13.600187943568319,4.099657216575977) -- (13.136766655931773,4.336514866839121) -- cycle;
\fill[color=zzttqq,fill=zzttqq,fill opacity=0.1] (6.410857931292018,2.457730030626145) -- (7.690857931292017,1.597730030626145) -- (8.550857931292017,2.877730030626145) -- (7.270857931292017,3.7377300306261443) -- cycle;
\draw [color=zzttqq] (-1.84,2.34)-- (-0.56,1.48);
\draw [color=zzttqq] (-0.56,1.48)-- (0.2999999999999996,2.76);
\draw [color=zzttqq] (0.2999999999999996,2.76)-- (-0.98,3.6199999999999997);
\draw [color=zzttqq] (-0.98,3.6199999999999997)-- (-1.84,2.34);
\draw (-0.56,1.48)-- (-0.98,3.6199999999999997);
\draw (-1.84,2.34)-- (0.2999999999999996,2.76);
\draw (-1.84,2.34)-- (-2.62,1.16);
\draw (-0.56,1.48)-- (0.7,0.68);
\draw (0.2999999999999996,2.76)-- (1.12,3.88);
\draw (-0.98,3.6199999999999997)-- (-2.12,4.28);
\draw [dotted] (11.136766655931773,0.8165148668391218)-- (13.776766655931773,0.8365148668391214);
\draw [dotted] (13.656766655931772,3.436514866839123)-- (13.776766655931773,0.8365148668391214);
\draw [dotted] (13.656766655931772,3.436514866839123)-- (11.116766655931773,3.41651486683912);
\draw [dotted] (11.116766655931773,3.41651486683912)-- (11.136766655931773,0.8165148668391218);
\draw [dotted] (13.656766655931772,3.436514866839123)-- (14.834444193620618,4.889892261298481);
\draw [dotted] (14.834444193620618,4.889892261298481)-- (14.996766655931772,1.9965148668391217);
\draw [dotted] (13.776766655931773,0.8365148668391214)-- (14.996766655931772,1.9965148668391217);
\draw [dotted] (11.116766655931773,3.41651486683912)-- (12.22431887289908,4.909665937970614);
\draw [dotted] (12.22431887289908,4.909665937970614)-- (14.834444193620618,4.889892261298481);
\draw (11.136766655931773,0.8165148668391218)-- (12.156766655931772,2.6165148668391227);
\draw (11.116766655931773,3.41651486683912)-- (12.156766655931772,2.6165148668391227);
\draw (12.156766655931772,2.6165148668391227)-- (12.796766655931773,2.196514866839123);
\draw (11.725584837749954,1.855605775930031)-- (13.776766655931773,0.8365148668391214);
\draw (12.796766655931773,2.196514866839123)-- (13.656766655931772,3.436514866839123);
\draw (12.796766655931773,2.196514866839123)-- (12.463397851523762,1.489037328276599);
\draw (13.776766655931773,0.8365148668391214)-- (14.256766655931774,2.9965148668391235);
\draw (13.656766655931772,3.436514866839123)-- (14.256766655931774,2.9965148668391235);
\draw (14.256766655931774,2.9965148668391235)-- (14.576766655931774,2.7565148668391215);
\draw (14.576766655931774,2.7565148668391215)-- (14.834444193620618,4.889892261298481);
\draw (14.576766655931774,2.7565148668391215)-- (14.436766655931773,2.236514866839122);
\draw (14.436766655931773,2.236514866839122)-- (14.10828957881953,2.328368019834027);
\draw (14.436766655931773,2.236514866839122)-- (14.996766655931772,1.9965148668391217);
\draw (13.656766655931772,3.436514866839123)-- (12.436766655931773,3.9765148668391217);
\draw (11.116766655931773,3.41651486683912)-- (12.436766655931773,3.9765148668391217);
\draw (12.436766655931773,3.9765148668391217)-- (13.136766655931773,4.336514866839121);
\draw (12.22431887289908,4.909665937970614)-- (13.136766655931773,4.336514866839121);
\draw (13.136766655931773,4.336514866839121)-- (13.600187943568319,4.099657216575977);
\draw (14.834444193620618,4.889892261298481)-- (13.600187943568319,4.099657216575977);
\draw (13.600187943568319,4.099657216575977)-- (12.977238120125843,3.75357398133016);
\draw (12.977238120125843,3.75357398133016)-- (13.136766655931773,4.336514866839121);
\draw (12.436766655931773,3.9765148668391217)-- (13.600187943568319,4.099657216575977);
\draw (14.256766655931774,2.9965148668391235)-- (14.436766655931773,2.236514866839122);
\draw (14.576766655931774,2.7565148668391215)-- (14.10828957881953,2.328368019834027);
\draw (12.156766655931772,2.6165148668391227)-- (12.463397851523762,1.489037328276599);
\draw (11.725584837749954,1.855605775930031)-- (12.796766655931773,2.196514866839123);
\draw [color=zzttqq] (11.725584837749954,1.855605775930031)-- (12.463397851523762,1.489037328276599);
\draw [color=zzttqq] (12.463397851523762,1.489037328276599)-- (12.796766655931773,2.196514866839123);
\draw [color=zzttqq] (12.796766655931773,2.196514866839123)-- (12.156766655931772,2.6165148668391227);
\draw [color=zzttqq] (12.156766655931772,2.6165148668391227)-- (11.725584837749954,1.855605775930031);
\draw [color=zzttqq] (14.10828957881953,2.328368019834027)-- (14.436766655931773,2.236514866839122);
\draw [color=zzttqq] (14.436766655931773,2.236514866839122)-- (14.576766655931774,2.7565148668391215);
\draw [color=zzttqq] (14.576766655931774,2.7565148668391215)-- (14.256766655931774,2.9965148668391235);
\draw [color=zzttqq] (14.256766655931774,2.9965148668391235)-- (14.10828957881953,2.328368019834027);
\draw [color=zzttqq] (12.436766655931773,3.9765148668391217)-- (12.977238120125843,3.75357398133016);
\draw [color=zzttqq] (12.977238120125843,3.75357398133016)-- (13.600187943568319,4.099657216575977);
\draw [color=zzttqq] (13.600187943568319,4.099657216575977)-- (13.136766655931773,4.336514866839121);
\draw [color=zzttqq] (13.136766655931773,4.336514866839121)-- (12.436766655931773,3.9765148668391217);
\draw [color=zzttqq] (6.410857931292018,2.457730030626145)-- (7.690857931292017,1.597730030626145);
\draw [color=zzttqq] (7.690857931292017,1.597730030626145)-- (8.550857931292017,2.877730030626145);
\draw [color=zzttqq] (8.550857931292017,2.877730030626145)-- (7.270857931292017,3.7377300306261443);
\draw [color=zzttqq] (7.270857931292017,3.7377300306261443)-- (6.410857931292018,2.457730030626145);
\draw (7.690857931292017,1.597730030626145)-- (7.270857931292017,3.7377300306261443);
\draw (6.410857931292018,2.457730030626145)-- (8.550857931292017,2.877730030626145);
\draw (6.410857931292018,2.457730030626145)-- (5.630857931292018,1.277730030626145);
\draw (7.690857931292017,1.597730030626145)-- (8.950857931292015,0.797730030626145);
\draw (8.550857931292017,2.877730030626145)-- (9.370857931292017,3.9977300306261507);
\draw (7.270857931292017,3.7377300306261443)-- (6.1308579312920175,4.3977300306261515);
\begin{scriptsize}
\draw [fill=qqqqff] (-1.84,2.34) circle (1.5pt);
\draw [fill=qqqqff] (-0.56,1.48) circle (1.5pt);
\draw [fill=uuuuuu] (0.2999999999999996,2.76) circle (1.5pt);
\draw [fill=uuuuuu] (-0.98,3.6199999999999997) circle (1.5pt);
\draw [fill=qqqqff] (-2.62,1.16) circle (1.5pt);
\draw [fill=qqqqff] (0.7,0.68) circle (1.5pt);
\draw [fill=qqqqff] (1.12,3.88) circle (1.5pt);
\draw [fill=qqqqff] (-2.12,4.28) circle (1.5pt);
\draw [fill=uuuuuu] (-0.7700000000000002,2.5500000000000003) circle (1.5pt);
\draw [fill=ffffqq] (11.136766655931773,0.8165148668391218) circle (1.5pt);
\draw [fill=ffffqq] (13.776766655931773,0.8365148668391214) circle (1.5pt);
\draw [fill=ffffqq] (13.656766655931772,3.436514866839123) circle (1.5pt);
\draw [fill=ffffqq] (11.116766655931773,3.41651486683912) circle (1.5pt);
\draw [fill=ffffqq] (14.834444193620618,4.889892261298481) circle (1.5pt);
\draw [fill=ffffqq] (14.996766655931772,1.9965148668391217) circle (1.5pt);
\draw [fill=ffffqq] (12.22431887289908,4.909665937970614) circle (1.5pt);
\draw [fill=qqqqff] (12.156766655931772,2.6165148668391227) circle (1.5pt);
\draw [fill=qqqqff] (12.796766655931773,2.196514866839123) circle (1.5pt);
\draw [fill=qqqqff] (11.725584837749954,1.855605775930031) circle (1.5pt);
\draw [fill=qqqqff] (12.463397851523762,1.489037328276599) circle (1.5pt);
\draw [fill=qqqqff] (14.256766655931774,2.9965148668391235) circle (1.5pt);
\draw [fill=qqqqff] (14.576766655931774,2.7565148668391215) circle (1.5pt);
\draw [fill=qqqqff] (14.436766655931773,2.236514866839122) circle (1.5pt);
\draw [fill=qqqqff] (14.10828957881953,2.328368019834027) circle (1.5pt);
\draw [fill=qqqqff] (12.436766655931773,3.9765148668391217) circle (1.5pt);
\draw [fill=qqqqff] (13.136766655931773,4.336514866839121) circle (1.5pt);
\draw [fill=qqqqff] (13.600187943568319,4.099657216575977) circle (1.5pt);
\draw [fill=qqqqff] (12.977238120125843,3.75357398133016) circle (1.5pt);
\draw [fill=qqqqff] (12.312873334353418,2.042513342955659) circle (1.5pt);
\draw [fill=qqqqff] (14.360434478515765,2.558806282595593) circle (1.5pt);
\draw [fill=qqqqff] (13.056190536950366,4.042077800062314) circle (1.5pt);
\draw [fill=qqqqff] (6.410857931292018,2.457730030626145) circle (1.5pt);
\draw [fill=qqqqff] (7.690857931292017,1.597730030626145) circle (1.5pt);
\draw [fill=uuuuuu] (8.550857931292017,2.877730030626145) circle (1.5pt);
\draw [fill=uuuuuu] (7.270857931292017,3.7377300306261443) circle (1.5pt);
\draw [fill=ffffqq] (5.630857931292018,1.277730030626145) circle (1.5pt);
\draw [fill=ffffqq] (8.950857931292015,0.797730030626145) circle (1.5pt);
\draw [fill=ffffqq] (9.370857931292017,3.9977300306261507) circle (1.5pt);
\draw [fill=ffffqq] (6.1308579312920175,4.3977300306261515) circle (1.5pt);
\draw [fill=uuuuuu] (7.480857931292016,2.6677300306261453) circle (1.5pt);
\end{scriptsize}
\end{tikzpicture}
\end{center}

\begin{center}
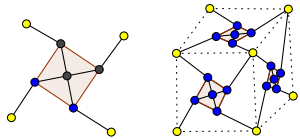
Figure 8. Platonic configuration $(39_{3})$ with rotational symmetry of the cube.
\end{center}

\begin{prob}\label{TOR3R}
For $P \in \lbrace T, O, I\rbrace$ construct $Z \in P_{3,R}$. Construct also $Z \in P_{3,R}$ for $P \in \lbrace C, D\rbrace$.
\end{prob} 

\textbf{Solution.} Place the copies of the same configuration $A = (9_{3}9_{1},12_{3}) \in Cyc_{3}$ (Figure  9) on each of the triangular faces of $P(v,e,f,d,m) \in \lbrace T, O, I\rbrace$, so that the 1-valent point lie on the edges  of $P$, which are included among the lines of configuration. The points on each edge must be equidistant and the middle one must be the center of the edge. We get 3-valent Platonic configurations with $p = 9f + 3e$ points and $n = 12f + e$ lines. Since $3f = 2e$, we have $p = 6e + 3e = 9e$ and $n = 8e + e = 9e$. Thus we get configurations $(36_{3}) \in T_{3,R}$, $(108_{3}) \in O_{3,R}$, and $(270_{3}) \in I_{3,R}$.

Triangulating the faces of the cube and dodecahedron 
and placing the same $A$ (as in Figure 9) into each of these triangles we obtain 3-valent configurations
with rotational symmetry of the cube and the dodecahedron. In the case of the cube (in which $4f = 2e$) the number of points is $p = 48 \times f + 3e = 24e + 3e = 27 \times 12 = 324$  and the number of lines is $52 \times f + e  = 26 e + e = 324$. Hence we have $(324_{3}) \in C_{3,R}$. In the case of the dodecahedron (in which $5f = 2e$) we have $p = 60f + 3e = 24 e + 3e = 27e = 810$ points and $n = 65f + e = 26e + e = 27e = 810$ lines. Hence we have $(810_{3}) \in C_{3,R}$.

\begin{center}
\definecolor{zzttqq}{rgb}{0.6,0.2,0.0}
\definecolor{qqqqff}{rgb}{0.0,0.0,1.0}
\definecolor{ffqqqq}{rgb}{1.0,0.0,0.0}
\begin{tikzpicture}[line cap=round,line join=round,>=triangle 45,x=0.3cm,y=0.35cm]
\clip(-23.436314079501244,-5.544309000454292) rectangle (15.837006733487597,10.415367186497363);
\fill[color=zzttqq,fill=zzttqq,fill opacity=0.1] (-14.177968236800782,-1.2002011018909178) -- (-12.4769342085851,0.2804748091929526) -- (-14.613436167373113,1.0145991202894618) -- cycle;
\fill[color=zzttqq,fill=zzttqq,fill opacity=0.1] (-14.394341323576095,-1.7396685291510645) -- (-11.91076061087006,0.36579135028584364) -- (-14.975259535567291,1.4632562378638136) -- cycle;
\fill[color=zzttqq,fill=zzttqq,fill opacity=0.1] (-14.619382815472687,-2.2817642827022655) -- (-11.33282755405132,0.4457285053639749) -- (-15.334826266572511,1.9298756611031396) -- cycle;
\fill[color=zzttqq,fill=zzttqq,fill opacity=0.1] (6.275563091526542,-1.4176563607455344) -- (7.976597119742227,0.06301955033833531) -- (5.840095160954209,0.7971438614348452) -- cycle;
\fill[color=zzttqq,fill=zzttqq,fill opacity=0.1] (6.059190004751227,-1.957123788005683) -- (8.542770717457264,0.14833609143122675) -- (5.478271792760035,1.2458009790091973) -- cycle;
\fill[color=zzttqq,fill=zzttqq,fill opacity=0.1] (5.834148512854638,-2.4992195415568785) -- (9.120703774276,0.22827324650936012) -- (5.118705061754811,1.7124204022485245) -- cycle;
\draw (-18.090819028834368,-2.4103164447867096)-- (-9.480819028834368,-2.500316444786712);
\draw (-19.499838266720836,0.09017646407463209)-- (-8.019838266720836,-0.029823535925373346);
\draw (-18.0388575046073,2.56066937293597)-- (-9.428857504607302,2.4706693729359666);
\draw (-13.707876742493768,5.001162281797307)-- (-18.090819028834368,-2.4103164447867096);
\draw (-10.83787674249377,4.971162281797306)-- (-16.681799790947903,-4.910809353648051);
\draw (-9.428857504607302,2.4706693729359666)-- (-13.811799790947902,-4.940809353648051);
\draw (-13.811799790947902,-4.940809353648051)-- (-18.0388575046073,2.56066937293597);
\draw (-10.941799790947902,-4.97080935364805)-- (-16.577876742493768,5.031162281797308);
\draw (-9.480819028834368,-2.500316444786712)-- (-13.707876742493768,5.001162281797307);
\draw (-17.386309409891133,-3.66056289921738)-- (-10.211309409891136,-3.7355628992173813);
\draw (-18.795328647777602,-1.1600699903560387)-- (-8.750328647777602,-1.2650699903560427);
\draw (-18.769347885664068,1.325422918505301)-- (-8.72434788566407,1.2204229185052966);
\draw (-17.308367123550532,3.795915827366639)-- (-10.133367123550535,3.7209158273666363);
\draw (-15.142876742493769,5.016162281797308)-- (-18.795328647777602,-1.1600699903560387);
\draw (-12.272876742493768,4.986162281797307)-- (-17.386309409891133,-3.66056289921738);
\draw (-10.133367123550535,3.7209158273666363)-- (-15.246799790947902,-4.9258093536480505);
\draw (-12.376799790947903,-4.955809353648051)-- (-8.72434788566407,1.2204229185052966);
\draw (-15.142876742493769,5.016162281797308)-- (-10.211309409891136,-3.7355628992173813);
\draw (-8.750328647777602,-1.2650699903560427)-- (-12.272876742493768,4.986162281797307);
\draw (-18.769347885664068,1.325422918505301)-- (-15.246799790947902,-4.9258093536480505);
\draw (-17.308367123550532,3.795915827366639)-- (-12.376799790947903,-4.955809353648051);
\draw (-17.034054600419516,-4.285686126432715)-- (-10.57655460041952,-4.353186126432716);
\draw (-17.73856421936275,-3.035439672002045)-- (-9.846064219362752,-3.1179396720020467);
\draw (-18.443073838305985,-1.7851932175713743)-- (-9.115573838305984,-1.8826932175713775);
\draw (-19.14758345724922,-0.5349467631407033)-- (-8.38508345724922,-0.647446763140708);
\draw (-19.134593076192452,0.7077996912899666)-- (-8.372093076192453,0.5952996912899616);
\draw (-18.404102695135684,1.9430461457206354)-- (-9.076602695135687,1.8455461457206317);
\draw (-17.673612314078916,3.1782926001513045)-- (-9.781112314078918,3.0957926001513014);
\draw (-16.943121933022148,4.413539054581974)-- (-10.485621933022152,4.346039054581971);
\draw (-19.134593076192452,0.7077996912899666)-- (-15.964299790947901,-4.91830935364805);
\draw (-18.404102695135684,1.9430461457206354)-- (-14.529299790947903,-4.933309353648051);
\draw (-17.673612314078916,3.1782926001513045)-- (-13.094299790947902,-4.948309353648051);
\draw (-16.943121933022148,4.413539054581974)-- (-11.659299790947902,-4.96330935364805);
\draw (-15.860376742493768,5.023662281797308)-- (-10.57655460041952,-4.353186126432716);
\draw (-14.425376742493768,5.008662281797307)-- (-9.846064219362752,-3.1179396720020467);
\draw (-12.990376742493769,4.993662281797307)-- (-9.115573838305984,-1.8826932175713775);
\draw (-11.555376742493769,4.978662281797306)-- (-8.38508345724922,-0.647446763140708);
\draw (-11.6592997909479,-4.963309353648049)-- (-8.372093076192453,0.5952996912899616);
\draw (-13.094299790947904,-4.94830935364805)-- (-9.076602695135687,1.8455461457206317);
\draw (-15.9642997909479,-4.918309353648051)-- (-10.485621933022152,4.346039054581971);
\draw (-14.529299790947904,-4.93330935364805)-- (-9.781112314078918,3.0957926001513014);
\draw (-17.034054600419516,-4.285686126432715)-- (-11.555376742493769,4.978662281797306);
\draw (-17.73856421936275,-3.035439672002045)-- (-12.990376742493769,4.993662281797307);
\draw (-18.443073838305985,-1.7851932175713743)-- (-14.425376742493768,5.008662281797307);
\draw (-19.14758345724922,-0.5349467631407033)-- (-15.860376742493768,5.023662281797308);
\draw (-18.038989708522067,2.5604458184290104)-- (-10.889838266720835,1.764640746280112E-4);
\draw (-9.428537057043567,2.4701006959163503)-- (-15.22081902883437,-2.44031644478671);
\draw (-13.813750365353826,-4.940788964368546)-- (-15.168857504607303,2.530669372935967);
\draw (-10.889838266720835,1.764640746280112E-4)-- (-9.102583457249217,-0.6399467631407084);
\draw (-17.858435470667843,2.8657598758423597)-- (-8.385083457249221,-0.6474467631407081);
\draw (-18.220249971258177,2.2539378852923093)-- (-9.820083457249218,-0.6324467631407076);
\draw (-14.168316200447062,-4.937082701283948)-- (-15.860376742493768,5.023662281797309);
\draw (-13.441739566499237,-4.94467757899072)-- (-15.155867123550536,3.773415827366639);
\draw (-15.22081902883437,-2.44031644478671)-- (-16.668809409891132,-3.6680628992173783);
\draw (-9.258059791675077,2.1675660095799185)-- (-17.034054600419516,-4.285686126432716);
\draw (-9.603907590234778,2.781319139835985)-- (-16.30356421936275,-3.0504396720020455);
\draw [color=zzttqq] (-14.177968236800782,-1.2002011018909178)-- (-12.4769342085851,0.2804748091929526);
\draw [color=zzttqq] (-12.4769342085851,0.2804748091929526)-- (-14.613436167373113,1.0145991202894618);
\draw [color=zzttqq] (-14.613436167373113,1.0145991202894618)-- (-14.177968236800782,-1.2002011018909178);
\draw (-15.334826266572511,1.9298756611031396)-- (-14.613436167373113,1.0145991202894618);
\draw (-12.4769342085851,0.2804748091929526)-- (-11.33282755405132,0.4457285053639749);
\draw (-14.619382815472687,-2.2817642827022655)-- (-14.177968236800782,-1.2002011018909178);
\draw [color=zzttqq] (-14.394341323576095,-1.7396685291510645)-- (-11.91076061087006,0.36579135028584364);
\draw [color=zzttqq] (-11.91076061087006,0.36579135028584364)-- (-14.975259535567291,1.4632562378638136);
\draw [color=zzttqq] (-14.975259535567291,1.4632562378638136)-- (-14.394341323576095,-1.7396685291510645);
\draw [color=zzttqq] (-14.619382815472687,-2.2817642827022655)-- (-11.33282755405132,0.4457285053639749);
\draw [color=zzttqq] (-11.33282755405132,0.4457285053639749)-- (-15.334826266572511,1.9298756611031396);
\draw [color=zzttqq] (-15.334826266572511,1.9298756611031396)-- (-14.619382815472687,-2.2817642827022655);
\draw [dotted] (-22.421799790947844,-4.850809353648041)-- (-5.201799790947904,-5.030809353648044);
\draw [dotted] (-13.655915218266731,9.972148099520036)-- (-5.201799790947904,-5.030809353648044);
\draw [dotted] (-13.655915218266731,9.972148099520036)-- (-22.421799790947844,-4.850809353648041);
\draw (-15.508121933022155,4.398539054581973)-- (-13.813750365353826,-4.940788964368546);
\draw (-19.49983826672082,0.09017646407463187)-- (-16.681799790947895,-4.91080935364804);
\draw (-10.941799790947908,-4.9708093536480416)-- (-8.019838266720834,-0.029823535925373373);
\draw (-16.577876742493775,5.031162281797323)-- (-10.83787674249377,4.971162281797307);
\draw [color=zzttqq] (6.275563091526542,-1.4176563607455344)-- (7.976597119742227,0.06301955033833531);
\draw [color=zzttqq] (7.976597119742227,0.06301955033833531)-- (5.840095160954209,0.7971438614348452);
\draw [color=zzttqq] (5.840095160954209,0.7971438614348452)-- (6.275563091526542,-1.4176563607455344);
\draw (5.118705061754811,1.7124204022485245)-- (5.840095160954209,0.7971438614348452);
\draw (7.976597119742227,0.06301955033833531)-- (9.120703774276,0.22827324650936012);
\draw (5.834148512854638,-2.4992195415568785)-- (6.275563091526542,-1.4176563607455344);
\draw [color=zzttqq] (6.059190004751227,-1.957123788005683)-- (8.542770717457264,0.14833609143122675);
\draw [color=zzttqq] (8.542770717457264,0.14833609143122675)-- (5.478271792760035,1.2458009790091973);
\draw [color=zzttqq] (5.478271792760035,1.2458009790091973)-- (6.059190004751227,-1.957123788005683);
\draw [color=zzttqq] (5.834148512854638,-2.4992195415568785)-- (9.120703774276,0.22827324650936012);
\draw [color=zzttqq] (9.120703774276,0.22827324650936012)-- (5.118705061754811,1.7124204022485245);
\draw [color=zzttqq] (5.118705061754811,1.7124204022485245)-- (5.834148512854638,-2.4992195415568785);
\draw [dotted] (-1.9682684626205864,-5.0682646125026825)-- (15.251731537379444,-5.2482646125026875);
\draw [dotted] (6.797616110060623,9.754692840665374)-- (15.251731537379444,-5.2482646125026875);
\draw [dotted] (6.797616110060623,9.754692840665374)-- (-1.9682684626205864,-5.0682646125026825);
\draw (5.11870506175481,1.7124204022485248)-- (2.5950958576594827,2.6483046169877436);
\draw (5.478271792760035,1.2458009790091973)-- (2.4145416198052563,2.3429905595743943);
\draw (5.840095160954209,0.7971438614348452)-- (2.2332813570691465,2.036482626437693);
\draw (6.059190004751227,-1.957123788005683)-- (6.639780962973498,-5.1582442232231625);
\draw (11.024994271283758,2.252645437061734)-- (8.542770717457264,0.14833609143122675);
\draw (11.195471536652246,1.9501107507253024)-- (9.120703774276,0.22827324650936012);
\draw (10.849623738092546,2.5638638809813687)-- (7.976597119742227,0.06301955033833531);
\begin{scriptsize}
\draw [fill=ffqqqq] (-9.428537057043567,2.4701006959163503) circle (1.5pt);
\draw [fill=ffqqqq] (-18.038989708522067,2.5604458184290104) circle (1.5pt);
\draw [fill=ffqqqq] (-13.813750365353826,-4.940788964368546) circle (1.5pt);
\draw [fill=ffqqqq] (-17.858435470667843,2.8657598758423597) circle (1.5pt);
\draw [fill=ffqqqq] (-18.220249971258177,2.2539378852923093) circle (1.5pt);
\draw [fill=ffqqqq] (-14.168316200447062,-4.937082701283948) circle (1.5pt);
\draw [fill=ffqqqq] (-13.441739566499237,-4.94467757899072) circle (1.5pt);
\draw [fill=ffqqqq] (-9.258059791675077,2.1675660095799185) circle (1.5pt);
\draw [fill=ffqqqq] (-9.603907590234778,2.781319139835985) circle (1.5pt);
\draw [fill=qqqqff] (-12.4769342085851,0.2804748091929526) circle (1.5pt);
\draw [fill=qqqqff] (-14.613436167373113,1.0145991202894618) circle (1.5pt);
\draw [fill=qqqqff] (-14.177968236800782,-1.2002011018909178) circle (1.5pt);
\draw [fill=qqqqff] (-15.334826266572511,1.9298756611031396) circle (1.5pt);
\draw [fill=qqqqff] (-11.91076061087006,0.36579135028584364) circle (1.5pt);
\draw [fill=qqqqff] (-11.33282755405132,0.4457285053639749) circle (1.5pt);
\draw [fill=qqqqff] (-14.394341323576095,-1.7396685291510645) circle (1.5pt);
\draw [fill=qqqqff] (-14.619382815472687,-2.2817642827022655) circle (1.5pt);
\draw [fill=qqqqff] (-14.975259535567291,1.4632562378638136) circle (1.5pt);
\draw [fill=ffqqqq] (11.024994271283758,2.252645437061734) circle (1.5pt);
\draw [fill=ffqqqq] (2.4145416198052563,2.3429905595743943) circle (1.5pt);
\draw [fill=ffqqqq] (6.639780962973498,-5.1582442232231625) circle (1.5pt);
\draw [fill=ffqqqq] (2.5950958576594827,2.6483046169877436) circle (1.5pt);
\draw [fill=ffqqqq] (2.2332813570691465,2.036482626437693) circle (1.5pt);
\draw [fill=ffqqqq] (6.285215127880263,-5.154537960138564) circle (1.5pt);
\draw [fill=ffqqqq] (7.0117917618280865,-5.162132837845336) circle (1.5pt);
\draw [fill=ffqqqq] (11.195471536652246,1.9501107507253024) circle (1.5pt);
\draw [fill=ffqqqq] (10.849623738092546,2.5638638809813687) circle (1.5pt);
\draw [fill=qqqqff] (7.976597119742227,0.06301955033833531) circle (1.5pt);
\draw [fill=qqqqff] (5.840095160954209,0.7971438614348452) circle (1.5pt);
\draw [fill=qqqqff] (6.275563091526542,-1.4176563607455344) circle (1.5pt);
\draw [fill=qqqqff] (8.542770717457264,0.14833609143122675) circle (1.5pt);
\draw [fill=qqqqff] (9.120703774276,0.22827324650936012) circle (1.5pt);
\draw [fill=qqqqff] (6.059190004751227,-1.957123788005683) circle (1.5pt);
\draw [fill=qqqqff] (5.834148512854638,-2.4992195415568785) circle (1.5pt);
\draw [fill=qqqqff] (5.478271792760035,1.2458009790091973) circle (1.5pt);
\draw [fill=qqqqff] (5.11870506175481,1.7124204022485248) circle (1.5pt);
\end{scriptsize}
\end{tikzpicture}
\end{center}

\begin{center}
Figure 9. A configuration $A = (9_{3}9_{1},12_{3}) \in Cyc_{3}$ 
\end{center}

\begin{remark}
Similarly, we easily find examples of Platonic configurations $P_{R}(p_{3}, n_{k})$ with rotational symmetry of the chosen Platonic solid $P$ for any $k \geq 3$ (we just add $k - 3$ concentric triangles and $3(k -3)$ lines to the pattern of $A$).
\end{remark}

Examples of Platonic $(n_{4})$ configurations may be found by the three-step pattern $A \rightarrow B \rightarrow Z$ from Algorithm 1: starting with some configuration $A$ we construct a Platonic configuration $B = (a_{4}b_{3},c_{4})$; if it is not balanced, we use the method of connecting antipodal points of valence 3 in concentric copies of $B$ to increase their valence to 4. In many cases we can take for $A$ the Pappus configuration $A = (9_{3})$, or some other configuration obtained from it by some small modification.

\begin{prob}\label{O4}
To construct a Platonic configuration $Z \in O_{4}$ with the full symmetry group of the octahedron.
\end{prob}

\textbf{Solution.} Place the copy of $A = (3_{1}6_{2}7_{3}, 9_{4})$ on each of the  faces of the octahedron in such a way, that their 1-valent points coincide with the vertices of the faces (see Figure 10) to obtain $B = (30_{4}56_{3},72_{4})$ with $2e + v = 24 + 6 = 30$ 4-valent points and $7 \times f = 56$ 3-valent points and $9f = 72$  4-valent lines.
Now take two concentric copies of $B$ and connect antipodal points of valence 3 with 4-valent lines. Now all the points and lines have valence 4. The number of lines is $72 \times 2 + \frac{56}{2} = 144 + 28 = 172$. The number of points is $2 \times 86 = 172$. So we have a $(172_{4})$ configuration.\\

\begin{center}
\definecolor{ffffqq}{rgb}{1.0,1.0,0.0}
\definecolor{ffqqqq}{rgb}{1.0,0.0,0.0}
\definecolor{xdxdff}{rgb}{0.49019607843137253,0.49019607843137253,1.0}
\definecolor{qqqqff}{rgb}{0.0,0.0,1.0}
\definecolor{zzttqq}{rgb}{0.6,0.2,0.0}
\begin{tikzpicture}[line cap=round,line join=round,>=triangle 45,x=0.25cm,y=0.25cm]
\clip(9.678420140988488,0.6121017704003622) rectangle (62.56569546196679,16.07963656276689);
\fill[color=zzttqq,fill=zzttqq,fill opacity=0.1] (16.252358518496393,2.0222919800730743) -- (21.992358518496392,1.9622919800730743) -- (24.914320042723457,6.903277797795752) -- (22.096281566950523,11.90426361551843) -- (16.356281566950525,11.96426361551843) -- (13.43432004272346,7.023277797795755) -- cycle;
\fill[color=zzttqq,fill=zzttqq,fill opacity=0.1] (39.54065020783629,-11.209691934324564) -- (45.28065020783629,-11.269691934324564) -- (48.20261173206335,-6.328706116601893) -- (45.38457325629042,-1.3277202988792176) -- (39.64457325629043,-1.2677202988792153) -- (36.72261173206336,-6.208706116601885) -- cycle;
\draw [color=zzttqq] (16.252358518496393,2.0222919800730743)-- (21.992358518496392,1.9622919800730743);
\draw [color=zzttqq] (21.992358518496392,1.9622919800730743)-- (24.914320042723457,6.903277797795752);
\draw [color=zzttqq] (24.914320042723457,6.903277797795752)-- (22.096281566950523,11.90426361551843);
\draw [color=zzttqq] (22.096281566950523,11.90426361551843)-- (16.356281566950525,11.96426361551843);
\draw [color=zzttqq] (16.356281566950525,11.96426361551843)-- (13.43432004272346,7.023277797795755);
\draw [color=zzttqq] (13.43432004272346,7.023277797795755)-- (16.252358518496393,2.0222919800730743);
\draw (14.843339280609927,4.522784888934415)-- (23.453339280609924,4.432784888934413);
\draw (13.43432004272346,7.023277797795755)-- (24.914320042723457,6.903277797795752);
\draw (14.895300804836992,9.493770706657092)-- (23.505300804836992,9.403770706657092);
\draw (19.226281566950526,11.934263615518429)-- (14.843339280609927,4.522784888934415);
\draw (22.096281566950523,11.90426361551843)-- (16.252358518496393,2.0222919800730743);
\draw (23.505300804836992,9.403770706657092)-- (19.12235851849639,1.9922919800730743);
\draw (19.12235851849639,1.9922919800730743)-- (14.895300804836992,9.493770706657092);
\draw (21.992358518496392,1.9622919800730743)-- (16.356281566950525,11.96426361551843);
\draw (23.453339280609924,4.432784888934413)-- (19.226281566950526,11.934263615518429);
\draw (15.54784889955316,3.2725384345037445)-- (22.722848899553156,3.197538434503744);
\draw (14.138829661666694,5.773031343365085)-- (24.183829661666692,5.668031343365083);
\draw (14.164810423780226,8.258524252226422)-- (24.209810423780226,8.153524252226422);
\draw (15.625791185893759,10.72901716108776)-- (22.800791185893758,10.654017161087761);
\draw (17.791281566950524,11.94926361551843)-- (14.138829661666694,5.773031343365085);
\draw (20.661281566950525,11.919263615518428)-- (15.54784889955316,3.2725384345037445);
\draw (22.800791185893758,10.654017161087761)-- (17.687358518496392,2.007291980073074);
\draw (20.557358518496393,1.9772919800730744)-- (24.209810423780226,8.153524252226422);
\draw (17.791281566950524,11.94926361551843)-- (22.722848899553156,3.197538434503744);
\draw (24.183829661666692,5.668031343365083)-- (20.661281566950525,11.919263615518428);
\draw (14.164810423780226,8.258524252226422)-- (17.687358518496392,2.007291980073074);
\draw (15.625791185893759,10.72901716108776)-- (20.557358518496393,1.9772919800730744);
\draw [color=zzttqq] (39.54065020783629,-11.209691934324564)-- (45.28065020783629,-11.269691934324564);
\draw [color=zzttqq] (45.28065020783629,-11.269691934324564)-- (48.20261173206335,-6.328706116601893);
\draw [color=zzttqq] (48.20261173206335,-6.328706116601893)-- (45.38457325629042,-1.3277202988792176);
\draw [color=zzttqq] (45.38457325629042,-1.3277202988792176)-- (39.64457325629043,-1.2677202988792153);
\draw [color=zzttqq] (39.64457325629043,-1.2677202988792153)-- (36.72261173206336,-6.208706116601885);
\draw [color=zzttqq] (36.72261173206336,-6.208706116601885)-- (39.54065020783629,-11.209691934324564);
\draw (38.131630969949825,-8.709199025463224)-- (46.741630969949824,-8.799199025463228);
\draw (36.72261173206336,-6.208706116601885)-- (48.20261173206335,-6.328706116601893);
\draw (38.18359249417689,-3.73821320774055)-- (46.79359249417689,-3.8282132077405553);
\draw (42.514573256290426,-1.2977202988792165)-- (38.131630969949825,-8.709199025463224);
\draw (45.38457325629042,-1.3277202988792176)-- (39.54065020783629,-11.209691934324564);
\draw (46.79359249417689,-3.8282132077405553)-- (42.41065020783629,-11.239691934324565);
\draw (42.41065020783629,-11.239691934324565)-- (38.18359249417689,-3.73821320774055);
\draw (45.28065020783629,-11.269691934324564)-- (39.64457325629043,-1.2677202988792153);
\draw (46.741630969949824,-8.799199025463228)-- (42.514573256290426,-1.2977202988792165);
\draw (38.83614058889306,-9.959445479893894)-- (46.011140588893056,-10.034445479893897);
\draw (37.42712135100659,-7.458952571032555)-- (47.47212135100659,-7.56395257103256);
\draw (37.453102113120124,-4.9734596621712175)-- (47.498102113120126,-5.078459662171224);
\draw (38.91408287523366,-2.5029667533098827)-- (46.08908287523366,-2.5779667533098864);
\draw (41.07957325629043,-1.282720298879216)-- (37.42712135100659,-7.458952571032555);
\draw (43.94957325629042,-1.312720298879217)-- (38.83614058889306,-9.959445479893894);
\draw (46.08908287523366,-2.5779667533098864)-- (40.975650207836296,-11.224691934324564);
\draw (43.845650207836286,-11.254691934324565)-- (47.498102113120126,-5.078459662171224);
\draw (41.07957325629043,-1.282720298879216)-- (46.011140588893056,-10.034445479893897);
\draw (47.47212135100659,-7.56395257103256)-- (43.94957325629042,-1.312720298879217);
\draw (37.453102113120124,-4.9734596621712175)-- (40.975650207836296,-11.224691934324564);
\draw (38.91408287523366,-2.5029667533098827)-- (43.845650207836286,-11.254691934324565);
\draw (14.950248964492088,7.007431850878312)-- (21.384696264337443,10.700985317991499);
\draw (21.384696264337443,10.700985317991499)-- (21.26403780782058,3.2547943983888614);
\draw (21.26403780782058,3.2547943983888614)-- (14.950248964492088,7.007431850878312);
\draw (17.062896378772784,10.710281215030923)-- (17.059505941649753,3.3871641515733835);
\draw (17.059505941649753,3.3871641515733835)-- (23.479320042723458,6.918277797795752);
\draw (23.479320042723458,6.918277797795752)-- (17.062896378772784,10.710281215030923);
\draw (21.35321294831601,8.159739681610596)-- (12.835495137193313,3.3008905644238844);
\draw (19.30369552080427,6.961925440882017)-- (19.216019976845974,14.484073421044528);
\draw (19.30369552080427,6.961925440882017)-- (25.592848899553143,3.167538434503745);
\draw (-10.734748363113829,-7.829376116581178)-- (-13.55757159818533,-10.122919995076769);
\draw [dotted] (12.835495137193313,3.3008905644238844)-- (25.592848899553143,3.167538434503745);
\draw [dotted] (19.216019976845974,14.484073421044528)-- (25.592848899553143,3.167538434503745);
\draw [dotted] (19.216019976845974,14.484073421044528)-- (12.835495137193313,3.3008905644238844);
\draw (31.88718837492109,7.095645076974299)-- (38.37453833323538,10.878656014025918);
\draw (38.37453833323538,10.878656014025918)-- (38.20097721824959,3.3430076244848475);
\draw (38.20097721824959,3.3430076244848475)-- (31.88718837492109,7.095645076974299);
\draw (34.01682501335818,10.768344792784408)-- (33.92326160826818,3.351624947235976);
\draw (33.92326160826818,3.351624947235976)-- (40.41625945315247,7.006491023891739);
\draw (40.41625945315247,7.006491023891739)-- (34.01682501335818,10.768344792784408);
\draw (38.249633896439015,8.319858333090353)-- (29.738136390109073,3.389462307845792);
\draw (36.27889609541216,7.049738724495298)-- (36.18726333190805,14.511409713742728);
\draw (36.27889609541216,7.049738724495298)-- (42.529788309982166,3.2557516605997323);
\draw [dotted] (36.18726333190805,14.511409713742728)-- (42.529788309982166,3.2557516605997323);
\draw [dotted] (36.18726333190805,14.511409713742728)-- (29.738136390109073,3.389462307845792);
\draw (36.27889609541216,7.049738724495298)-- (36.29909459219123,4.563763105610689);
\draw (36.27889609541216,7.049738724495298)-- (34.036670773597315,8.426517651035184);
\draw [dotted] (29.738136390109073,3.389462307845792)-- (42.529788309982166,3.2557516605997323);
\draw (19.303125208733007,4.476166917560305)-- (19.30369552080427,6.961925440882017);
\draw [dotted] (56.246482195544104,4.802123029628403)-- (52.81419372053058,12.427077698184148);
\draw [dotted] (56.246482195544104,4.802123029628403)-- (47.89696756574237,5.649279484827429);
\draw [dotted] (47.89696756574237,5.649279484827429)-- (52.81419372053058,12.427077698184148);
\draw [dotted] (58.66170590460305,9.038178591505785)-- (52.81419372053058,12.427077698184148);
\draw [dotted] (58.66170590460305,9.038178591505785)-- (56.246482195544104,4.802123029628403);
\draw [dotted] (56.246482195544104,4.802123029628403)-- (52.983449060242094,1.8862210789329936);
\draw [dotted] (47.89696756574237,5.649279484827429)-- (52.983449060242094,1.8862210789329936);
\draw (51.25583056895894,10.279063624396207)-- (50.77292696903647,5.357479580112634);
\draw (53.95163648993808,9.900206889408512)-- (53.36782789867874,5.094196362989796);
\draw (53.95163648993808,9.900206889408512)-- (49.83297744184828,8.317833638378813);
\draw (49.83297744184828,8.317833638378813)-- (53.36782789867874,5.094196362989796);
\draw (50.77292696903647,5.357479580112634)-- (54.96320219984322,7.652976582736324);
\draw (54.96320219984322,7.652976582736324)-- (51.25583056895894,10.279063624396207);
\draw (50.77292696903647,5.357479580112634)-- (51.651822796183744,2.871378974797398);
\draw (51.651822796183744,2.871378974797398)-- (54.81618856569558,3.5239882965722793);
\draw (54.81618856569558,3.5239882965722793)-- (50.77292696903647,5.357479580112634);
\draw (53.36782789867874,5.094196362989796)-- (49.89744835806506,4.169292568622072);
\draw (49.89744835806506,4.169292568622072)-- (53.71587072846738,2.540725548410907);
\draw (53.71587072846738,2.540725548410907)-- (53.36782789867874,5.094196362989796);
\draw (54.96320219984322,7.652976582736324)-- (58.12605171157044,8.098695878112393);
\draw (58.12605171157044,8.098695878112393)-- (55.325554256379675,10.971630114907969);
\draw (55.325554256379675,10.971630114907969)-- (54.96320219984322,7.652976582736324);
\draw (53.95163648993808,9.900206889408512)-- (57.20559161394702,6.484302933310147);
\draw (57.385233398046296,9.777952430532997)-- (57.20559161394702,6.484302933310147);
\draw (57.385233398046296,9.777952430532997)-- (53.95163648993808,9.900206889408512);
\draw (51.081426240736384,4.484838929852787)-- (53.61559933324358,3.276382367994172);
\draw (52.099154287608016,4.756077427305306)-- (52.52053460579036,3.050539497372712);
\draw (52.52053460579036,3.050539497372712)-- (52.983449060242094,1.8862210789329936);
\draw (51.081426240736384,4.484838929852787)-- (47.89696756574237,5.649279484827429);
\draw (51.42300034772199,3.518640698040508)-- (53.5005312678595,4.120597279206977);
\draw (53.5005312678595,4.120597279206977)-- (56.246482195544104,4.802123029628403);
\draw (52.550590931796826,9.36193047740375)-- (52.21366992318365,6.146740655013277);
\draw (52.550590931796826,9.36193047740375)-- (52.81419372053058,12.427077698184148);
\draw (53.58707989665431,6.899115690976821)-- (51.1115958522866,8.809074002615818);
\draw (53.58707989665431,6.899115690976821)-- (56.246482195544104,4.802123029628403);
\draw (53.78042748837801,8.490785589286041)-- (50.962339190527146,7.287902478658929);
\draw (50.962339190527146,7.287902478658929)-- (47.89696756574237,5.649279484827429);
\draw (55.07931417724401,8.716404861885048)-- (57.337751104744775,8.907386419158374);
\draw (56.45691541347224,9.8110055278842)-- (55.95868882575286,7.7932638761442385);
\draw (57.28719598417677,7.980481532994249)-- (55.20370178040246,9.855626652223101);
\draw (55.20370178040246,9.855626652223101)-- (52.81419372053058,12.427077698184148);
\draw (57.337751104744775,8.907386419158374)-- (58.66170590460305,9.038178591505785);
\begin{scriptsize}
\draw [fill=qqqqff] (19.30369552080427,6.961925440882017) circle (1.5pt);
\draw [fill=qqqqff] (16.99404818912852,5.769263041464302) circle (1.5pt);
\draw [fill=qqqqff] (21.400727381363467,5.844973867095878) circle (1.5pt);
\draw [fill=qqqqff] (21.35321294831601,8.159739681610596) circle (1.5pt);
\draw [fill=qqqqff] (17.09973136316831,8.338304424939203) circle (1.5pt);
\draw [fill=qqqqff] (19.303125208733007,4.476166917560305) circle (1.5pt);
\draw [fill=xdxdff] (42.368465850734836,-6.26772201331274) circle (1.5pt);
\draw [fill=xdxdff] (42.368630928385016,-3.78195925409116) circle (1.5pt);
\draw [fill=xdxdff] (40.27945292459008,-7.457597587525383) circle (1.5pt);
\draw [fill=xdxdff] (44.63818582386217,-7.472968178382072) circle (1.5pt);
\draw [fill=xdxdff] (44.53678225468092,-4.886400371335734) circle (1.5pt);
\draw [fill=xdxdff] (40.33718980566702,-4.979637620538738) circle (1.5pt);
\draw [fill=xdxdff] (42.54470795025088,-8.75532874999599) circle (1.5pt);
\draw [fill=ffqqqq] (17.059505941649753,3.3871641515733835) circle (1.5pt);
\draw [fill=ffqqqq] (17.062896378772784,10.710281215030923) circle (1.5pt);
\draw [fill=ffqqqq] (21.384696264337443,10.700985317991499) circle (1.5pt);
\draw [fill=ffqqqq] (21.26403780782058,3.2547943983888614) circle (1.5pt);
\draw [fill=ffqqqq] (14.950248964492088,7.007431850878312) circle (1.5pt);
\draw [fill=ffqqqq] (23.479320042723458,6.918277797795752) circle (1.5pt);
\draw [fill=ffffqq] (19.216019976845974,14.484073421044528) circle (1.5pt);
\draw [fill=ffffqq] (12.835495137193313,3.3008905644238844) circle (1.5pt);
\draw [fill=ffffqq] (25.592848899553143,3.167538434503745) circle (1.5pt);
\draw [fill=qqqqff] (-0.7666538142675855,-10.652199351652675) circle (1.5pt);
\draw [fill=qqqqff] (-10.734748363113829,-7.829376116581178) circle (1.5pt);
\draw [fill=qqqqff] (-13.55757159818533,-10.122919995076769) circle (1.5pt);
\draw [fill=qqqqff] (36.27889609541216,7.049738724495298) circle (1.5pt);
\draw [fill=qqqqff] (36.25321043763535,9.535771700463505) circle (1.5pt);
\draw [fill=qqqqff] (33.93098759955753,5.857476267560285) circle (1.5pt);
\draw [fill=qqqqff] (38.337666791792486,5.933187093191874) circle (1.5pt);
\draw [fill=qqqqff] (38.249633896439015,8.319858333090353) circle (1.5pt);
\draw [fill=qqqqff] (34.036670773597315,8.426517651035184) circle (1.5pt);
\draw [fill=qqqqff] (36.29909459219123,4.563763105610689) circle (1.5pt);
\draw [fill=ffqqqq] (34.01682501335818,10.768344792784408) circle (1.5pt);
\draw [fill=ffqqqq] (38.37453833323538,10.878656014025918) circle (1.5pt);
\draw [fill=ffqqqq] (38.20097721824959,3.3430076244848475) circle (1.5pt);
\draw [fill=ffqqqq] (31.88718837492109,7.095645076974299) circle (1.5pt);
\draw [fill=ffqqqq] (40.41625945315247,7.006491023891739) circle (1.5pt);
\draw [fill=ffffqq] (36.18726333190805,14.511409713742728) circle (1.5pt);
\draw [fill=ffffqq] (29.738136390109073,3.389462307845792) circle (1.5pt);
\draw [fill=ffffqq] (42.529788309982166,3.2557516605997323) circle (1.5pt);
\draw [fill=ffqqqq] (34.01694868475073,3.34473604692619) circle (1.5pt);
\draw [fill=ffffqq] (56.246482195544104,4.802123029628403) circle (1.5pt);
\draw [fill=ffffqq] (52.81419372053058,12.427077698184148) circle (1.5pt);
\draw [fill=ffffqq] (47.89696756574237,5.649279484827429) circle (1.5pt);
\draw [fill=ffffqq] (58.66170590460305,9.038178591505785) circle (1.5pt);
\draw [fill=ffffqq] (52.983449060242094,1.8862210789329936) circle (1.5pt);
\draw [fill=ffqqqq] (51.25583056895894,10.279063624396207) circle (1.5pt);
\draw [fill=ffqqqq] (50.77292696903647,5.357479580112634) circle (1.5pt);
\draw [fill=ffqqqq] (53.95163648993808,9.900206889408512) circle (1.5pt);
\draw [fill=ffqqqq] (53.36782789867874,5.094196362989796) circle (1.5pt);
\draw [fill=ffqqqq] (49.83297744184828,8.317833638378813) circle (1.5pt);
\draw [fill=ffqqqq] (54.96320219984322,7.652976582736324) circle (1.5pt);
\draw [fill=ffqqqq] (51.651822796183744,2.871378974797398) circle (1.5pt);
\draw [fill=ffqqqq] (54.81618856569558,3.5239882965722793) circle (1.5pt);
\draw [fill=ffqqqq] (49.89744835806506,4.169292568622072) circle (1.5pt);
\draw [fill=ffqqqq] (53.71587072846738,2.540725548410907) circle (1.5pt);
\draw [fill=ffqqqq] (58.12605171157044,8.098695878112393) circle (1.5pt);
\draw [fill=ffqqqq] (55.325554256379675,10.971630114907969) circle (1.5pt);
\draw [fill=ffqqqq] (57.20559161394702,6.484302933310147) circle (1.5pt);
\draw [fill=ffqqqq] (57.385233398046296,9.777952430532997) circle (1.5pt);
\draw [fill=qqqqff] (51.1115958522866,8.809074002615818) circle (1.5pt);
\draw [fill=qqqqff] (52.550590931796826,9.36193047740375) circle (1.5pt);
\draw [fill=qqqqff] (50.962339190527146,7.287902478658929) circle (1.5pt);
\draw [fill=qqqqff] (52.21366992318365,6.146740655013277) circle (1.5pt);
\draw [fill=qqqqff] (51.081426240736384,4.484838929852787) circle (1.5pt);
\draw [fill=qqqqff] (52.099154287608016,4.756077427305306) circle (1.5pt);
\draw [fill=qqqqff] (51.42300034772199,3.518640698040508) circle (1.5pt);
\draw [fill=qqqqff] (53.5005312678595,4.120597279206977) circle (1.5pt);
\draw [fill=qqqqff] (53.61559933324358,3.276382367994172) circle (1.5pt);
\draw [fill=qqqqff] (52.52053460579036,3.050539497372712) circle (1.5pt);
\draw [fill=qqqqff] (55.95868882575286,7.7932638761442385) circle (1.5pt);
\draw [fill=qqqqff] (57.28719598417677,7.980481532994249) circle (1.5pt);
\draw [fill=qqqqff] (57.337751104744775,8.907386419158374) circle (1.5pt);
\draw [fill=qqqqff] (56.45691541347224,9.8110055278842) circle (1.5pt);
\draw [fill=qqqqff] (55.20370178040246,9.855626652223101) circle (1.5pt);
\draw [fill=qqqqff] (55.07931417724401,8.716404861885048) circle (1.5pt);
\draw [fill=qqqqff] (53.78042748837801,8.490785589286041) circle (1.5pt);
\draw [fill=qqqqff] (53.78042748837801,8.490785589286041) circle (1.5pt);
\draw [fill=qqqqff] (53.58707989665431,6.899115690976821) circle (1.5pt);
\draw [fill=qqqqff] (52.339288942631434,3.784131788742063) circle (1.5pt);
\draw [fill=qqqqff] (56.23806597747249,8.924698708134134) circle (1.5pt);
\draw [fill=qqqqff] (52.389349280835454,7.823224068888821) circle (1.5pt);
\draw [fill=qqqqff] (19.200300804836992,9.448770706657092) circle (1.5pt);
\draw [fill=qqqqff] (29.063943715182283,-1.98537499226316) circle (1.5pt);
\end{scriptsize}
\end{tikzpicture}
\end{center}

\begin{center}
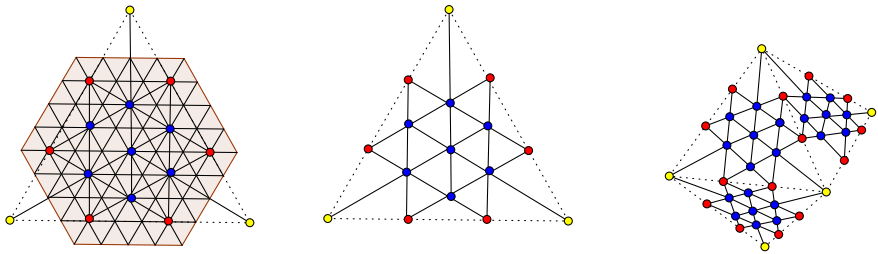
Figure 10. Construction of a 1-layer Platonic configuration on octahedron with 4-valent points in the vertices and 3-valent points in the interiors of faces.
\end{center}

\textbf{Another solution.} Place the copies of  $A = (15_{3}3_{1}, 12_{4})$ on each of the 8 faces of the octahedron and identify the 1-valent points of $A$ with the vertices (Figure 11). Now these 1-valent points produce 6 4-valent points, and we have a configuration $(6_{4}120_{3},96_{4})$. To get a balanced $(n_{4})$ configuration, take two concentric copies of the octahedron, and connect 60 antipodal pairs of 3-valent points with 60-valent lines. Thus we get a $(252_{4})$ configuration.

\begin{center}

\definecolor{ffffqq}{rgb}{1.0,1.0,0.0}
\definecolor{zzttqq}{rgb}{0.6,0.2,0.0}
\definecolor{xdxdff}{rgb}{0.49019607843137253,0.49019607843137253,1.0}
\definecolor{qqqqff}{rgb}{0.0,0.0,1.0}
\definecolor{ffqqqq}{rgb}{1.0,0.0,0.0}
\begin{tikzpicture}[line cap=round,line join=round,>=triangle 45,x=0.3cm,y=0.3cm]
\clip(-16.125942440181074,8.18782953673407) rectangle (5.83548785425175,18.6925495579874);
\fill[color=zzttqq,fill=zzttqq,fill opacity=0.1] (-2.043052161688277,13.407362297294329) -- (1.8969478383117215,13.387362297294327) -- (-0.055731653612587626,16.80950238820502) -- cycle;
\fill[color=zzttqq,fill=zzttqq,fill opacity=0.1] (-12.711908433273729,12.652335545766745) -- (-8.771908433273735,12.632335545766745) -- (-10.724587925198042,16.054475636677427) -- cycle;
\draw (3.6767086155662945,2.0924973326786414)-- (0.8567086155662946,2.1324973326786405);
\draw (-0.20329138443370542,2.5724973326786427)-- (0.8567086155662946,2.1324973326786405);
\draw (3.6767086155662945,2.0924973326786414)-- (4.056708615566295,2.7524973326786433);
\draw (4.056708615566295,2.7524973326786433)-- (2.3767086155662946,3.232497332678641);
\draw (2.3767086155662946,3.232497332678641)-- (-0.20329138443370542,2.5724973326786427);
\draw (2.5767086155662944,2.932497332678643)-- (1.1567086155662947,2.2924973326786406);
\draw (1.1567086155662947,2.2924973326786406)-- (3.576708615566295,2.652497332678641);
\draw (3.576708615566295,2.652497332678641)-- (0.3367086155662946,2.5924973326786422);
\draw (0.3367086155662946,2.5924973326786422)-- (3.096708615566295,2.2524973326786415);
\draw (3.096708615566295,2.2524973326786415)-- (2.5767086155662944,2.932497332678643);
\draw (0.3367086155662946,2.5924973326786422)-- (-0.20329138443370542,2.5724973326786427);
\draw (2.5767086155662944,2.932497332678643)-- (2.3767086155662946,3.232497332678641);
\draw (1.1567086155662947,2.2924973326786406)-- (0.8567086155662946,2.1324973326786405);
\draw (3.096708615566295,2.2524973326786415)-- (3.6767086155662945,2.0924973326786414);
\draw (3.576708615566295,2.652497332678641)-- (4.056708615566295,2.7524973326786433);
\draw (1.1567086155662947,2.2924973326786406)-- (2.3167086155662946,-2.267502667321359);
\draw (3.096708615566295,2.2524973326786415)-- (2.3167086155662946,-2.267502667321359);
\draw (3.576708615566295,2.652497332678641)-- (2.3167086155662946,-2.267502667321359);
\draw (2.5767086155662944,2.932497332678643)-- (2.3167086155662946,-2.267502667321359);
\draw (0.3367086155662946,2.5924973326786422)-- (2.3167086155662946,-2.267502667321359);
\draw (2.497339538259561,1.3451157865439796)-- (1.623606229261588,0.4571067133247284);
\draw (1.623606229261588,0.4571067133247284)-- (3.161587515367155,1.0315482747581908);
\draw (3.161587515367155,1.0315482747581908)-- (1.0406127549998456,0.8647326267962891);
\draw (1.0406127549998456,0.8647326267962891)-- (2.806560032106249,0.5711234900640165);
\draw (2.806560032106249,0.5711234900640165)-- (2.497339538259561,1.3451157865439796);
\draw (2.8186119396518703,-0.30768968755863524)-- (1.4577404238451614,-0.15912619673312234);
\draw (2.582349145553784,-0.7281498525220624)-- (1.4577404238451614,-0.15912619673312234);
\draw (2.582349145553784,-0.7281498525220624)-- (2.4393345507283883,0.18501603592051685);
\draw (2.4393345507283883,0.18501603592051685)-- (1.944100556066529,-0.8027675368740046);
\draw (1.944100556066529,-0.8027675368740046)-- (2.8186119396518703,-0.30768968755863524);
\draw [color=zzttqq] (-2.043052161688277,13.407362297294329)-- (1.8969478383117215,13.387362297294327);
\draw [color=zzttqq] (1.8969478383117215,13.387362297294327)-- (-0.055731653612587626,16.80950238820502);
\draw [color=zzttqq] (-0.055731653612587626,16.80950238820502)-- (-2.043052161688277,13.407362297294329);
\draw (-0.055731653612587626,16.80950238820502)-- (-0.07305216168827766,13.397362297294329);
\draw (-1.0493919076504323,15.108432342749673)-- (1.8969478383117215,13.387362297294327);
\draw (0.920608092349567,15.098432342749673)-- (-2.043052161688277,13.407362297294329);
\draw (-1.5462220346693547,14.257897320022)-- (0.9119478383117219,13.392362297294328);
\draw (-0.55256178063151,15.958967365477346)-- (-1.0580521616882772,13.402362297294328);
\draw (0.4324382193684897,15.953967365477347)-- (0.9119478383117219,13.392362297294328);
\draw (-1.5462220346693547,14.257897320022)-- (0.4324382193684897,15.953967365477347);
\draw (-0.55256178063151,15.958967365477346)-- (1.4087779653306443,14.242897320022);
\draw (1.4087779653306443,14.242897320022)-- (-1.0580521616882772,13.402362297294328);
\draw (-0.07132011088070872,13.738576306385395)-- (-0.07201390257723628,8.644913138211654);
\draw (0.6242420669457823,14.929325338204134)-- (4.670627129659983,17.768468344295137);
\draw (-0.7547579330542168,14.93632533820414)-- (-4.8548468079690155,17.728276471140582);
\draw [dotted] (-4.8548468079690155,17.728276471140582)-- (-0.07201390257723628,8.644913138211654);
\draw [dotted] (-0.07201390257723628,8.644913138211654)-- (4.670627129659983,17.768468344295137);
\draw [dotted] (-4.8548468079690155,17.728276471140582)-- (4.670627129659983,17.768468344295137);
\draw [color=zzttqq] (-12.711908433273729,12.652335545766745)-- (-8.771908433273735,12.632335545766745);
\draw [color=zzttqq] (-8.771908433273735,12.632335545766745)-- (-10.724587925198042,16.054475636677427);
\draw [color=zzttqq] (-10.724587925198042,16.054475636677427)-- (-12.711908433273729,12.652335545766745);
\draw (-10.724587925198042,16.054475636677427)-- (-10.741908433273732,12.642335545766745);
\draw (-11.718248179235886,14.353405591222085)-- (-8.771908433273735,12.632335545766745);
\draw (-9.748248179235889,14.343405591222087)-- (-12.711908433273729,12.652335545766745);
\draw (-12.215078306254807,13.502870568494416)-- (-9.756908433273733,12.637335545766746);
\draw (-11.221418052216965,15.203940613949756)-- (-11.726908433273731,12.647335545766744);
\draw (-10.236418052216965,15.198940613949757)-- (-9.756908433273733,12.637335545766746);
\draw (-12.215078306254807,13.502870568494416)-- (-10.236418052216965,15.198940613949757);
\draw (-11.221418052216965,15.203940613949756)-- (-9.260078306254812,13.487870568494415);
\draw (-9.260078306254812,13.487870568494415)-- (-11.726908433273731,12.647335545766744);
\draw (-10.74017638246617,12.983549554857824)-- (-10.709446179222306,9.599401289589986);
\draw (-10.044614204639668,14.174298586676542)-- (-6.343421920388875,16.624432803803163);
\draw (-11.423614204639676,14.181298586676554)-- (-15.108297688122155,16.394642053338245);
\begin{scriptsize}
\draw [fill=ffqqqq] (3.6767086155662945,2.0924973326786414) circle (1.5pt);
\draw [fill=ffqqqq] (0.8567086155662946,2.1324973326786405) circle (1.5pt);
\draw [fill=ffqqqq] (-0.20329138443370542,2.5724973326786427) circle (1.5pt);
\draw [fill=ffqqqq] (4.056708615566295,2.7524973326786433) circle (1.5pt);
\draw [fill=ffqqqq] (2.3767086155662946,3.232497332678641) circle (1.5pt);
\draw [fill=qqqqff] (2.3167086155662946,-2.267502667321359) circle (1.5pt);
\draw [fill=qqqqff] (2.5767086155662944,2.932497332678643) circle (1.5pt);
\draw [fill=qqqqff] (1.1567086155662947,2.2924973326786406) circle (1.5pt);
\draw [fill=qqqqff] (3.576708615566295,2.652497332678641) circle (1.5pt);
\draw [fill=qqqqff] (0.3367086155662946,2.5924973326786422) circle (1.5pt);
\draw [fill=qqqqff] (3.096708615566295,2.2524973326786415) circle (1.5pt);
\draw [fill=xdxdff] (2.497339538259561,1.3451157865439796) circle (1.5pt);
\draw [fill=xdxdff] (1.623606229261588,0.4571067133247284) circle (1.5pt);
\draw [fill=xdxdff] (3.161587515367155,1.0315482747581908) circle (1.5pt);
\draw [fill=xdxdff] (1.0406127549998456,0.8647326267962891) circle (1.5pt);
\draw [fill=xdxdff] (2.806560032106249,0.5711234900640165) circle (1.5pt);
\draw [fill=xdxdff] (2.8186119396518703,-0.30768968755863524) circle (1.5pt);
\draw [fill=xdxdff] (1.4577404238451614,-0.15912619673312234) circle (1.5pt);
\draw [fill=xdxdff] (2.582349145553784,-0.7281498525220624) circle (1.5pt);
\draw [fill=xdxdff] (2.4393345507283883,0.18501603592051685) circle (1.5pt);
\draw [fill=xdxdff] (1.944100556066529,-0.8027675368740046) circle (1.5pt);
\draw [fill=qqqqff] (-2.043052161688277,13.407362297294329) circle (1.5pt);
\draw [fill=qqqqff] (1.8969478383117215,13.387362297294327) circle (1.5pt);
\draw [fill=qqqqff] (-0.055731653612587626,16.80950238820502) circle (1.5pt);
\draw [fill=qqqqff] (-1.5462220346693547,14.257897320022) circle (1.5pt);
\draw [fill=qqqqff] (-0.55256178063151,15.958967365477346) circle (1.5pt);
\draw [fill=qqqqff] (0.4324382193684897,15.953967365477347) circle (1.5pt);
\draw [fill=qqqqff] (1.4087779653306443,14.242897320022) circle (1.5pt);
\draw [fill=qqqqff] (0.9119478383117219,13.392362297294328) circle (1.5pt);
\draw [fill=qqqqff] (-1.0580521616882772,13.402362297294328) circle (1.5pt);
\draw [fill=qqqqff] (-0.7547579330542168,14.93632533820414) circle (1.5pt);
\draw [fill=qqqqff] (-0.9316795664240857,14.041513564340082) circle (1.5pt);
\draw [fill=qqqqff] (-0.07132011088070872,13.738576306385395) circle (1.5pt);
\draw [fill=qqqqff] (0.7920704335759136,14.032763564340081) circle (1.5pt);
\draw [fill=qqqqff] (0.6242420669457823,14.929325338204134) circle (1.5pt);
\draw [fill=qqqqff] (-0.06222684414097143,15.529949854113509) circle (1.5pt);
\draw [fill=ffffqq] (-0.07201390257723628,8.644913138211654) circle (1.5pt);
\draw [fill=ffffqq] (4.670627129659983,17.768468344295137) circle (1.5pt);
\draw [fill=ffffqq] (-4.8548468079690155,17.728276471140582) circle (1.5pt);
\draw [fill=qqqqff] (-12.711908433273729,12.652335545766745) circle (1.5pt);
\draw [fill=qqqqff] (-8.771908433273735,12.632335545766745) circle (1.5pt);
\draw [fill=qqqqff] (-10.724587925198042,16.054475636677427) circle (1.5pt);
\draw [fill=qqqqff] (-12.215078306254807,13.502870568494416) circle (1.5pt);
\draw [fill=qqqqff] (-11.221418052216965,15.203940613949756) circle (1.5pt);
\draw [fill=qqqqff] (-10.236418052216965,15.198940613949757) circle (1.5pt);
\draw [fill=qqqqff] (-9.260078306254812,13.487870568494415) circle (1.5pt);
\draw [fill=qqqqff] (-9.756908433273733,12.637335545766746) circle (1.5pt);
\draw [fill=qqqqff] (-11.726908433273731,12.647335545766744) circle (1.5pt);
\draw [fill=qqqqff] (-11.423614204639676,14.181298586676554) circle (1.5pt);
\draw [fill=qqqqff] (-11.600535838009545,13.286486812812495) circle (1.5pt);
\draw [fill=qqqqff] (-10.74017638246617,12.983549554857824) circle (1.5pt);
\draw [fill=qqqqff] (-9.876785838009539,13.277736812812499) circle (1.5pt);
\draw [fill=qqqqff] (-10.044614204639668,14.174298586676542) circle (1.5pt);
\draw [fill=qqqqff] (-10.731083115726433,14.774923102585923) circle (1.5pt);
\draw [fill=ffffqq] (-10.709446179222306,9.599401289589986) circle (1.5pt);
\draw [fill=ffffqq] (-6.343421920388875,16.624432803803163) circle (1.5pt);
\draw [fill=ffffqq] (-15.108297688122155,16.394642053338245) circle (1.5pt);
\end{scriptsize}
\end{tikzpicture}

\end{center}
\vspace{-5mm}
\begin{center}

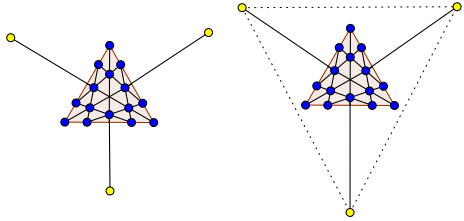
Figure 11. Placing 1-valent points into 4-valent vertices of the octahedron.
\end{center}

\begin{prob}\label{TOI4}
To find $((23f)_{4})$ configurations in the classes $T_{4},O_{4}$ and $I_{4}$.
\end{prob}

\textbf{Solution.} Place a copy of $A = (10_{3}6_{2},9_{4})$ on each face of $T,O$ or $I$ in such a way that the 2-valent vertices are on the edges to get a connected configuration $B = ((2e)_{4}(10f)_{3}, (9f)_{4})$ (see Figure 12) having ten 3-valent points on each face and two 4-valent points on each edge. Hence the number of points of $B$ is $10f + 2e$ and the number of lines of $B$ is $9f$. To get a balanced configuration $C$, take two concentric copies of $B$ and connect 3-valent points with $\frac{10f}{2} = 5f$ 4-valent lines passing through the center of $B$. Now all points in $C$ are 4-valent and their number is $p = 20f + 4e = 23f$; the number of lines in $C$ is $n = 18f + 5f = 23f$.

\begin{center}
\definecolor{qqqqff}{rgb}{0.0,0.0,1.0}
\definecolor{ffqqqq}{rgb}{1.0,0.0,0.0}
\definecolor{uuuuuu}{rgb}{0.26666666666666666,0.26666666666666666,0.26666666666666666}
\begin{tikzpicture}[line cap=round,line join=round,>=triangle 45,x=0.75cm,y=0.75cm]
\clip(-3.58,-0.6) rectangle (9.859999999999998,6.12);
\draw (-0.7499999999999987,4.383739178243871)-- (-1.5974999999999997,-0.02);
\draw (0.0975000000000002,-0.02)-- (-0.7499999999999987,4.383739178243871);
\draw (1.3687500000000008,2.1818695891219355)-- (-2.86875,0.7139565297073118);
\draw (-2.86875,0.7139565297073118)-- (0.5212500000000011,3.6497826485365588);
\draw (1.3687500000000004,0.713956529707312)-- (-2.8687499999999995,2.1818695891219355);
\draw (1.3687500000000004,0.713956529707312)-- (-2.021249999999999,3.6497826485365588);
\draw (-2.8687499999999995,2.1818695891219355)-- (1.3687500000000008,2.1818695891219355);
\draw (-1.5974999999999997,-0.02)-- (0.5212500000000011,3.6497826485365588);
\draw (-2.021249999999999,3.6497826485365588)-- (0.0975000000000002,-0.02);
\draw (-0.7499999999999993,1.9372174125528319)-- (-0.32624999999999915,2.1818695891219355);
\draw (-2.86875,0.7139565297073118)-- (-0.7499999999999993,1.9372174125528319);
\draw (-1.1737499999999992,2.181869589121936)-- (1.3687500000000004,0.713956529707312);
\draw (-0.7499999999999987,4.383739178243871)-- (-0.7499999999999992,1.4479130594146237);
\draw (6.110000000000002,4.16373917824387)-- (5.2625,-0.24);
\draw (6.9575,-0.24)-- (6.110000000000002,4.16373917824387);
\draw (8.228750000000002,1.9618695891219353)-- (3.9912499999999995,0.4939565297073121);
\draw (3.9912499999999995,0.4939565297073121)-- (7.381250000000001,3.429782648536559);
\draw (8.22875,0.4939565297073118)-- (3.991250000000001,1.9618695891219353);
\draw (8.22875,0.4939565297073118)-- (4.838750000000001,3.429782648536559);
\draw (3.991250000000001,1.9618695891219353)-- (8.228750000000002,1.9618695891219353);
\draw (5.2625,-0.24)-- (7.381250000000001,3.429782648536559);
\draw (4.838750000000001,3.429782648536559)-- (6.9575,-0.24);
\draw (6.110000000000002,1.7172174125528312)-- (6.53375,1.9618695891219355);
\draw (3.9912499999999995,0.4939565297073121)-- (6.110000000000002,1.7172174125528312);
\draw (5.68625,1.9618695891219358)-- (8.22875,0.4939565297073118);
\draw (6.110000000000002,4.16373917824387)-- (6.110000000000002,1.2279130594146237);
\draw [dotted] (2.72,-0.24)-- (9.5,-0.24);
\draw [dotted] (9.5,-0.24)-- (6.110000000000002,5.6316522376584945);
\draw [dotted] (3.991250000000001,1.9618695891219353)-- (6.110000000000002,5.6316522376584945);
\draw [dotted] (3.991250000000001,1.9618695891219353)-- (2.72,-0.24);
\begin{scriptsize}
\draw [fill=uuuuuu] (-2.8687499999999995,2.1818695891219355) circle (1.5pt);
\draw [fill=uuuuuu] (-2.021249999999999,3.6497826485365588) circle (1.5pt);
\draw [fill=uuuuuu] (0.5212500000000011,3.6497826485365588) circle (1.5pt);
\draw [fill=uuuuuu] (1.3687500000000008,2.1818695891219355) circle (1.5pt);
\draw [fill=uuuuuu] (0.0975000000000002,-0.02) circle (1.5pt);
\draw [fill=uuuuuu] (-1.5974999999999997,-0.02) circle (1.5pt);
\draw [fill=uuuuuu] (-0.7499999999999993,1.9372174125528319) circle (1.5pt);
\draw [fill=uuuuuu] (-0.7499999999999987,4.383739178243871) circle (1.5pt);
\draw [fill=uuuuuu] (-2.86875,0.7139565297073118) circle (1.5pt);
\draw [fill=uuuuuu] (1.3687500000000004,0.713956529707312) circle (1.5pt);
\draw [fill=uuuuuu] (0.5212500000000012,3.649782648536559) circle (1.5pt);
\draw [fill=uuuuuu] (-0.7499999999999992,2.5488478539755914) circle (1.5pt);
\draw [fill=uuuuuu] (-1.1737499999999992,2.181869589121936) circle (1.5pt);
\draw [fill=uuuuuu] (-1.2796874999999992,1.6314021918414514) circle (1.5pt);
\draw [fill=uuuuuu] (-0.7499999999999992,1.4479130594146237) circle (1.5pt);
\draw [fill=uuuuuu] (-0.22031249999999938,1.6314021918414516) circle (1.5pt);
\draw [fill=uuuuuu] (-0.32624999999999915,2.1818695891219355) circle (1.5pt);
\draw [fill=uuuuuu] (-0.7499999999999989,1.9372174125528314) circle (1.5pt);
\draw [fill=ffqqqq] (4.838750000000001,3.429782648536559) circle (1.5pt);
\draw [fill=ffqqqq] (7.381250000000001,3.429782648536559) circle (1.5pt);
\draw [fill=ffqqqq] (8.228750000000002,1.9618695891219353) circle (1.5pt);
\draw [fill=ffqqqq] (6.9575,-0.24) circle (1.5pt);
\draw [fill=ffqqqq] (5.2625,-0.24) circle (1.5pt);
\draw [fill=qqqqff] (6.110000000000002,1.7172174125528312) circle (1.5pt);
\draw [fill=qqqqff] (6.110000000000002,4.16373917824387) circle (1.5pt);
\draw [fill=qqqqff] (3.9912499999999995,0.4939565297073121) circle (1.5pt);
\draw [fill=qqqqff] (8.22875,0.4939565297073118) circle (1.5pt);
\draw [fill=qqqqff] (6.11,2.3288478539755912) circle (1.5pt);
\draw [fill=qqqqff] (5.68625,1.9618695891219358) circle (1.5pt);
\draw [fill=qqqqff] (5.5803125,1.4114021918414519) circle (1.5pt);
\draw [fill=ffqqqq] (3.991250000000001,1.9618695891219355) circle (1.5pt);
\draw [fill=qqqqff] (6.110000000000002,1.2279130594146237) circle (1.5pt);
\draw [fill=qqqqff] (6.639687500000001,1.4114021918414512) circle (1.5pt);
\draw [fill=qqqqff] (6.53375,1.9618695891219355) circle (1.5pt);
\draw [fill=qqqqff] (6.11,1.7172174125528314) circle (1.5pt);
\end{scriptsize}
\end{tikzpicture}
\end{center}

\begin{center}
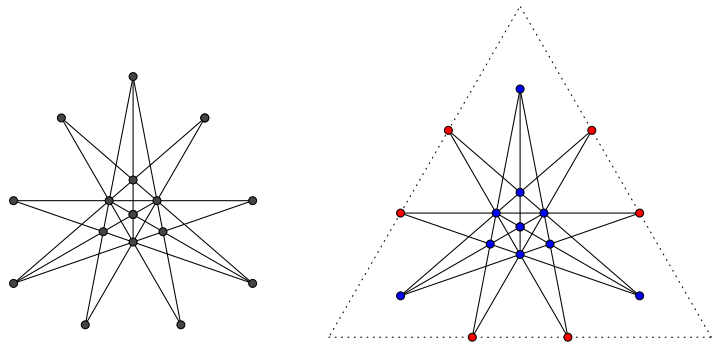
Figure 12. Face of a $(n_{4}) = ((23f)_{4})$ configuration on tetrahedron, octahedron or icosahedron. 
\end{center}

\begin{remark}
It is obvious how these solutions may be generalised to obtain examples of Platonic configurations with the full symmetry group of the cube and the dodecahedron (use the barycentric subdivision of square and pentagonal faces into 4 or 5 triangles and continue as in Example 2). 

It is equally obvious that Platonic configurations with the rotational symmetry of any Platonic solid can be obtained by exactly the same procedure, but starting from the configuration $A = (16_{3}12_{2},21_{4})$, shown in Figure 13 left. Note, however, that the configuration $A' = (24_{3}12_{2}, 27_{4})$ in Figure 13 right produces a spatial "Platonic" configuration only if its copies are placed on the triangular faces of the octahedron (and not on tetrahedron or icosahedron), since any two adjacent faces must be (because of the non-symmetrical arrangement of points along the edges) the mirror images of each other! Consequently the symmetry group of this configuration $B = O(48_{4})$ is smaller both from $Sym(O)$ and from $Sym_{R}(O)$, therefore the obtained configuration is \textit{not} Platonic (see Definition 1). However, the same configuration $A'$ placed on the triangulated faces of the cube produces a Platonic configuration with the full symmetry group of a cube.
\end{remark}

\vspace{-0mm}

\begin{center}
\definecolor{qqqqff}{rgb}{0.0,0.0,1.0}
\definecolor{ffqqqq}{rgb}{1.0,0.0,0.0}
\definecolor{uuuuuu}{rgb}{0.26666666666666666,0.26666666666666666,0.26666666666666666}
\definecolor{zzttqq}{rgb}{0.6,0.2,0.0}

\end{center}

\vspace{-0mm}

\begin{center}
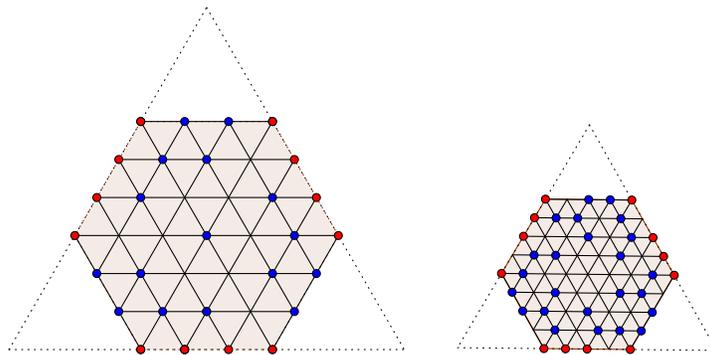
Figure 13. Left: the configuration $A = (16_{3}12_{2},21_{4})$ with 3-fold rotational symmetry (left) produces an octahedral 1-layer Platonic configuration $B \in O_{R}(48_{4}96_{3},168_{4})$. Right: the configuration $A' = (24_{3}12_{2}, 27_{4})$ produces a 1-layer Platonic configuration $B'$ with the full symmetry of a cube. 
\end{center}

\begin{example}\label{examplesesti}
An example of $(p_{3},n_{5})$ configuration is shown in Figure 14. Here the Pappus configuration is useful again. We exclude the edges. The points on vertices have now valence 4. Now we use radial projection with 5 concentric solids to raise the valence of interior points to 4, too, and we get a $(p_{4},n_{5})$ configuration. 
\end{example}

\begin{center}
\definecolor{ffqqqq}{rgb}{1.0,0.0,0.0}
\definecolor{qqqqff}{rgb}{0.0,0.0,1.0}
\definecolor{zzttqq}{rgb}{0.6,0.2,0.0}
\begin{tikzpicture}[line cap=round,line join=round,>=triangle 45,x=0.25cm,y=0.25cm]
\clip(-7.756513701756573,-5.246307358212073) rectangle (34.251011831058044,10.940989714716506);
\fill[color=zzttqq,fill=zzttqq,fill opacity=0.1] (1.1018462351790044,-4.540801646714921) -- (6.841846235179005,-4.600801646714921) -- (9.76380775940607,0.34018417100775666) -- (6.945769283633139,5.341169988730436) -- (1.2057692836331388,5.401169988730438) -- (-1.71619224059393,0.4601841710077621) -- cycle;
\fill[color=zzttqq,fill=zzttqq,fill opacity=0.1] (-4.638153764821004,-4.480801646714921) -- (12.581846235178993,-4.660801646714917) -- (4.127730807860206,10.342155806453121) -- cycle;
\fill[color=zzttqq,fill=zzttqq,fill opacity=0.1] (14.87598043659288,-4.6297645032142585) -- (32.0959804365929,-4.8097645032142715) -- (23.64186500927409,10.193192949953774) -- cycle;
\draw [color=zzttqq] (1.1018462351790044,-4.540801646714921)-- (6.841846235179005,-4.600801646714921);
\draw [color=zzttqq] (6.841846235179005,-4.600801646714921)-- (9.76380775940607,0.34018417100775666);
\draw [color=zzttqq] (9.76380775940607,0.34018417100775666)-- (6.945769283633139,5.341169988730436);
\draw [color=zzttqq] (6.945769283633139,5.341169988730436)-- (1.2057692836331388,5.401169988730438);
\draw [color=zzttqq] (1.2057692836331388,5.401169988730438)-- (-1.71619224059393,0.4601841710077621);
\draw [color=zzttqq] (-1.71619224059393,0.4601841710077621)-- (1.1018462351790044,-4.540801646714921);
\draw (-0.30717300270746284,-2.0403087378535796)-- (8.302826997292538,-2.130308737853582);
\draw (-1.71619224059393,0.4601841710077621)-- (9.76380775940607,0.34018417100775666);
\draw (-0.2552114784803956,2.9306770798691)-- (8.354788521519605,2.8406770798690966);
\draw (4.0757692836331385,5.371169988730437)-- (-0.30717300270746284,-2.0403087378535796);
\draw (6.945769283633139,5.341169988730436)-- (1.1018462351790044,-4.540801646714921);
\draw (8.354788521519605,2.8406770798690966)-- (3.9718462351790045,-4.570801646714921);
\draw (3.9718462351790045,-4.570801646714921)-- (-0.2552114784803956,2.9306770798691);
\draw (6.841846235179005,-4.600801646714921)-- (1.2057692836331388,5.401169988730438);
\draw (8.302826997292538,-2.130308737853582)-- (4.0757692836331385,5.371169988730437);
\draw (0.39733661623577077,-3.29055519228425)-- (7.572336616235772,-3.3655551922842517);
\draw (-1.0116826216506964,-0.7900622834229087)-- (9.033317378349304,-0.8950622834229127);
\draw (-0.9857018595371628,1.695430625438431)-- (9.059298140462838,1.5904306254384266);
\draw (0.4752789025763716,4.165923534299769)-- (7.650278902576372,4.090923534299766);
\draw (2.640769283633139,5.386169988730438)-- (-1.0116826216506964,-0.7900622834229087);
\draw (5.510769283633138,5.356169988730437)-- (0.39733661623577077,-3.29055519228425);
\draw (7.650278902576372,4.090923534299766)-- (2.5368462351790044,-4.5558016467149205);
\draw (5.406846235179005,-4.585801646714922)-- (9.059298140462838,1.5904306254384266);
\draw (2.640769283633139,5.386169988730438)-- (7.572336616235772,-3.3655551922842517);
\draw (9.033317378349304,-0.8950622834229127)-- (5.510769283633138,5.356169988730437);
\draw (-0.9857018595371628,1.695430625438431)-- (2.5368462351790044,-4.5558016467149205);
\draw (0.4752789025763716,4.165923534299769)-- (5.406846235179005,-4.585801646714922);
\draw [color=zzttqq] (-4.638153764821004,-4.480801646714921)-- (12.581846235178993,-4.660801646714917);
\draw [color=zzttqq] (12.581846235178993,-4.660801646714917)-- (4.127730807860206,10.342155806453121);
\draw [color=zzttqq] (4.127730807860206,10.342155806453121)-- (-4.638153764821004,-4.480801646714921);
\draw (-4.053761459975588,-3.492604483170383)-- (9.059298140462838,1.5904306254384264);
\draw (-0.307173002707463,-2.040308737853579)-- (8.354788521519605,2.840677079869097);
\draw (-3.490153764821002,-4.492801646714921)-- (7.650278902576372,4.090923534299766);
\draw (8.302826997292536,-2.130308737853582)-- (-0.25521147848039577,2.9306770798691);
\draw (12.018238540024415,-3.6606044831703866)-- (-0.9857018595371632,1.6954306254384313);
\draw (11.433846235178997,-4.648801646714918)-- (0.47527890257637095,4.165923534299769);
\draw (3.543338503014792,9.353958642908577)-- (5.406846235179005,-4.58580164671492);
\draw (2.536846235179004,-4.5558016467149205)-- (4.6913385030147925,9.341958642908581);
\draw (3.9718462351790054,-4.57080164671492)-- (4.0757692836331385,5.371169988730437);
\draw [color=zzttqq] (14.87598043659288,-4.6297645032142585)-- (32.0959804365929,-4.8097645032142715);
\draw [color=zzttqq] (32.0959804365929,-4.8097645032142715)-- (23.64186500927409,10.193192949953774);
\draw [color=zzttqq] (23.64186500927409,10.193192949953774)-- (14.87598043659288,-4.6297645032142585);
\draw (15.460372741438283,-3.64156733966973)-- (28.57343234187672,1.4414677689390856);
\draw (19.206961198706416,-2.1892715943529213)-- (27.868922722933487,2.6917142233697544);
\draw (16.023980436592883,-4.64176450321426)-- (27.164413103990253,3.9419606778004277);
\draw (27.816961198706423,-2.279271594352925)-- (19.258922722933487,2.78171422336976);
\draw (31.53237274143831,-3.8095673396697314)-- (18.52843234187672,1.5464677689390909);
\draw (30.947980436592893,-4.797764503214271)-- (19.98941310399026,4.016960677800432);
\draw (23.057472704428676,9.204995786409222)-- (24.920980436592888,-4.734764503214265);
\draw (22.050980436592884,-4.704764503214263)-- (24.205472704428672,9.192995786409236);
\draw (23.48598043659289,-4.7197645032142646)-- (23.58990348504702,5.222207132231098);
\begin{scriptsize}
\draw [fill=qqqqff] (-0.307173002707463,-2.040308737853579) circle (1.5pt);
\draw [fill=ffqqqq] (8.354788521519605,2.840677079869097) circle (1.5pt);
\draw [fill=qqqqff] (8.302826997292536,-2.130308737853582) circle (1.5pt);
\draw [fill=ffqqqq] (-0.25521147848039577,2.9306770798691) circle (1.5pt);
\draw [fill=qqqqff] (4.0757692836331385,5.371169988730437) circle (1.5pt);
\draw [fill=ffqqqq] (3.9718462351790054,-4.57080164671492) circle (1.5pt);
\draw [fill=ffqqqq] (7.650278902576372,4.090923534299766) circle (1.5pt);
\draw [fill=ffqqqq] (9.059298140462838,1.5904306254384264) circle (1.5pt);
\draw [fill=ffqqqq] (0.47527890257637095,4.165923534299769) circle (1.5pt);
\draw [fill=ffqqqq] (-0.9857018595371632,1.6954306254384313) circle (1.5pt);
\draw [fill=ffqqqq] (2.536846235179004,-4.5558016467149205) circle (1.5pt);
\draw [fill=ffqqqq] (5.406846235179005,-4.58580164671492) circle (1.5pt);
\draw [fill=qqqqff] (4.02380775940607,0.40018417100775844) circle (1.5pt);
\draw [fill=qqqqff] (4.6901124920389226,0.7756446185248875) circle (1.5pt);
\draw [fill=qqqqff] (3.2363567117490653,-0.043541812421575496) circle (1.5pt);
\draw [fill=qqqqff] (4.0158136787557535,-0.3645828778726531) circle (1.5pt);
\draw [fill=qqqqff] (3.3654971074235376,0.7894907723710418) circle (1.5pt);
\draw [fill=qqqqff] (4.033255309265538,1.3039997742300649) circle (1.5pt);
\draw [fill=qqqqff] (4.801811257203611,-0.05990544878521281) circle (1.5pt);
\draw [fill=ffqqqq] (-3.490153764821002,-4.492801646714921) circle (1.5pt);
\draw [fill=ffqqqq] (-4.053761459975588,-3.492604483170383) circle (1.5pt);
\draw [fill=ffqqqq] (11.433846235178997,-4.648801646714918) circle (1.5pt);
\draw [fill=ffqqqq] (12.018238540024415,-3.6606044831703866) circle (1.5pt);
\draw [fill=ffqqqq] (4.6913385030147925,9.341958642908581) circle (1.5pt);
\draw [fill=ffqqqq] (3.543338503014792,9.353958642908577) circle (1.5pt);
\draw [fill=ffqqqq] (27.868922722933487,2.6917142233697544) circle (1.5pt);
\draw [fill=qqqqff] (27.816961198706423,-2.279271594352925) circle (1.5pt);
\draw [fill=ffqqqq] (19.258922722933487,2.78171422336976) circle (1.5pt);
\draw [fill=qqqqff] (23.58990348504702,5.222207132231098) circle (1.5pt);
\draw [fill=ffqqqq] (23.48598043659289,-4.7197645032142646) circle (1.5pt);
\draw [fill=ffqqqq] (27.164413103990253,3.9419606778004277) circle (1.5pt);
\draw [fill=ffqqqq] (28.57343234187672,1.4414677689390856) circle (1.5pt);
\draw [fill=ffqqqq] (19.98941310399026,4.016960677800432) circle (1.5pt);
\draw [fill=ffqqqq] (18.52843234187672,1.5464677689390909) circle (1.5pt);
\draw [fill=ffqqqq] (22.050980436592884,-4.704764503214263) circle (1.5pt);
\draw [fill=ffqqqq] (24.920980436592888,-4.734764503214265) circle (1.5pt);
\draw [fill=qqqqff] (23.537941960819957,0.2512213145084194) circle (1.5pt);
\draw [fill=qqqqff] (24.20424669345281,0.6266817620255479) circle (1.5pt);
\draw [fill=qqqqff] (22.750490913162945,-0.1925046689209164) circle (1.5pt);
\draw [fill=qqqqff] (23.529947880169637,-0.5135457343719942) circle (1.5pt);
\draw [fill=qqqqff] (22.879631308837418,0.6405279158717031) circle (1.5pt);
\draw [fill=qqqqff] (23.54738951067942,1.1550369177307278) circle (1.5pt);
\draw [fill=qqqqff] (24.315945458617495,-0.2088683052845524) circle (1.5pt);
\draw [fill=ffqqqq] (16.023980436592883,-4.64176450321426) circle (1.5pt);
\draw [fill=ffqqqq] (15.460372741438283,-3.64156733966973) circle (1.5pt);
\draw [fill=ffqqqq] (30.947980436592893,-4.797764503214271) circle (1.5pt);
\draw [fill=ffqqqq] (31.53237274143831,-3.8095673396697314) circle (1.5pt);
\draw [fill=ffqqqq] (24.205472704428672,9.192995786409236) circle (1.5pt);
\draw [fill=ffqqqq] (23.057472704428676,9.204995786409222) circle (1.5pt);
\draw [fill=qqqqff] (19.206961198706406,-2.189271594352926) circle (1.5pt);
\end{scriptsize}
\end{tikzpicture}
\end{center}

\begin{center}
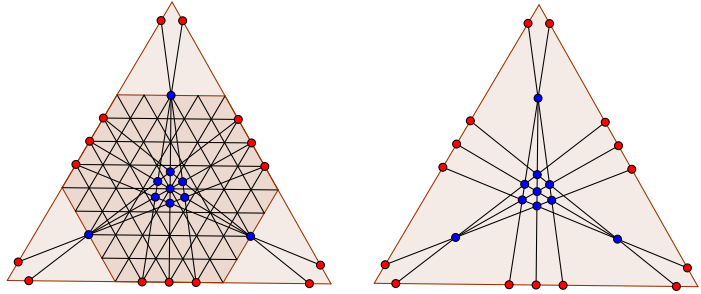
Figure 14. Pappus configuration being useful once again.  
\end{center}

The same trick can be used to obtain rotational symmetric Platonic configurations, for example, from the configuration in the Figure 15 we can get a tetrahedral 1-layer configuration with 4-valent points on edges, 3-valent points in the interiors of faces and 5-valent lines. Using radial projection and connecting 3-valent points in 5 copies of the tetrahedron we get a $(p_{4},n_{5})$ configuration.

\begin{center}
\definecolor{qqqqff}{rgb}{0.0,0.0,1.0}
\definecolor{ffqqqq}{rgb}{1.0,0.0,0.0}
\definecolor{uuuuuu}{rgb}{0.26666666666666666,0.26666666666666666,0.26666666666666666}
\definecolor{zzttqq}{rgb}{0.6,0.2,0.0}

\end{center}

\begin{center}
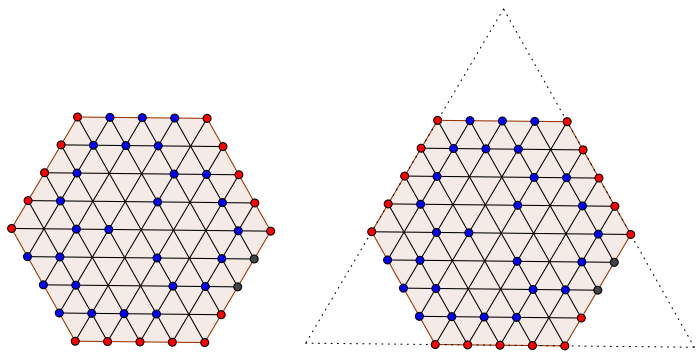
Figure 15. A $(p_{3},n_{5})$ configuration with 3-fold rotational symmetry.
\end{center}

Another interesting theme is shown in Example \ref{exampleseven}. Note how the variations of valences in edge points opens new possibilities in constructing Platonic configurations. 

\begin{example}\label{exampleseven}
Place the copies of the configuration $A = (12_{3}6_{2},12_{3})$ with 3-fold rotational symmetry (Figure 16) on each of the (triangulated) faces of any Platonic solid $P = (v,e,f,d,m)$ to get a Platonic 3-configuration $Z \in P_{R}$. The 3-valent "edge points" have valence 2 in one triangle and valence 1 in the adjacent triangle. If $P \in \lbrace T, O, I\rbrace$ then $y \in P_{R}(2e + 12f)_{3},(15f)_{3})$; thus we get configurations $Y$ in classes $T_{R}(60_{3})$, $O(120_{3})$ and $I(300_{3})$. If $P = C$ then $Z \in C_{R}((2e + 2mf + 12mf)_{3},15mf) = (360_{3})$.  If $P = D$ then $Z \in D_{R}((2e + 2mf + 12mf)_{3}, (15mf)_{3}) = D_{R}(900_{3})$. 
\end{example}

\begin{center}
\definecolor{qqqqff}{rgb}{0.0,0.0,1.0}
\definecolor{zzttqq}{rgb}{0.6,0.2,0.0}
\definecolor{ffqqqq}{rgb}{1.0,0.0,0.0}
\definecolor{ffffqq}{rgb}{1.0,1.0,0.0}
\definecolor{ffxfqq}{rgb}{1.0,0.4980392156862745,0.0}
\definecolor{qqzzqq}{rgb}{0.0,0.6,0.0}
\definecolor{qqffqq}{rgb}{0.0,1.0,0.0}
\definecolor{uuuuuu}{rgb}{0.26666666666666666,0.26666666666666666,0.26666666666666666}
\definecolor{xdxdff}{rgb}{0.49019607843137253,0.49019607843137253,1.0}
\begin{tikzpicture}[line cap=round,line join=round,>=triangle 45,x=0.3cm,y=0.3cm]
\clip(-18.764659340457058,8.571443560341388) rectangle (29.00608357926479,21.67213701530917);
\fill[color=zzttqq,fill=zzttqq,fill opacity=0.1] (-0.23631892977108748,13.385774275856496) -- (1.0848166265067734,13.250201395581563) -- (0.5416584067501612,14.462124789298544) -- cycle;
\fill[color=zzttqq,fill=zzttqq,fill opacity=0.1] (10.266464013475089,12.253033141857586) -- (12.96891924007515,13.341590798950623) -- (10.67497304224848,15.13770684922979) -- cycle;
\fill[color=zzttqq,fill=zzttqq,fill opacity=0.1] (15.39000851335711,11.212237439319527) -- (11.30651296015448,18.285059209948088) -- (7.223017406951823,11.212237439319514) -- cycle;
\draw (-14.717165954814133,-9.548093937233627)-- (-15.833737744050689,-6.83176302698025);
\draw (-13.666473437397878,-7.9141369322664366)-- (-14.717165954814133,-9.548093937233627);
\draw (-17.29441549218825,-6.102272178188465)-- (-15.833737744050689,-6.83176302698025);
\draw (-15.833737744050689,-6.83176302698025)-- (-13.666473437397878,-7.9141369322664366);
\draw (-16.586279982511577,-7.948378744252181)-- (-15.383591614983978,-7.926852349764274);
\draw (-16.586279982511577,-7.948378744252181)-- (-15.833737744050689,-6.83176302698025);
\draw (-14.968564735916427,15.372008068025778)-- (-11.518036270868523,14.872998777228817);
\draw (-14.223383105888763,13.374883803053251)-- (-10.104983718351882,15.655491861227912);
\draw (-13.529633547173201,9.91351643572761)-- (-14.646205336409754,12.629847345980988);
\draw (-13.232465318034244,15.1209364573882)-- (-11.613880200289943,12.887774201758388);
\draw (-12.134809045377708,13.60649984798208)-- (-9.182314252978014,14.089643000406397);
\draw (-11.156505062932034,11.557266200244975)-- (-9.182314252978014,14.089643000406397);
\draw (-14.196059207343053,11.53475802319696)-- (-12.134809045377708,13.60649984798208);
\draw (-11.615575694250003,12.890113475161284)-- (-12.47894102975695,11.5474734406948);
\draw (-12.47894102975695,11.5474734406948)-- (-13.529633547173201,9.91351643572761);
\draw (-11.346145881887148,9.932896503762693)-- (-11.156505062932034,11.557266200244975);
\draw (-11.156505062932034,11.557266200244975)-- (-10.893548972187975,13.809618345018242);
\draw (-12.134809045377708,13.60649984798208)-- (-11.615575694250003,12.890113475161284);
\draw (-12.134809045377708,13.60649984798208)-- (-12.709657039335461,13.028725905764135);
\draw (-10.252109637221219,12.71737189149552)-- (-11.615575694250003,12.890113475161284);
\draw (-11.615575694250003,12.890113475161284)-- (-12.709657039335461,13.028725905764135);
\draw (-16.106883084547317,13.35933819477277)-- (-14.646205336409754,12.629847345980988);
\draw (-14.646205336409754,12.629847345980988)-- (-12.47894102975695,11.5474734406948);
\draw (-14.968564735916427,15.372008068025778)-- (-14.223383105888763,13.374883803053251);
\draw [dotted] (-18.07859351155319,9.873141051428373)-- (-6.7569712358693685,9.973628823106038);
\draw [dotted] (-6.7569712358693685,9.973628823106038)-- (-12.504807336753831,19.728197440061177);
\draw [dotted] (-12.504807336753831,19.728197440061177)-- (-18.07859351155319,9.873141051428373);
\draw (-14.196059207343053,11.53475802319696)-- (-11.156505062932034,11.557266200244975);
\draw (-10.252109637221219,12.71737189149552)-- (-11.518036270868523,14.872998777228817);
\draw (-13.232465318034244,15.1209364573882)-- (-14.646205336409754,12.629847345980988);
\draw (-15.398747574870631,11.513231628709063)-- (-14.196059207343053,11.53475802319696);
\draw (-15.398747574870631,11.513231628709063)-- (-14.646205336409754,12.629847345980988);
\draw (-11.156505062932034,11.557266200244975)-- (-9.734194404118165,11.581204614849288);
\draw (-10.252109637221219,12.71737189149552)-- (-9.734194404118165,11.581204614849288);
\draw (-13.232465318034244,15.1209364573882)-- (-12.41534862767361,16.737270429378977);
\draw (-11.518036270868523,14.872998777228817)-- (-12.41534862767361,16.737270429378977);
\draw (-2.392046446806132,14.21066699221791)-- (0.4507067276222281,12.410696432803816);
\draw (1.4793883405735793,11.013848483529747)-- (1.6168319838760397,14.375730228980714);
\draw (2.3121969866959606,15.965018612839472)-- (-0.6679998310335904,14.403107426804986);
\draw (-1.3791710663870163,15.965018612839472)-- (1.0739375449667967,15.316050980318124);
\draw (-0.5463624202646373,11.013848483529747)-- (-1.2108942699439231,13.462786675465694);
\draw (1.5364956054402221,12.410696432803816)-- (3.325072367115073,14.210666992217911);
\draw(0.46651296015446964,13.729844696195707) circle (0.7825325038927978cm);
\draw [color=zzttqq] (-0.23631892977108748,13.385774275856496)-- (1.0848166265067734,13.250201395581563);
\draw [color=zzttqq] (1.0848166265067734,13.250201395581563)-- (0.5416584067501612,14.462124789298544);
\draw [color=zzttqq] (0.5416584067501612,14.462124789298544)-- (-0.23631892977108748,13.385774275856496);
\draw (8.447953553193877,14.050666992217927)-- (10.676984622908348,11.212237439319525);
\draw (12.319388340573578,10.853848483529747)-- (13.66302490537884,14.203460792176504);
\draw (13.152196986695959,15.805018612839474)-- (9.579529352176241,15.293835857091182);
\draw (9.460828933612971,15.805018612839453)-- (12.769942413647545,15.750325043205322);
\draw (10.293637579735357,10.853848483529747)-- (8.686446860444912,13.746971606062306);
\draw (12.463149606370957,11.21223743931953)-- (14.165072367115075,14.05066699221791);
\draw [color=zzttqq] (10.266464013475089,12.253033141857586)-- (12.96891924007515,13.341590798950623);
\draw [color=zzttqq] (12.96891924007515,13.341590798950623)-- (10.67497304224848,15.13770684922979);
\draw [color=zzttqq] (10.67497304224848,15.13770684922979)-- (10.266464013475089,12.253033141857586);
\draw [dotted] (6.602269526652369,10.853848483529745)-- (16.010756393656557,10.853848483529744);
\draw [dotted] (16.010756393656557,10.853848483529744)-- (11.30651296015447,19.001837121527636);
\draw [dotted] (11.30651296015447,19.001837121527636)-- (6.602269526652369,10.853848483529745);
\draw (10.094544980061512,11.212237439319527)-- (10.67497304224848,15.13770684922979);
\draw (10.266464013475089,12.253033141857586)-- (13.954244726802271,13.699053265299984);
\draw (12.96891924007515,13.341590798950623)-- (9.870749173599618,15.798243383967602);
\draw (10.67497304224848,15.13770684922979)-- (14.165072367115075,14.05066699221791);
\draw (10.266464013475089,12.253033141857586)-- (9.460828933612971,15.805018612839453);
\draw (10.293637579735357,10.853848483529747)-- (12.96891924007515,13.341590798950623);
\draw (-0.23631892977108748,13.385774275856496)-- (-0.9254042247938115,12.410696432803816);
\draw (1.0848166265067734,13.250201395581563)-- (2.304887460084114,13.183983185762335);
\draw (0.5416584067501612,14.462124789298544)-- (0.02005564517323702,15.594854470021208);
\draw (-0.9254042247938115,12.410696432803816)-- (1.5364956054402221,12.410696432803816);
\draw (2.304887460084114,13.183983185762335)-- (1.0739375449667967,15.316050980318124);
\draw (0.02005564517323702,15.594854470021208)-- (-1.2108942699439231,13.462786675465694);
\draw (0.4507067276222281,12.410696432803816)-- (-0.5463624202646373,11.013848483529747);
\draw (0.4507067276222281,12.410696432803816)-- (1.0848166265067734,13.250201395581563);
\draw (3.325072367115073,14.210666992217911)-- (1.6168319838760397,14.375730228980714);
\draw (1.6168319838760397,14.375730228980714)-- (0.5416584067501612,14.462124789298544);
\draw (-1.3791710663870163,15.965018612839472)-- (-0.6679998310335904,14.403107426804986);
\draw (-0.6679998310335904,14.403107426804986)-- (-0.23631892977108748,13.385774275856496);
\draw [dotted] (0.4665129601544773,19.161837121527636)-- (5.170756393656562,11.013848483529745);
\draw [dotted] (-4.2377304733476215,11.013848483529747)-- (5.170756393656562,11.013848483529745);
\draw [dotted] (-4.2377304733476215,11.013848483529747)-- (0.4665129601544773,19.161837121527636);
\draw (9.870749173599618,15.798243383967602)-- (8.686446860444912,13.746971606062306);
\draw (10.094544980061512,11.212237439319527)-- (12.463149606370957,11.21223743931953);
\draw (13.954244726802271,13.699053265299984)-- (12.769942413647545,15.750325043205322);
\draw (20.10856953007282,14.059187682797484)-- (22.337600599787326,11.22075812989911);
\draw (23.98000431745254,10.862369174109332)-- (25.323640882257767,14.211981482756205);
\draw (24.81281296357493,15.813539303419057)-- (21.24014532905524,15.3023565476708);
\draw (21.121444910491952,15.813539303419073)-- (24.430558390526542,15.758845733784911);
\draw (21.95425355661428,10.862369174109334)-- (20.34706283732387,13.755492296641858);
\draw (24.123765583249888,11.220758129899105)-- (25.825688343994045,14.059187682797498);
\draw [dotted] (18.26288550353132,10.862369174109336)-- (27.67137237053551,10.862369174109325);
\draw [dotted] (27.67137237053551,10.862369174109325)-- (22.967128937033415,19.01035781210722);
\draw [dotted] (22.967128937033415,19.01035781210722)-- (18.26288550353132,10.862369174109336);
\draw (21.75516095694043,11.220758129899112)-- (22.335589019127443,15.146227539809372);
\draw (21.927079990354052,12.261553832437176)-- (25.61486070368125,13.70757395587959);
\draw (24.629535216954107,13.35011148953021)-- (21.53136515047862,15.806764074547223);
\draw (22.335589019127443,15.146227539809372)-- (25.825688343994045,14.059187682797498);
\draw (21.927079990354052,12.261553832437176)-- (21.121444910491952,15.813539303419073);
\draw (21.95425355661428,10.862369174109334)-- (24.629535216954107,13.35011148953021);
\draw (21.53136515047862,15.806764074547223)-- (20.34706283732387,13.755492296641858);
\draw (21.75516095694043,11.220758129899112)-- (24.123765583249888,11.220758129899105);
\draw (25.61486070368125,13.70757395587959)-- (24.430558390526542,15.758845733784911);
\draw [color=zzttqq] (15.39000851335711,11.212237439319527)-- (11.30651296015448,18.285059209948088);
\draw [color=zzttqq] (11.30651296015448,18.285059209948088)-- (7.223017406951823,11.212237439319514);
\draw [color=zzttqq] (7.223017406951823,11.212237439319514)-- (15.39000851335711,11.212237439319527);
\begin{scriptsize}
\draw [fill=xdxdff] (-17.29441549218825,-6.102272178188465) circle (1.5pt);
\draw [fill=xdxdff] (-14.717165954814133,-9.548093937233627) circle (1.5pt);
\draw [fill=xdxdff] (-12.533678289528075,-9.528713869198544) circle (1.5pt);
\draw [fill=xdxdff] (-10.369846660618943,-5.3719673725548365) circle (1.5pt);
\draw [fill=xdxdff] (-15.833737744050689,-6.83176302698025) circle (1.5pt);
\draw [fill=uuuuuu] (-15.383591614983978,-7.926852349764274) circle (1.5pt);
\draw [fill=uuuuuu] (-13.666473437397878,-7.9141369322664366) circle (1.5pt);
\draw [fill=qqffqq] (-16.106883084547317,13.35933819477277) circle (1.5pt);
\draw [fill=qqzzqq] (-14.968564735916427,15.372008068025778) circle (1.5pt);
\draw [fill=qqzzqq] (-13.529633547173201,9.91351643572761) circle (1.5pt);
\draw [fill=qqffqq] (-11.346145881887148,9.932896503762693) circle (1.5pt);
\draw [fill=qqffqq] (-10.104983718351882,15.655491861227912) circle (1.5pt);
\draw [fill=qqzzqq] (-9.182314252978014,14.089643000406397) circle (1.5pt);
\draw [fill=ffxfqq] (-11.518036270868523,14.872998777228817) circle (1.5pt);
\draw [fill=qqffqq] (-14.223383105888763,13.374883803053251) circle (1.5pt);
\draw [fill=ffffqq] (-13.232465318034244,15.1209364573882) circle (1.5pt);
\draw [fill=ffxfqq] (-14.646205336409754,12.629847345980988) circle (1.5pt);
\draw [fill=ffffqq] (-14.196059207343053,11.53475802319696) circle (1.5pt);
\draw [fill=qqffqq] (-12.47894102975695,11.5474734406948) circle (1.5pt);
\draw [fill=ffqqqq] (-12.134809045377708,13.60649984798208) circle (1.5pt);
\draw [fill=qqffqq] (-10.893548972187975,13.809618345018242) circle (1.5pt);
\draw [fill=ffqqqq] (-11.615575694250003,12.890113475161284) circle (1.5pt);
\draw [fill=ffxfqq] (-11.156505062932034,11.557266200244975) circle (1.5pt);
\draw [fill=ffffqq] (-10.252109637221219,12.71737189149552) circle (1.5pt);
\draw [fill=ffqqqq] (-12.709657039335461,13.028725905764135) circle (1.5pt);
\draw [fill=qqffqq] (1.4793883405735793,11.013848483529747) circle (1.5pt);
\draw [fill=qqffqq] (1.4793883405735793,11.013848483529747) circle (1.5pt);
\draw [fill=qqzzqq] (-0.5463624202646373,11.013848483529747) circle (1.5pt);
\draw [fill=qqzzqq] (3.325072367115073,14.210666992217911) circle (1.5pt);
\draw [fill=qqffqq] (2.3121969866959606,15.965018612839472) circle (1.5pt);
\draw [fill=qqzzqq] (-1.3791710663870163,15.965018612839472) circle (1.5pt);
\draw [fill=qqffqq] (-2.392046446806132,14.21066699221791) circle (1.5pt);
\draw [fill=qqffqq] (0.4507067276222281,12.410696432803816) circle (1.5pt);
\draw [fill=ffxfqq] (-1.2108942699439231,13.462786675465694) circle (1.5pt);
\draw [fill=qqffqq] (-0.6679998310335904,14.403107426804986) circle (1.5pt);
\draw [fill=ffxfqq] (1.0739375449667967,15.316050980318124) circle (1.5pt);
\draw [fill=ffxfqq] (1.5364956054402221,12.410696432803816) circle (1.5pt);
\draw [fill=ffffqq] (0.02005564517323702,15.594854470021208) circle (1.5pt);
\draw [fill=ffffqq] (-0.9254042247938115,12.410696432803816) circle (1.5pt);
\draw [fill=ffffqq] (2.304887460084114,13.183983185762335) circle (1.5pt);
\draw [fill=ffqqqq] (-0.23631892977108748,13.385774275856496) circle (1.5pt);
\draw [fill=ffqqqq] (1.0848166265067734,13.250201395581563) circle (1.5pt);
\draw [fill=ffqqqq] (0.5416584067501612,14.462124789298544) circle (1.5pt);
\draw [fill=qqffqq] (12.319388340573578,10.853848483529747) circle (1.5pt);
\draw [fill=qqffqq] (12.319388340573578,10.853848483529747) circle (1.5pt);
\draw [fill=qqzzqq] (10.293637579735357,10.853848483529747) circle (1.5pt);
\draw [fill=qqzzqq] (14.165072367115075,14.05066699221791) circle (1.5pt);
\draw [fill=qqffqq] (13.152196986695959,15.805018612839474) circle (1.5pt);
\draw [fill=qqzzqq] (9.460828933612971,15.805018612839453) circle (1.5pt);
\draw [fill=qqffqq] (8.447953553193877,14.050666992217927) circle (1.5pt);
\draw [fill=qqffqq] (10.676984622908348,11.212237439319525) circle (1.5pt);
\draw [fill=ffxfqq] (8.686446860444912,13.746971606062306) circle (1.5pt);
\draw [fill=qqffqq] (13.66302490537884,14.203460792176504) circle (1.5pt);
\draw [fill=qqffqq] (9.579529352176241,15.293835857091182) circle (1.5pt);
\draw [fill=ffxfqq] (12.769942413647545,15.750325043205322) circle (1.5pt);
\draw [fill=ffxfqq] (12.463149606370957,11.21223743931953) circle (1.5pt);
\draw [fill=ffffqq] (9.870749173599618,15.798243383967602) circle (1.5pt);
\draw [fill=ffffqq] (10.094544980061512,11.212237439319527) circle (1.5pt);
\draw [fill=ffffqq] (13.954244726802271,13.699053265299984) circle (1.5pt);
\draw [fill=ffqqqq] (12.967893713817443,13.334242335245024) circle (1.5pt);
\draw [fill=ffqqqq] (10.266464013475089,12.253033141857586) circle (1.5pt);
\draw [fill=ffqqqq] (12.96891924007515,13.341590798950623) circle (1.5pt);
\draw [fill=ffqqqq] (10.67497304224848,15.13770684922979) circle (1.5pt);
\draw [fill=qqffqq] (23.98000431745254,10.862369174109332) circle (1.5pt);
\draw [fill=qqffqq] (23.98000431745254,10.862369174109332) circle (1.5pt);
\draw [fill=qqzzqq] (21.95425355661428,10.862369174109334) circle (1.5pt);
\draw [fill=qqzzqq] (25.825688343994045,14.059187682797498) circle (1.5pt);
\draw [fill=qqffqq] (24.81281296357493,15.813539303419057) circle (1.5pt);
\draw [fill=qqzzqq] (21.121444910491952,15.813539303419073) circle (1.5pt);
\draw [fill=qqffqq] (20.10856953007282,14.059187682797484) circle (1.5pt);
\draw [fill=qqffqq] (22.337600599787326,11.22075812989911) circle (1.5pt);
\draw [fill=ffxfqq] (20.34706283732387,13.755492296641858) circle (1.5pt);
\draw [fill=qqffqq] (25.323640882257767,14.211981482756205) circle (1.5pt);
\draw [fill=qqffqq] (21.24014532905524,15.3023565476708) circle (1.5pt);
\draw [fill=ffxfqq] (24.430558390526542,15.758845733784911) circle (1.5pt);
\draw [fill=ffxfqq] (24.123765583249888,11.220758129899105) circle (1.5pt);
\draw [fill=ffffqq] (21.53136515047862,15.806764074547223) circle (1.5pt);
\draw [fill=ffffqq] (21.75516095694043,11.220758129899112) circle (1.5pt);
\draw [fill=ffffqq] (25.61486070368125,13.70757395587959) circle (1.5pt);
\draw [fill=qqqqff] (24.6285096906964,13.342763025824608) circle (1.5pt);
\draw [fill=ffqqqq] (21.927079990354052,12.261553832437176) circle (1.5pt);
\draw [fill=ffqqqq] (24.629535216954107,13.35011148953021) circle (1.5pt);
\draw [fill=ffqqqq] (22.335589019127443,15.146227539809372) circle (1.5pt);
\draw [fill=qqffqq] (1.6168319838760405,14.375730228980716) circle (1.5pt);
\end{scriptsize}
\end{tikzpicture}
\end{center}

\begin{center}
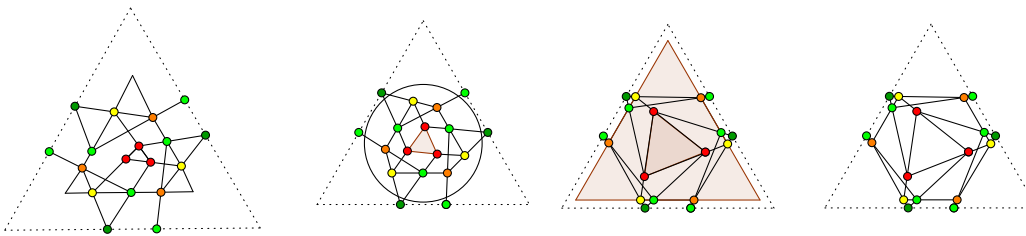
Figure 16.  Two-valent and one-valent vertices of $A$ on the edges may produce 3-valent vertices in $P$.
\end{center}

We just briefly mention another important construction; using the axes through the vertices (of degree $m \in \lbrace 3,4,5\rbrace$) of the chosen Platonic solid $P$ (see Figure 17) as the "meeting places" of 1-valent points we may obtain examples of "spider-web configurations". 
\begin{example}\label{osi}
Take any two "adjacent" axes. Place in the plane spanned by these axes a copy of symmetric configuration $A$ with $m$-valent lines, and with $m$-valent points and 1-valent points; placing the 1-valent points on the axes these points become $m$-valent and we get a balanced $(n_{m})$ configuration. 
\end{example}

\begin{center}\definecolor{qqqqff}{rgb}{0.0,0.0,1.0}
\definecolor{uuuuuu}{rgb}{0.26666666666666666,0.26666666666666666,0.26666666666666666}
\definecolor{xfqqff}{rgb}{0.4980392156862745,0.0,1.0}
\definecolor{qqqqff}{rgb}{0.0,0.0,1.0}
\definecolor{ffxfqq}{rgb}{1.0,0.4980392156862745,0.0}
\definecolor{qqffqq}{rgb}{0.0,1.0,0.0}
\begin{tikzpicture}[line cap=round,line join=round,>=triangle 45,x=0.4cm,y=0.4cm]
\clip(-8.974466967713285,0.5077959837588182) rectangle (19.17049890721249,6.223771393235991);
\draw [color=qqffqq] (-7.8433852801,5.055058319095684)-- (-6.403385280099999,0.8550583190956851);
\draw (-7.361281554753032,3.648922453500363)-- (-5.730543789961213,2.8175126653338123);
\draw (-7.078330495585757,2.823648530929142)-- (-5.435320999603288,3.6785791372110945);
\draw (-6.417266861068892,3.1676292171214047)-- (-5.607475857163551,3.176460802660328);
\draw (-6.417266861068892,3.1676292171214047)-- (-7.188951386747188,3.1462927968166525);
\draw (-3.4433852800999993,1.575058319095685)-- (-0.3833852800999993,1.575058319095685);
\draw (-0.3833852800999993,1.575058319095685)-- (-2.0233852800999994,4.375058319095684);
\draw (-2.0233852800999994,4.375058319095684)-- (-3.4433852800999993,1.575058319095685);
\draw [color=qqffqq] (-1.8233852800999995,2.5550583190956853)-- (-2.1633852800999995,5.795058319095684);
\draw [color=qqffqq] (-1.8233852800999995,2.5550583190956853)-- (-4.7233852801,0.7350583190956851);
\draw [color=qqffqq] (-1.8233852800999995,2.5550583190956853)-- (1.1566147199000003,0.575058319095685);
\draw (-2.0233852800999994,4.375058319095684)-- (-0.6200407646670156,2.8827776764902);
\draw (-0.6200407646670156,2.8827776764902)-- (-0.3833852800999993,1.575058319095685);
\draw [color=qqffqq] (-1.8233852800999995,2.5550583190956857)-- (2.376614719900001,3.735058319095685);
\draw [color=qqffqq] (-6.403385280099999,0.8550583190956846)-- (-5.0033852801,5.035058319095684);
\draw (1.284178663490929,3.4281358080093267)-- (-0.021529146108372643,1.3578518944972215);
\draw (0.3482471178751265,1.1121616251390576)-- (0.6255922529489071,3.2431043879046637);
\draw (0.12571625202417014,3.1026630352639044)-- (1.0657462149926313,0.6354340371213866);
\draw (4.176614719900001,1.6750583190956851)-- (7.2366147199,1.6750583190956851);
\draw (7.2366147199,1.6750583190956851)-- (5.5966147199,4.475058319095684);
\draw (5.5966147199,4.475058319095684)-- (4.176614719900001,1.6750583190956851);
\draw [color=qqffqq] (5.796614719900001,2.655058319095685)-- (5.456614719900001,5.895058319095684);
\draw [color=qqffqq] (5.796614719900001,2.655058319095685)-- (2.8966147199000005,0.8350583190956851);
\draw (5.796614719900001,2.655058319095685)-- (8.776614719900001,0.6750583190956851);
\draw (5.5966147199,4.475058319095684)-- (6.999959235332986,2.9827776764902);
\draw (6.999959235332986,2.9827776764902)-- (7.2366147199,1.6750583190956851);
\draw [color=qqffqq] (5.7966147199,2.655058319095685)-- (9.9966147199,3.835058319095685);
\draw (8.904178663490928,3.5281358080093264)-- (7.598470853891627,1.4578518944972225);
\draw (7.968247117875128,1.2121616251390572)-- (8.245592252948907,3.3431043879046634);
\draw (7.74571625202417,3.202663035263904)-- (8.685746214992639,0.7354340371213808);
\draw (3.794441103575697,1.3985217736783635)-- (8.685746214992639,0.7354340371213808);
\draw (7.598470853891627,1.4578518944972225)-- (2.9953726534289706,0.8970374360000732);
\draw (3.3920194047565975,1.1459674661436183)-- (7.968247117875128,1.2121616251390572);
\draw [dotted] (-3.4433852801,1.575058319095685)-- (-0.6200407646670156,2.8827776764902);
\draw [dotted] (4.1766147199,1.675058319095685)-- (6.999959235332988,2.982777676490196);
\draw [color=qqffqq] (14.06017268506215,2.7041871892727776)-- (13.720172685062154,5.944187189272778);
\draw [color=qqffqq] (14.06017268506215,2.7041871892727776)-- (11.160172685062152,0.8841871892727755);
\draw [color=qqffqq] (14.06017268506215,2.7041871892727776)-- (17.040172685062167,0.7241871892727754);
\draw [color=qqffqq] (14.060172685062149,2.704187189272777)-- (18.260172685062162,3.884187189272779);
\draw (17.405071163992897,3.643944380972366)-- (15.862028819053787,1.5069807646743132);
\draw (16.23180508303729,1.2612904953161475)-- (16.50915021811106,3.392233258081756);
\draw (16.009274217186327,3.2517919054409976)-- (16.762884900002167,0.9084253820307623);
\draw (12.057999068737852,1.4476506438554568)-- (16.762884900002167,0.9084253820307623);
\draw (15.862028819053787,1.5069807646743132)-- (11.160172685062152,0.8841871892727755);
\draw (11.655577369918747,1.195096336320708)-- (16.23180508303729,1.2612904953161475);
\draw (13.735728383640332,5.795950532233668)-- (12.057999068737852,1.4476506438554568);
\draw (13.82241100636786,4.969916127418394)-- (11.160172685062152,0.8841871892727755);
\draw (13.780916835788437,5.365331164704669)-- (11.655577369918747,1.195096336320708);
\draw (13.82241100636786,4.969916127418394)-- (17.405071163992897,3.643944380972366);
\draw (16.009274217186327,3.2517919054409976)-- (13.735728383640332,5.795950532233668);
\draw (16.50915021811106,3.392233258081756)-- (13.780916835788437,5.365331164704669);
\draw (13.82241100636786,4.969916127418394)-- (16.762884900002167,0.9084253820307623);
\draw (15.862028819053787,1.5069807646743132)-- (13.735728383640332,5.795950532233668);
\draw (16.23180508303729,1.2612904953161475)-- (13.780916835788437,5.365331164704669);
\draw (12.057999068737852,1.4476506438554568)-- (17.405071163992897,3.643944380972366);
\draw (16.009274217186327,3.2517919054409976)-- (11.160172685062152,0.8841871892727755);
\draw (-8.519107992134925,-3.4511344154408325)-- (-4.99467526362131,-1.6743542795786024);
\draw (-6.509307510585839,-1.5869716499460336)-- (-3.8587010783979134,-3.888047563603676);
\draw (-8.16957747360465,-2.8103284648019953)-- (-3.8295735351870572,-2.8685835512237077);
\draw (9.306948452909143,-0.8879106128854841)-- (5.025199600913265,-3.655027217916826);
\draw (5.025199600913265,-3.655027217916826)-- (10.239029835656545,-3.3346242425974073);
\draw (5.025199600913265,-3.655027217916826)-- (10.363157614056579,-2.1274376805396042);
\draw (7.268020428149202,-0.6378526140289199)-- (5.869898354028098,-4.674491230296794);
\draw (7.268020428149202,-0.6378526140289199)-- (6.901513929411405,-4.693828688121362);
\draw (7.268020428149202,-0.6378526140289199)-- (8.234445693039145,-4.753512199925589);
\draw (6.354784992965808,-3.274533267056196)-- (7.00638367893577,-3.533278475691978);
\draw (7.00638367893577,-3.533278475691978)-- (9.96526753536174,-4.8131957117298185);
\draw (9.96526753536174,-4.8131957117298185)-- (7.1175709902943955,-2.302814415255554);
\draw (7.1175709902943955,-2.302814415255554)-- (5.409426134305727,-0.7746114129771279);
\draw (6.354784992965808,-3.274533267056196)-- (4.553962465111805,-2.3462772238217697);
\draw (9.96526753536174,-4.8131957117298185)-- (5.01153605561088,-1.7096530979100162);
\draw (16.50915021811106,3.392233258081756)-- (11.655577369918747,1.195096336320708);
\draw (14.0974737656347,-2.6419428370232705)-- (17.706845663885613,-5.11986216283387);
\begin{scriptsize}
\draw [fill=ffxfqq] (-7.361281554753032,3.648922453500363) circle (1.5pt);
\draw [fill=ffxfqq] (-5.730543789961213,2.8175126653338123) circle (1.5pt);
\draw [fill=ffxfqq] (-7.078330495585757,2.823648530929142) circle (1.5pt);
\draw [fill=ffxfqq] (-5.435320999603288,3.6785791372110945) circle (1.5pt);
\draw [fill=qqqqff] (-6.417266861068892,3.1676292171214047) circle (1.5pt);
\draw [fill=ffxfqq] (-5.607475857163551,3.176460802660328) circle (1.5pt);
\draw [fill=ffxfqq] (-7.188951386747188,3.1462927968166525) circle (1.5pt);
\draw [fill=ffxfqq] (0.12571625202417014,3.1026630352639044) circle (1.5pt);
\draw [fill=ffxfqq] (0.6255922529489071,3.2431043879046637) circle (1.5pt);
\draw [fill=ffxfqq] (1.284178663490929,3.4281358080093267) circle (1.5pt);
\draw [fill=ffxfqq] (-0.021529146108372643,1.3578518944972215) circle (1.5pt);
\draw [fill=ffxfqq] (0.3482471178751265,1.1121616251390576) circle (1.5pt);
\draw [fill=xfqqff] (0.48468892726752666,2.1604932866482205) circle (1.5pt);
\draw [fill=ffxfqq] (1.0657462149926313,0.6354340371213866) circle (1.5pt);
\draw [fill=ffxfqq] (7.74571625202417,3.202663035263904) circle (1.5pt);
\draw [fill=ffxfqq] (8.245592252948907,3.3431043879046634) circle (1.5pt);
\draw [fill=ffxfqq] (8.904178663490928,3.5281358080093264) circle (1.5pt);
\draw [fill=ffxfqq] (7.598470853891627,1.4578518944972225) circle (1.5pt);
\draw [fill=ffxfqq] (7.968247117875128,1.2121616251390572) circle (1.5pt);
\draw [fill=xfqqff] (8.104688927267528,2.260493286648225) circle (1.5pt);
\draw [fill=ffxfqq] (8.685746214992639,0.7354340371213808) circle (1.5pt);
\draw [fill=ffxfqq] (3.794441103575697,1.3985217736783635) circle (1.5pt);
\draw [fill=ffxfqq] (3.3920194047565975,1.1459674661436183) circle (1.5pt);
\draw [fill=ffxfqq] (2.9953726534289706,0.8970374360000732) circle (1.5pt);
\draw [fill=ffxfqq] (5.558853041205708,4.920787257241299) circle (1.5pt);
\draw [fill=ffxfqq] (5.472170418478179,5.746821662056574) circle (1.5pt);
\draw [fill=ffxfqq] (5.517358870626284,5.316202294527574) circle (1.5pt);
\draw [fill=xfqqff] (5.260382529638784,1.1729929407272217) circle (1.5pt);
\draw [fill=ffxfqq] (16.009274217186327,3.2517919054409976) circle (1.5pt);
\draw [fill=ffxfqq] (16.50915021811106,3.392233258081756) circle (1.5pt);
\draw [fill=ffxfqq] (17.405071163992897,3.643944380972366) circle (1.5pt);
\draw [fill=ffxfqq] (15.862028819053787,1.5069807646743132) circle (1.5pt);
\draw [fill=ffxfqq] (16.23180508303729,1.2612904953161475) circle (1.5pt);
\draw [fill=xfqqff] (16.352119361104915,2.1857085134558174) circle (1.5pt);
\draw [fill=ffxfqq] (16.762884900002167,0.9084253820307623) circle (1.5pt);
\draw [fill=ffxfqq] (12.057999068737852,1.4476506438554568) circle (1.5pt);
\draw [fill=ffxfqq] (11.655577369918747,1.195096336320708) circle (1.5pt);
\draw [fill=ffxfqq] (11.160172685062152,0.8841871892727755) circle (1.5pt);
\draw [fill=ffxfqq] (13.82241100636786,4.969916127418394) circle (1.5pt);
\draw [fill=ffxfqq] (13.735728383640332,5.795950532233668) circle (1.5pt);
\draw [fill=ffxfqq] (13.780916835788437,5.365331164704669) circle (1.5pt);
\draw [fill=xfqqff] (13.734438528179359,1.2251666201146871) circle (1.5pt);
\draw [fill=xfqqff] (12.828449769865843,3.4444868843371363) circle (1.5pt);
\draw [fill=ffxfqq] (11.655577369918749,1.1950963363207086) circle (1.5pt);
\draw [fill=xfqqff] (14.79586022697563,4.609634630894462) circle (1.5pt);
\draw [fill=ffxfqq] (13.82241100636786,4.969916127418392) circle (1.5pt);
\draw [fill=ffxfqq] (16.76288490000217,0.9084253820307607) circle (1.5pt);
\draw [fill=ffxfqq] (13.780916835788435,5.36533116470467) circle (1.5pt);
\draw [fill=xfqqff] (13.671919756803161,2.1105636991229826) circle (1.5pt);
\draw [fill=xfqqff] (14.84650250852443,3.555403260140974) circle (1.5pt);
\draw [fill=qqqqff] (-8.519107992134925,-3.4511344154408325) circle (1.5pt);
\draw [fill=qqqqff] (-4.99467526362131,-1.6743542795786024) circle (1.5pt);
\draw [fill=qqqqff] (-6.509307510585839,-1.5869716499460336) circle (1.5pt);
\draw [fill=qqqqff] (-3.8587010783979134,-3.888047563603676) circle (1.5pt);
\draw [fill=qqqqff] (-8.16957747360465,-2.8103284648019953) circle (1.5pt);
\draw [fill=qqqqff] (-3.8295735351870572,-2.8685835512237077) circle (1.5pt);
\draw [fill=qqqqff] (7.268020428149202,-0.6378526140289199) circle (1.5pt);
\draw [fill=uuuuuu] (6.5684295193754005,-2.6577017604072837) circle (1.5pt);
\draw [fill=uuuuuu] (7.2680204281492005,-0.6378526140289226) circle (1.5pt);
\draw [fill=uuuuuu] (7.2680204281492005,-0.6378526140289226) circle (1.5pt);
\draw [fill=uuuuuu] (7.2680204281492005,-0.6378526140289226) circle (1.5pt);
\draw [fill=uuuuuu] (7.2680204281492005,-0.6378526140289226) circle (1.5pt);
\draw [fill=uuuuuu] (7.1175709902943955,-2.302814415255554) circle (1.5pt);
\draw [fill=uuuuuu] (6.354784992965808,-3.274533267056196) circle (1.5pt);
\draw [fill=uuuuuu] (7.00638367893577,-3.533278475691978) circle (1.5pt);
\draw [fill=qqqqff] (9.96526753536174,-4.8131957117298185) circle (1.5pt);
\draw [fill=uuuuuu] (7.800309439630175,-2.904681154380263) circle (1.5pt);
\draw [fill=uuuuuu] (9.96526753536174,-4.813195711729819) circle (1.5pt);
\draw [fill=uuuuuu] (7.952328861541782,-3.552077506445031) circle (1.5pt);
\draw [fill=uuuuuu] (5.025199600913265,-3.655027217916825) circle (1.5pt);
\draw [fill=qqqqff] (14.0974737656347,-2.6419428370232705) circle (1.5pt);
\draw [fill=qqqqff] (17.706845663885613,-5.11986216283387) circle (1.5pt);
\end{scriptsize}
\end{tikzpicture}
\end{center}
\vspace{-0mm}

\begin{center}
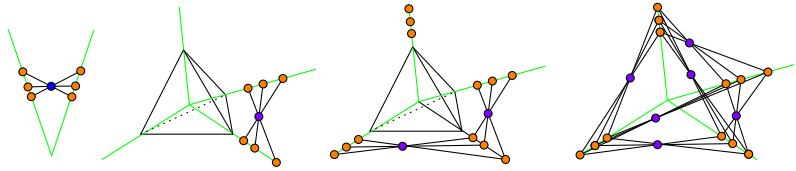
Figure 17. Using the 3-fold axes of tetrahedron as "linking places" to construct Platonic configuration $(24_{3})$ 
\end{center}

Using such a "spider-web" construction we can obtain examples of $(n_{5})$ configurations with full symmetry group of icosahedron. We start with a planar $(n_{5})$ configuration $A$ with dihedral symmetry (whose existence is guaranteed by  \cite{BermBoko}) and "split" two $5$-valent points in every copy of $A$ into $5$ 1-valent points (similarly as we replaced in the solution of Problem 3 some 3-valent points of the Pappus configuration with 3 1-valent points). Now we place these 1-valent points on the pair of axes going through the adjacent vertices of the icosahedron; thus their valence is increased to 5. We do this for all pairs of adjacent axes. We include these axes among the 5-valent lines of configuration. Now all the lines and all the points have valence 5, and the obtained spider-web configuration is balanced.

The last example illustrates the concept of a "helical" configuration. This construction works for any Platonic solid $P$ and produces a balanced 3-configuration with the full symmetry group of $P$.
\begin{example} Take a Platonic solid $P$. If the degree of its vertices $d$ is bigger than 3, truncate $P$ in each vertex; so you obtain the solid $P'$ with 3-valent vertices.
Place a copy of the planar configuration $ A = (1_{3}6_{1},3_{3})$ on each edge of $P'$ so that the 3-valent point is in the center of the edge. 
All the 1-valent points must be turned into 3-valent points. This is obtained by connecting them to two other points around the chosen vertex and adding three new 3-valent points around each vertex, as it is shown in Figure 18. So we see that the whole configuration is made of copies of $A$ -- each one is glued to two other such copies.

\begin{center}
\definecolor{ffxfqq}{rgb}{1.0,0.4980392156862745,0.0}
\definecolor{xfqqff}{rgb}{0.4980392156862745,0.0,1.0}
\definecolor{xdxdff}{rgb}{0.49019607843137253,0.49019607843137253,1.0}
\definecolor{ffqqqq}{rgb}{1.0,0.0,0.0}
\definecolor{uuuuuu}{rgb}{0.26666666666666666,0.26666666666666666,0.26666666666666666}
\definecolor{qqqqff}{rgb}{0.0,0.0,1.0}
\begin{tikzpicture}[line cap=round,line join=round,>=triangle 45,x=1.0cm,y=1.0cm]
\clip(3.053181823624341,0.8060172122750571) rectangle (9.107339734539185,5.038570088047371);
\draw (-2.8,0.74)-- (-2.82,3.3);
\draw (-3.62,0.76)-- (-1.94,3.3);
\draw (-3.66,3.32)-- (-1.92,0.74);
\draw (11.62,3.46)-- (11.66,0.84);
\draw (11.62,3.46)-- (9.96,5.0);
\draw (11.62,3.46)-- (12.84,5.0);
\draw (11.66,0.84)-- (9.98,-0.96);
\draw (11.66,0.84)-- (13.04,-1.06);
\draw (11.368524191785584,2.808524191785584)-- (12.048524191785583,1.7485241917855843);
\draw (11.281200688245544,1.8083085579690314)-- (11.841215286604793,2.8038571769057428);
\draw (10.654868219158892,4.781102289315581)-- (11.011420616478699,3.597348330213822);
\draw (11.24,4.24)-- (10.212743246482331,4.15357007003272);
\draw (12.34,4.74)-- (12.152388287902083,3.654396713784991);
\draw (11.94,4.22)-- (12.82,4.4);
\draw (10.28,-0.08)-- (11.38,-0.06);
\draw (10.84,0.34)-- (10.8,-0.66);
\draw (12.0,0.04)-- (12.76,-0.3);
\draw (12.467304651784184,0.21753383039956012)-- (12.26,-0.52);
\draw (11.111173349976594,3.932044000624123)-- (12.065255705343276,4.022044087072661);
\draw (12.065255705343276,4.022044087072661)-- (11.62987777582714,2.813005683322278);
\draw (11.62987777582714,2.813005683322278)-- (11.111173349976594,3.932044000624123);
\draw (11.24,4.24)-- (11.339016592404569,3.440501407369559);
\draw (11.339016592404569,3.440501407369559)-- (11.368524191785584,2.808524191785584);
\draw (11.011420616478699,3.597348330213822)-- (11.339016592404569,3.440501407369559);
\draw (11.339016592404569,3.440501407369559)-- (11.841215286604793,2.8038571769057428);
\draw (11.841215286604793,2.8038571769057428)-- (11.842001869866536,3.402071317263379);
\draw (11.842001869866536,3.402071317263379)-- (11.94,4.22);
\draw (11.368524191785584,2.808524191785584)-- (11.842001869866536,3.402071317263379);
\draw (11.842001869866536,3.402071317263379)-- (12.152388287902083,3.654396713784991);
\draw (11.011420616478699,3.597348330213822)-- (11.58190445219039,3.9681628234264212);
\draw (11.58190445219039,3.9681628234264212)-- (11.94,4.22);
\draw (11.24,4.24)-- (11.58190445219039,3.9681628234264212);
\draw (11.58190445219039,3.9681628234264212)-- (12.152388287902083,3.654396713784991);
\draw (7.42639262953199,3.3075258095039475)-- (5.76639262953199,4.847525809503951);
\draw (7.42639262953199,3.3075258095039475)-- (8.646392629531984,4.847525809503951);
\draw (7.174916821317573,2.6560500012895316)-- (7.8549168213175715,1.5960500012895318);
\draw (7.087593317777533,1.655834367472979)-- (7.6476079161367805,2.6513829864096903);
\draw (6.461260848690881,4.628628098819535)-- (6.824975778742969,3.4205562038913726);
\draw (7.046392629531989,4.087525809503951)-- (6.01913587601432,4.0010958795366705);
\draw (8.146392629531988,4.587525809503953)-- (7.958780917434071,3.501922523288939);
\draw (7.746392629531987,4.06752580950395)-- (8.626392629531988,4.247525809503951);
\draw (6.917565979508584,3.7795698101280717)-- (7.871648334875265,3.86956989657661);
\draw (7.871648334875265,3.86956989657661)-- (7.436270405359131,2.6605314928262254);
\draw (7.436270405359131,2.6605314928262254)-- (6.917565979508584,3.7795698101280717);
\draw (7.046392629531989,4.087525809503951)-- (7.225151636660757,3.1159931775567395);
\draw (7.225151636660757,3.1159931775567395)-- (7.174916821317573,2.6560500012895316);
\draw (6.824975778742969,3.4205562038913726)-- (7.225151636660757,3.1159931775567395);
\draw (7.225151636660757,3.1159931775567395)-- (7.6476079161367805,2.6513829864096903);
\draw (7.6476079161367805,2.6513829864096903)-- (7.648394499398525,3.249597126767328);
\draw (7.648394499398525,3.249597126767328)-- (7.746392629531987,4.06752580950395);
\draw (7.174916821317573,2.6560500012895316)-- (7.648394499398525,3.249597126767328);
\draw (7.648394499398525,3.249597126767328)-- (7.958780917434071,3.501922523288939);
\draw (6.824975778742969,3.4205562038913726)-- (7.388297081722379,3.815688632930369);
\draw (7.388297081722379,3.815688632930369)-- (7.746392629531987,4.06752580950395);
\draw (7.046392629531989,4.087525809503951)-- (7.388297081722379,3.815688632930369);
\draw (7.388297081722379,3.815688632930369)-- (7.958780917434071,3.501922523288939);
\draw (7.42639262953199,3.3075258095039466)-- (7.452247984426925,1.6140000638857046);
\draw (8.50942467006298,0.4491969318618469)-- (5.2326786715536775,0.5670654929593038);
\draw (3.888977075042668,3.4312715276275045)-- (3.8675464275704035,1.9739874995134927);
\draw (3.3746415357083106,3.420556203891372)-- (4.4676045567938205,1.9418415283050954);
\draw (3.310349593291516,1.9954181469857577)-- (4.424743261849291,3.4098408801552393);
\draw (7.452247984426925,1.6140000638857046)-- (7.4571798791747845,1.0310390107338383);
\begin{scriptsize}
\draw [fill=qqqqff] (-2.8,0.74) circle (1.5pt);
\draw [fill=qqqqff] (-2.82,3.3) circle (1.5pt);
\draw [fill=qqqqff] (-3.62,0.76) circle (1.5pt);
\draw [fill=qqqqff] (-1.94,3.3) circle (1.5pt);
\draw [fill=qqqqff] (-3.66,3.32) circle (1.5pt);
\draw [fill=qqqqff] (-1.92,0.74) circle (1.5pt);
\draw [fill=uuuuuu] (-2.7849513365318708,2.0225140507196713) circle (1.5pt);
\draw [fill=qqqqff] (13.04,-1.06) circle (1.5pt);
\draw [fill=qqqqff] (11.368524191785584,2.808524191785584) circle (1.5pt);
\draw [fill=qqqqff] (12.048524191785583,1.7485241917855843) circle (1.5pt);
\draw [fill=qqqqff] (11.281200688245544,1.8083085579690314) circle (1.5pt);
\draw [fill=qqqqff] (11.841215286604793,2.8038571769057428) circle (1.5pt);
\draw [fill=qqqqff] (10.654868219158892,4.781102289315581) circle (1.5pt);
\draw [fill=qqqqff] (11.011420616478699,3.597348330213822) circle (1.5pt);
\draw [fill=qqqqff] (11.24,4.24) circle (1.5pt);
\draw [fill=qqqqff] (10.212743246482331,4.15357007003272) circle (1.5pt);
\draw [fill=qqqqff] (12.34,4.74) circle (1.5pt);
\draw [fill=qqqqff] (12.152388287902083,3.654396713784991) circle (1.5pt);
\draw [fill=qqqqff] (11.94,4.22) circle (1.5pt);
\draw [fill=qqqqff] (12.82,4.4) circle (1.5pt);
\draw [fill=qqqqff] (10.28,-0.08) circle (1.5pt);
\draw [fill=qqqqff] (11.38,-0.06) circle (1.5pt);
\draw [fill=qqqqff] (10.84,0.34) circle (1.5pt);
\draw [fill=qqqqff] (10.8,-0.66) circle (1.5pt);
\draw [fill=qqqqff] (12.0,0.04) circle (1.5pt);
\draw [fill=qqqqff] (12.76,-0.3) circle (1.5pt);
\draw [fill=qqqqff] (12.467304651784184,0.21753383039956012) circle (1.5pt);
\draw [fill=qqqqff] (12.26,-0.52) circle (1.5pt);
\draw [fill=ffqqqq] (10.810332922318127,-0.07035758323057957) circle (1.5pt);
\draw [fill=uuuuuu] (12.37077987184845,-0.12587520582693798) circle (1.5pt);
\draw [fill=xdxdff] (11.111173349976594,3.932044000624123) circle (1.5pt);
\draw [fill=xdxdff] (12.065255705343276,4.022044087072661) circle (1.5pt);
\draw [fill=xdxdff] (12.540254011302427,4.621632112627653) circle (1.5pt);
\draw [fill=xdxdff] (10.455618016773526,4.540209791667934) circle (1.5pt);
\draw [fill=xdxdff] (11.62987777582714,2.813005683322278) circle (1.5pt);
\draw [fill=xdxdff] (11.645855354894934,1.766474254381757) circle (1.5pt);
\draw [fill=xdxdff] (11.096123423991592,0.23584652570527764) circle (1.5pt);
\draw [fill=xdxdff] (12.166190115029703,0.14307158075620474) circle (1.5pt);
\draw [fill=xdxdff] (10.502623980897477,-0.400045734752703) circle (1.5pt);
\draw [fill=xdxdff] (12.665348568099173,-0.5441755647742238) circle (1.5pt);
\draw [fill=ffqqqq] (11.636319409068376,2.3910787060212253) circle (1.5pt);
\draw [fill=xdxdff] (11.339016592404569,3.440501407369559) circle (1.5pt);
\draw [fill=xdxdff] (11.842001869866536,3.402071317263379) circle (1.5pt);
\draw [fill=qqqqff] (11.58190445219039,3.9681628234264212) circle (1.5pt);
\draw [fill=ffqqqq] (10.832831267549446,4.190264968658947) circle (1.5pt);
\draw [fill=ffqqqq] (12.257962880762863,4.265297406864603) circle (1.5pt);
\draw [fill=xfqqff] (7.174916821317573,2.6560500012895316) circle (1.5pt);
\draw [fill=xfqqff] (7.8549168213175715,1.5960500012895318) circle (1.5pt);
\draw [fill=xfqqff] (7.087593317777533,1.655834367472979) circle (1.5pt);
\draw [fill=xfqqff] (7.6476079161367805,2.6513829864096903) circle (1.5pt);
\draw [fill=xfqqff] (6.461260848690881,4.628628098819535) circle (1.5pt);
\draw [fill=xfqqff] (6.824975778742969,3.4205562038913726) circle (1.5pt);
\draw [fill=xfqqff] (7.046392629531989,4.087525809503951) circle (1.5pt);
\draw [fill=xfqqff] (6.01913587601432,4.0010958795366705) circle (1.5pt);
\draw [fill=xfqqff] (8.146392629531988,4.587525809503953) circle (1.5pt);
\draw [fill=xfqqff] (7.958780917434071,3.501922523288939) circle (1.5pt);
\draw [fill=xfqqff] (7.746392629531987,4.06752580950395) circle (1.5pt);
\draw [fill=xfqqff] (8.626392629531988,4.247525809503951) circle (1.5pt);
\draw [fill=xfqqff] (6.917565979508584,3.7795698101280717) circle (1.5pt);
\draw [fill=xfqqff] (7.871648334875265,3.86956989657661) circle (1.5pt);
\draw [fill=xfqqff] (8.346646640834413,4.469157922131605) circle (1.5pt);
\draw [fill=xfqqff] (6.262010646305516,4.387735601171884) circle (1.5pt);
\draw [fill=xfqqff] (7.436270405359131,2.6605314928262254) circle (1.5pt);
\draw [fill=xfqqff] (7.452247984426925,1.6140000638857046) circle (1.5pt);
\draw [fill=ffqqqq] (7.442712038600367,2.238604515525178) circle (1.5pt);
\draw [fill=ffxfqq] (7.225151636660757,3.1159931775567395) circle (1.5pt);
\draw [fill=ffxfqq] (7.648394499398525,3.249597126767328) circle (1.5pt);
\draw [fill=ffxfqq] (7.388297081722379,3.815688632930369) circle (1.5pt);
\draw [fill=ffqqqq] (6.639113861608572,4.037892859264222) circle (1.5pt);
\draw [fill=ffqqqq] (8.064355510294844,4.11282321636854) circle (1.5pt);
\draw [fill=qqqqff] (8.50942467006298,0.4491969318618469) circle (1.5pt);
\draw [fill=qqqqff] (5.2326786715536775,0.5670654929593038) circle (1.5pt);
\draw [fill=qqqqff] (3.888977075042668,3.4312715276275045) circle (1.5pt);
\draw [fill=qqqqff] (3.8675464275704035,1.9739874995134927) circle (1.5pt);
\draw [fill=qqqqff] (3.3746415357083106,3.420556203891372) circle (1.5pt);
\draw [fill=qqqqff] (4.4676045567938205,1.9418415283050954) circle (1.5pt);
\draw [fill=qqqqff] (3.310349593291516,1.9954181469857577) circle (1.5pt);
\draw [fill=qqqqff] (4.424743261849291,3.4098408801552393) circle (1.5pt);
\draw [fill=ffqqqq] (3.8787888834326565,2.738474498146671) circle (1.5pt);
\end{scriptsize}
\end{tikzpicture}
\end{center}

\begin{center}
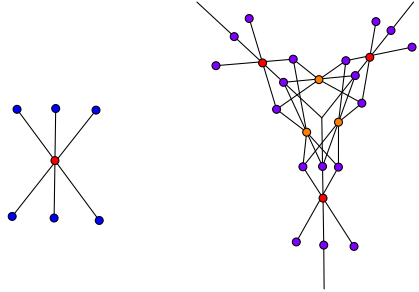
 Figure 18. Construction of a helical configuration $n_{3}$ from the skeleton of any solid with 3-valent vertices.
 \end{center}

\end{example}

\section{Conclusion and open problems}

In this paper we proved the existence of Platonic configurations of points and lines for most of the cases of Problem 1. We did this by constructing concrete examples rather than trying to develop a theory of Platonic configurations. In particular:

\begin{itemize}
\item for all the Platonic solids $P$ we were able to find 1-layer Platonic configurations $(n_{3})$ (see Problems $3,4,7$ and Examples $2, 4$) 
\item the examples of $(n_{4})$ cases we constructed as 2-layer or 4-layer configurations (Problems 8,9)
\item the examples of $(n_{5})$ configurations may be constructed using the axes as "meting places" of 1-valent points (as in Example 5, where we used this trick for construction of a tetrahedral $(n_{3})$ configuration).
\end{itemize}

The key ingredients of the solution (of Problem 1) are: 

1) systematic progress from $k = 3$ to $k = 4 $ to $k = 5$; 

2) observing that the solutions of some cases (subproblems of Problem 1) may be obtained by the same procedure; 

3) observing how the differences in the parameters $d$ and $m$ of Platonic solids reflect in variations of the general procedure (given in Algorithm 1);

4) using some known tricks and methods, e.g. from the Gr\"unbaum's "incidence calculus" (e.g. "parallel shift", "radial projection" etc.);

5) using the vertices, edges, face centers and axes of solids as "meeting places" where we link together two or more copies of the same configuration;

6)  using triangular grid to construct examples of (symmetrical) planar configurations of the desired class $(p_{q},n_{k})$, especially with 3-valent lines (using either the greed on computer screen or the paper grid and coins turns the search for configurations into an interesting game!);

7) using GeoGebra to draw examples of configurations.

Our methods allow to produce many more Platonic configurations $(p_{q},n_{k})$. 
%However, we haven't given all the possible examples of \textit{connected} Platonic $(n_{5})$ configurations for all the five Platonic solids. 

Our main construction of highly symmetrical configurations, based on the idea of placing identical copies of some configuration along each of the faces of the chosen Platonic solid, may be viewed as a natural generalization of the constructions of planar configurations with a cyclical or dihedral group of symmetries, starting from a regular polygon, explored by several authors (as explained in \cite{Grunbaum}). Likewise, the idea to use smaller planar configurations to construct other planar configurations may be found e.g. in the papers \cite{BermBoko, BBGT}.  Constructions of configurations from identical copies of smaller configurations may be found also in the incidence calculus of Gr\"unbaum \cite{Grunbaum, PisaServ}.  

We just mention that some Platonic configurations may be formed also directly from the vertices and edges of various uniform polyhedra. For example, the vertices and edges of the uniform polyhedron, called Great cubicuboctahedron, form a $(n_{4})$ configuration. Platonic configurations may be obtained also from compounds (compound bodies) and from stellations of some polyhedra.  %(DODAJ REFERENCO IZ WIKIPEDIJE). 
This theme will be explored in another paper. 

Future work on Platonic configurations may be focused also on:

\begin{itemize}
\item finding their Levi graphs
\item finding other examples of spider-web and helical configurations
\item finding examples of Platonic configurations with some points on (the intersection of) their diagonals
\item finding minimal (with minimal number of points) configurations of a given class $(P(p_{q},n_{k})$.
 
\end{itemize}

%We also just mention that Platonic $(n_{5})$ configurations may be obtained from 1-layer Platonic configurations $(p_{4},n_{5})$.

The purpose of the rest of the paper is to show that the concept of Platonic configurations is not limited to configurations of points of lines. It is not difficult to obtain examples of connected Platonic configurations whose blocks are spheres or circles or ellipses (instead of lines). 

\begin{example}\label{examplecube}
Take a cube. Each of its 8 vertices $v_{i}$ is the top of the pyramid whose base is an equilateral triangle whose vertices are adjacent to $v_{i}$ in the cube. Thus we get a configuration $(8_{4})$ of 8 points and 8 pyramids whose symmetry group is isomorphic to the full symmetry group of the cube. Each pyramid may be, of course, represented also by a circumscribed sphere.
\end{example}

An identical procedure of ``truncating'' pyramids at each vertex and taking them as ``blocks'' of configurations may be applied to other solids:
\begin{prop}\label{prviizrek} For each of the five Platonic solids and each of the 75 uniform polyhedra with parameters $P(v,e,f,m,n)$ there is a polyhedral configuration $(v_{d})$ of points and spheres, whose symmetry group $Sym(C)$ is isomorphic to the symmetry group of $P$.
\end{prop}

Proof. Take a Platonic solid or a uniform polyhedron $P = P(v,e,f,m,n)$. Choose any of its vertices $v_{i}$ as the top of a pyramid whose other $m$ base vertices are adjacent to $v_{i}$ in $P$. To each pyramid corresponds a sphere with the center in $v_{i}$, passing through the vertices adjacent to $v_{i}$. The symmetry group $Sym(C)$ of the obtained configuration of $v$ points and $v$ spheres is isomorphic to the full symmetry group of $P$, by the construction. \qed

\begin{example}\label{exampleicosahedron}
Take the icosahedron $I = P(12,30,20,5,3)$ and truncate from it two antipodal pyramids to get a pentagonal antiprism incident with 10 vertices of $P$. This antiprism may be circumscribed a sphere; take this sphere as one of the $\frac{v}{2} = 6$ identical blocks of a configuration $C$. Since each of the 6 blocks is incident with 10 vertices -- points of the configuration $C$ -- each of the 12 vertices must be (by symmetry) incident with $\frac{6 \times 10}{12}$ blocks and  we get a configuration $(12_{5}, 6_{10})$ of 12 points and 6 spheres.
\end{example}

Acknowledgement. The authors appreciate the possibility to apply in their work on configurations the software GeoGebra http://www.geogebra.org.

\end{document}